%% file: essay.tex
\documentclass{report}

\input{diagrams}
\diagramstyle[PostScript=dvips,nohug]
\usepackage{makeidx}
\usepackage{epsfig}
\input{essab}

\input{essfor}

\input{essthms}

\input{esshdpics}

\makeindex

\begin{document}
\include{esstitle}

\include{preface}

\tableofcontents
\include{essintro}

\include{essprelims}
\include{essom}
\include{essglob}

\include{essgray}

\include{essopap}

\include{essrefs}
\input{esssee}
\printindex
\end{document}

%% file: essab.tex
\usepackage{amssymb}
\usepackage{latexsym}

\newcommand{\mcm}[3]{\newcommand{#1}[#2]{{\ensuremath{#3}}}}

\mcm{\mc}{1}{\mathcal{#1}}
\mcm{\mr}{1}{\mathrm{#1}}
\mcm{\mi}{1}{\mathit{#1}}
\mcm{\mb}{1}{\mathbf{#1}}
\mcm{\cat}{1}{\mc{#1}}
\mcm{\scat}{1}{\mathbb{#1}}
\mcm{\fcat}{1}{\mb{#1}}  
\mcm{\twid}{1}{\widetilde{#1}}
\newcommand{\url}[1]{\mbox{\tt #1}}

\mcm{\pr}{2}{(#1,#2)}
\mcm{\triple}{3}{(#1,#2,#3)}
\mcm{\range}{2}{#1,\,\ldots\,,#2}
\mcm{\tuplebts}{1}{( #1 )}
\mcm{\tuple}{3}{\tuplebts{\range{#1,#2}{#3}}}
\mcm{\bftuple}{2}{\tuplebts{\range{#1}{#2}}}
\mcm{\rttuple}{1}{\tuplebts{\,\ldots\,,#1}}
\mcm{\atuplebts}{1}{\langle #1 \rangle}
\mcm{\atuple}{3}{\atuplebts{\range{#1,#2}{#3}}}
\mcm{\abftuple}{2}{\atuplebts{\range{#1}{#2}}}
\mcm{\arttuple}{1}{\atuplebts{\,\ldots\,,#1}}

\mcm{\elt}{0}{\in}
\mcm{\such}{0}{\:|\:}
\mcm{\nat}{0}{\mathbb{N}}	
\mcm{\without}{0}{\setminus}
\mcm{\sub}{0}{\subseteq}

\mcm{\iso}{0}{\cong}
\mcm{\eqv}{0}{\simeq}
\mcm{\id}{0}{\mi{id}}
\mcm{\emptybk}{0}{\:\:}
\mcm{\blank}{0}{(\:\:)}
\mcm{\op}{0}{\mr{op}}
\mcm{\ftrcat}{2}{[#1,#2]}
\mcm{\Cat}{0}{\fcat{Cat}}
\mcm{\Sets}{0}{\fcat{Set}}
\mcm{\Hom}{0}{\mr{Hom}}
\mcm{\End}{0}{\mr{End}}
\mcm{\homset}{3}{#1(#2,#3)}
\mcm{\mtihom}{4}{#1(\range{#2}{#3};#4)} 
\mcm{\mtiendo}{2}{#1(#2;#2)}
\mcm{\dom}{0}{\mr{dom}}
\mcm{\cod}{0}{\mr{cod}}
\mcm{\Mon}{0}{\mb{Mon}}

\newcommand{\done}{\hfill\ensuremath{\Box}}

\mcm{\go}{0}{\rTo}
\mcm{\og}{0}{\lTo}
\mcm{\goby}{1}{\rTo^{#1}}
\mcm{\ogby}{1}{\lTo^{#1}}
\mcm{\goesto}{0}{\,\longmapsto\,}
\mcm{\slob}{3}{(#1 \goby{#2} #3)}
\mcm{\bktdslob}{3}{\slob{#1}{#2}{#3}}
\mcm{\vslob}{3}
	{\left.
	\begin{diagram}[height=1.5em]
	#1		\\
	\dTo>{\,#2}	\\
	#3		\\
	\end{diagram}
	\right.}
\mcm{\bktdvslob}{3}{\left( \vslob{#1}{#2}{#3} \right)}
\newarrow{Get}....>
\newarrow{Goesto}|---{->}

\newenvironment{slopeydiag}
	{\begin{diagram}[size=2em]}
	{\end{diagram}}

\newenvironment{ntdiag}
	{\begin{diagram}[size=1.5em,noPS]}
	{\end{diagram}}

\newcommand{\pullshape}
	{\setlength{\unitlength}{1em}
	\begin{picture}(2,5)(-1,-5)
	\put(0,-5){\line(1,1){1}}
	\put(0,-5){\line(-1,1){1}}
	\end{picture}}
\newcommand{\Spbk}{\overprint{\raisebox{-2.5em}{\pullshape}}}
\mcm{\diso}{0}{\sim}     
\mcm{\vdiso}{0}{\wr}

\mcm{\node}{0}{\bullet}
\mcm{\enode}{0}{\circ}
\newcommand{\nodelt}[1]{\raisebox{1ex}{\ensuremath{#1}}}
\newenvironment{tree}
	{\begin{diagram}[height=1em,width=.75em,abut,noPS,tight]}	
	{\end{diagram}}
\newcommand{\lt}[1]{\ldLine(#1,2)}
\newcommand{\rt}[1]{\rdLine(#1,2)}
\newcommand{\dn}{\dLine}

\newenvironment{opetope}
	{\begin{diagram}[size=1em,abut,tight,noPS]}
	{\end{diagram}}
\newarrow{Edge}---->
\newcommand{\cnr}{}

\mcm{\Sym}{0}{\mr{Sym}}
\mcm{\bee}{0}{\cat{B}}
\mcm{\trip}{2}{\homset{\bee}{#1}{#2}}
\mcm{\beep}{0}{\cat{B'}}
\mcm{\tripp}{2}{\homset{\beep}{#1}{#2}}
\mcm{\goiso}{0}{\goby{\diso}}
\mcm{\gobyciso}{3}{{#1}: {#2} \goiso {#3}}
\newarrow{Equals}=====
\mcm{\rep}{1}{\trip{\dashbk}{#1}}
\mcm{\ftrbi}{2}{\ftrcat{#1}{#2}}
\mcm{\Pshf}{0}{\ftrbi{\bee^{\op}}{\Cat}}

\mcm{\ess}{0}{\cat{S}}
\mcm{\blob}{0}{\scriptscriptstyle{\bullet}}
\mcm{\ust}{0}{{}^{*}}
\mcm{\ubl}{0}{{}^{\blob}}
\mcm{\stbk}{0}{\blank^{*}}
\mcm{\blbk}{0}{\blank^{\blob}}
\mcm{\Mnd}{0}{\triple{\stbk}{\eta}{\mu}}
\mcm{\Cartpr}{0}{\pr{\cat{S}}{\ust}}

\mcm{\spn}{3}{#2 \og #1 \go #3}
\mcm{\spaan}{5}{#2 \ogby{#4} #1 \goby{#5} #3}
\mcm{\gph}{2}{\spn{#1}{#2^{*}}{#2}}
\mcm{\graph}{4}{\spaan{#1}{#2^{*}}{#2}{#3}{#4}}
\mcm{\comp}{0}{\mi{comp}}
\mcm{\ids}{0}{\mi{ids}}

\mcm{\Slice}{0}{\cat{S}/C_0}
\mcm{\unit}{0}{\mi{unit}}
\mcm{\mult}{0}{\mi{mult}}
\mcm{\Imnd}{0}{\triple{\blbk}{\unit}{\mult}}
\mcm{\Icartpr}{0}{\pr{\Slice}{\ubl}}
\mcm{\Alg}{0}{\mb{Alg}}

\mcm{\Multicat}{0}{\fcat{Multicat}}
\mcm{\Graph}{0}{\fcat{Graph}}
\mcm{\ush}{0}{{}^{\sharp}}

\mcm{\nl}{1}{\stackrel{\textstyle #1}{\node}}
\newlength{\nllwidth}
\newlength{\nllheight}
\newcommand{\stackbelow}[2]{%
\settowidth{\nllwidth}{\ensuremath{#1}\ensuremath{#2}}%
\settoheight{\nllheight}{\ensuremath{#2}}%
\addtolength{\nllheight}{.3ex}%
\mbox{%
\ensuremath{#1}%
\hspace{-.5\nllwidth}%
\raisebox{-1\nllheight}{\ensuremath{#2}}}}
\mcm{\nll}{1}{\stackbelow{\node}{#1}}
\mcm{\nlal}{2}{\stackbelow{\nl{#1}}{#2}}
\mcm{\Set}{0}{\fcat{Set}}
\mcm{\bdry}{0}{\partial}
\mcm{\Tr}{0}{\fcat{Tr}}
\mcm{\wun}{0}{\fcat{1}}
\mcm{\Cpn}{1}{\pr{\Set/S_{#1}}{T_{#1}}}
\mcm{\Operadof}{1}{#1\hyph\fcat{Operad}}
\mcm{\Strucof}{1}{#1\hyph\fcat{Struc}}
\mcm{\Catof}{1}{#1\hyph\Cat}
\mcm{\CatOF}{1}{\Catof{\mb{#1}}}
\mcm{\BicatOF}{1}{\mb{#1}\hyph\Bicat}
\mcm{\PD}{1}{\fcat{PD}_{#1}}
\mcm{\pd}{1}{\fcat{pd}_{#1}}
\mcm{\TR}{1}{\fcat{TR}(#1)}
\mcm{\Gy}{0}{\fcat{Gy}}
\mcm{\Ssq}{0}{\fcat{Ssq}}
\mcm{\tooc}{0}{\fcat{2\hyph\Cat}}
\mcm{\Bicat}{0}{\fcat{Bicat}}
\mcm{\tee}{0}{\cat{T}}
\mcm{\teep}{0}{\cat{T'}}
\mcm{\Cay}{0}{\mb{Cay}}
\mcm{\stee}{0}{\scat{T}}
\mcm{\steep}{0}{\scat{T'}}
\mcm{\parpair}{2}{\stackrel{\rTo^{#1}}{\rTo_{#2}}}
\mcm{\thomset}{2}{\homset{\stee}{#1}{#2}}
\mcm{\cayset}{1}{\coprod\thomset{B}{#1}}
\mcm{\dcayset}{1}{\coprod_{B\elt\stee_0}\thomset{B}{#1}}
\mcm{\bkpr}{1}{{}^{(#1)}}

\newcommand{\ucontents}[2]{\addcontentsline{toc}{#1}{\numberline{}{#2}}}
\mcm{\oppair}{2}{\stackrel{\rTo^{#1}}{\lTo_{#2}}}
\mcm{\hyph}{0}{\mbox{-}}
\mcm{\ob}{0}{\mr{ob}\,}
\mcm{\dashbk}{0}{\mbox{---}}
\newcommand{\from}{:}
\mcm{\gobyc}{3}{{#1}\from {#2} \go {#3}}	
\mcm{\of}{0}{\raisebox{0.2mm}{\ensuremath{\scriptstyle\circ}}}
\mcm{\diagspace}{0}{\mbox{\hspace{2em}}}
\mcm{\Multicatof}{1}{#1\hyph\Multicat}
\mcm{\Graphof}{1}{#1\hyph\Graph}
\usepackage{amssymb}

%% file: essfor.tex
\newcommand{\EmptyOne}
{\begin{tree}
\enode	\\
\dn	\\
\node	\\
\end{tree}}

\newcommand{\Oak}[5]
{\begin{tree}
 & & & &\enode& & & & & & & & \\
 & & & &\dn& & & & & & & & \\
\nl{#1}& &\nl{#2}& &\node& & & &\nl{#4}& & & &\nl{#5}\\
 &\rt{2}&\dn&\lt{2}& & & & & &\rt{2}& &\lt{2}& \\
 & &\node& & & &\nl{#3}& & & &\node& & \\
 & & &\rt{4}& & &\dn& & &\lt{4}& & & \\
 & & & & & &\node& & & & & & \\
\end{tree}}

\newcommand{\Pear}[3]
{\begin{tree}
	&	&	&	&\enode	&&	&	&	\\
	&	&	&	&\dn	&&	&	&	\\
\nl{#1}&	&	&	&\node	&&	&	&	\\
	&\rt{2}	&	&\lt{2}	&	&&	&	&	\\
	&	&\node	&	&	&&	&	&\nl{#2}\\
	&	&	&\rt{4}	&	&&	&\lt{2}	&	\\
	&	&	&	&	&&\nll{#3}&	&	\\
\end{tree}}

\newcommand{\Orange}[3]
{\begin{tree}
	&	&	&	&\nl{#2}\\
	&	&	&	&\dn	\\
\nl{#1}&	&	&	&\node	\\
	&\rt{2}	&	&\lt{2}	&	\\
	&	&\nll{#3}&	&	\\
\end{tree}}

\newcommand{\Apple}[3]
{\begin{tree}
	&	&\enode	&	&	\\
	&	&\dn	&	&	\\
\nl{#1}	&	&\node	&	&\nl{#2}\\
	&\rt{2}	&\dn	&\lt{2}	&	\\
	&	&\node	&	&	\\
	&	&\dn	&	&	\\
	&	&\nll{#3}&	&	\\
\end{tree}}

\newcommand{\MixedFruit}[5]
{\begin{tree}
 & & & &\nl{#2}& & & & & & &\enode& & \\
 & & & &\dn    & & & & & & &\dn   & & \\
\nl{#1}& & & &\node& & &\enode& &\nl{#3}& &\node& &\nl{#4}\\
 &\rt{2}& &\lt{2}& & & &\dn& & &\rt{2}&\dn&\lt{2}& \\
 & &\node& & & & &\node& & & &\node& & \\
 & & &\rt{3}& & &\lt{2}& & & & &\dn& & \\
 & & & & &\node& & & & & &\node& & \\
 & & & & & &\rt{3}& & & &\lt{3}& & & \\
 & & & & & & & &\nll{#5}& & & & & \\
\end{tree}}

\newcommand{\treedc}{
\begin{tree}
\node &\node            &\node &      &      &\node \\
      &\rt{1} \dn \lt{1}&      &      &      & \dn  \\
      &\node            &      &\node &      &\node \\
      &                 &\rt{2}&\dn   &\lt{2}&      \\
      &                 &      &\node &      &      \\
\end{tree}}

\newcommand{\treedcbang}{
\begin{tree}
\node &                 &\node &      &      &\node \\
      &\rt{1}     \lt{1}&      &      &      & \dn  \\
      &\node            &      &\node &      &\node \\
      &                 &\rt{2}&\dn   &\lt{2}&      \\
      &                 &      &\node &      &      \\
\end{tree}}

\newcommand{\treec}{
\begin{tree}
\node	&	&\node	&	&\node	\\
	&\rt{2}	&\dn	&\lt{2}	&	\\
	&	&\node	&	&	\\
\end{tree}}

\newcommand{\treeccbang}{
\begin{tree}
	&	&	&\node	&			&\node	\\
	&	&	&	&\rt{1}  \lt{1}		&	\\
\node	&	&\node	&	&\node			&	\\
	&\rt{2}	&\dn	&\lt{2}	&			&	\\
	&	&\node	&	&			&	\\
\end{tree}}

\newcommand{\treegcbang}{
\begin{tree}
\node & 	&\node 	&	&\node 	&\node 	&\node 		& 	\\
&\rt{1}  \lt{1} &	&	&	&\rt{1}\dn\lt{1}&	&	\\
	&\node	&	&\node	&	&\node	&		&	\\
	&	&\rt{2}	&\dn	&\lt{2}	&	&		&	\\
	&	&	&\node	&	&	&		&	\\
\end{tree}}

\newcommand{\treebaa}{
\begin{tree}
\node	&		&\node	\\
	&\rt{1}\lt{1}	&	\\
	&\node		&	\\
	&\dn		&	\\
	&\node		&	\\
	&\dn		&	\\
	&\node		&	\\
\end{tree}}

\newcommand{\treeba}{
\begin{tree}
\node	&		&\node	\\
	&\rt{1}\lt{1}	&	\\
	&\node		&	\\
	&\dn		&	\\
	&\node		&	\\
\end{tree}}

\newcommand{\treebcb}{
\begin{tree}
\node	&	&\node	&	&	&	\\
	&\rt{1}\lt{1}&	&	&	&	\\
	&\node	&	&\node	&	&\node	\\
	&\dn	&	&	&\rt{1}\lt{1}&	\\
	&\node	&	&	&\node	&	\\
	&	&\rt{1}	&\lt{2}	&	&	\\
	&	&\node	&	&	&	\\
\end{tree}}

\newcommand{\treeb}{
\begin{tree}
\node	&		&\node	\\
	&\rt{1}\lt{1}	&	\\
	&\node		&	\\
\end{tree}}

\newcommand{\treedb}{
\begin{tree}
\node	&	&	&\node	&\node		&\node	\\
\dn	&	&	&	&\rt{1}\dn\lt{1}&	\\
\node	&	&	&	&\node		&	\\
	&\rt{2}	&	&\lt{2}	&		&	\\
	&	&\node	&	&		&	\\
\end{tree}}

\newcommand{\treebb}{
\begin{tree}
	&	&	&\node	&		&\node	\\
	&	&	&	&\rt{1}\lt{1}	&	\\
\node	&	&	&	&\node		&	\\
	&\rt{2}	&	&\lt{2}	&		&	\\
	&	&\node	&	&		&	\\
\end{tree}}

\newcommand{\treefg}{
\begin{tree}
\node	&\node	&\node	&\node	&	&\node	&	&\node	&
	&	&	&	&	\\
\dn	&	&\rt{1}\dn\lt{1}&&	&	&\rt{1}\lt{1}&	&
	&	&	&	&	\\
\node	&	&\node	&	&\node	&	&\node	&	&
\node	&	&\node	&	&\node	\\
	&\rt{6}	&	&\rt{4}	&	&\rt{2}	&\dn	&\lt{2}	&
	&\lt{4}	&	&\lt{6}	&	\\
	&	&	&	&	&	&\node	&	&
	&	&	&	&	\\
\end{tree}}

\newcommand{\treeec}{
\begin{tree}
	&	&\node	&\node		&\node	&\node	&	&\node	\\
	&	&	&\rt{1} \dn \lt{1}&	&	&\rt{1}\lt{1}&	\\
\node	&	&	&\node		&	&	&\node	&	\\
	&\rt{3}	&	&\dn		&	&\lt{3}	&	&	\\
	&	&	&\node		&	&	&	&	\\
\end{tree}}

\newcommand{\treeaaq}{
\begin{tree}
\node	\\
\dn	\\
\node	\\
\dn	\\
\node	\\
\vdots	\\
\node	\\
\dn	\\
\node	\\
\end{tree}}

\newcommand{\treeq}{
\begin{tree}
\node	&\node	&	&\ldots	&	&\node	\\
	&\rt{2}	&\rt{1}	&	&\lt{3}	&	\\
	&	&\node	&	&	&	\\
\end{tree}}

\newcommand{\treea}{
\begin{tree}
\node 	\\
\dn	\\
\node	\\
\end{tree}}

\newlength{\volt}
\newcommand{\transistor}[5]
{\setlength{\unitlength}{1\volt}
\begin{picture}(18,12)(-5,-6)
\put(0,6){\line(0,-1){12}}
\put(0,-6){\line(3,2){9}}
\put(0,6){\line(3,-2){9}}
\put(9,0){\line(1,0){2}}
\put(-2,4){\line(1,0){2}}
\put(-2,2){\line(1,0){2}}	
\put(-2,-4){\line(1,0){2}}
\put(12,-0.5){\ensuremath{#5}}
\put(-5,3.5){\ensuremath{#2}}
\put(-5,1.5){\ensuremath{#3}}		
\put(-5,-4.5){\ensuremath{#4}}
\thicklines
\put(-1.5,0){\line(1,0){.1}}	
\put(-1.5,-1){\line(1,0){.1}}		
\put(-1.5,-2){\line(1,0){.1}}		
\thinlines
\put(2,-0.5){\ensuremath{#1}}
\end{picture}}
\newcommand{\ctransistor}[5]
	{\raisebox{-6\volt}{\transistor{#1}{#2}{#3}{#4}{#5}}}

\newcommand{\bftransistor}[4]
{\setlength{\unitlength}{1\volt}
\begin{picture}(18,12)(-5,-6)
\put(0,6){\line(0,-1){12}}
\put(0,-6){\line(3,2){9}}
\put(0,6){\line(3,-2){9}}
\put(9,0){\line(1,0){2}}
\put(-2,4){\line(1,0){2}}
\put(-2,-4){\line(1,0){2}}
\put(12,-0.5){\ensuremath{#4}}
\put(-5,3.5){\ensuremath{#2}}
\put(-5,-4.5){\ensuremath{#3}}
\thicklines
\put(-1.5,1){\line(1,0){.1}}		
\put(-1.5,0){\line(1,0){.1}}		
\put(-1.5,-1){\line(1,0){.1}}
\thinlines		
\put(2,-0.5){\ensuremath{#1}}
\end{picture}}
\newcommand{\comptrans}[4]
{\setlength{\volt}{.5ex}
\setlength{\unitlength}{1\volt}
\begin{picture}(80,72)(0,-36)
\put(0,23){\bftransistor{#1}{}{}{}}
\put(0,6){\bftransistor{#2}{}{}{}}
\put(0,-35){\bftransistor{#3}{}{}{}}
\put(63,-6){\transistor{#4}{}{}{}}
\put(17,29){\line(2,-1){50}}
\put(17,12){\line(5,-1){50}}
\put(17,-29){\line(2,1){50}}
\thicklines
\put(24,-5){\line(1,0){.4}}		
\put(24,-9){\line(1,0){.4}}		
\put(24,-13){\line(1,0){.4}}
\thinlines		
\setlength{\volt}{1ex}
\end{picture}}

\newcommand{\threetransistor}[5]
{\setlength{\unitlength}{1\volt}
\begin{picture}(18,12)(-5,-6)
\put(0,6){\line(0,-1){12}}
\put(0,-6){\line(3,2){9}}
\put(0,6){\line(3,-2){9}}
\put(9,0){\line(1,0){2}}
\put(-2,4){\line(1,0){2}}
\put(-2,0){\line(1,0){2}}	
\put(-2,-4){\line(1,0){2}}
\put(12,-0.5){\ensuremath{#5}}
\put(-5,3.5){\ensuremath{#2}}
\put(-5,-0.5){\ensuremath{#3}}		
\put(-5,-4.5){\ensuremath{#4}}
\put(2,-0.5){\ensuremath{#1}}
\end{picture}}
\newcommand{\twotransistor}[4]
{\setlength{\unitlength}{1\volt}
\begin{picture}(18,12)(-5,-6)
\put(0,6){\line(0,-1){12}}
\put(0,-6){\line(3,2){9}}
\put(0,6){\line(3,-2){9}}
\put(9,0){\line(1,0){2}}
\put(-2,4){\line(1,0){2}}
\put(-2,-4){\line(1,0){2}}
\put(12,-0.5){\ensuremath{#4}}
\put(-5,3.5){\ensuremath{#2}}
\put(-5,-4.5){\ensuremath{#3}}
\put(2,-0.5){\ensuremath{#1}}
\end{picture}}
\newcommand{\notransistor}[2]
{\setlength{\unitlength}{1\volt}
\begin{picture}(18,12)(-5,-6)
\put(0,6){\line(0,-1){12}}
\put(0,-6){\line(3,2){9}}
\put(0,6){\line(3,-2){9}}
\put(9,0){\line(1,0){2}}
\put(12,-0.5){\ensuremath{#2}}
\put(2,-0.5){\ensuremath{#1}}
\end{picture}}
\newcommand{\freemoncatpic}
{\setlength{\volt}{.6ex}
\setlength{\unitlength}{1\volt}
\begin{picture}(18,48)
\put(0,3){\twotransistor{a_3}{s_4}{s_5}{s'_3}}
\put(0,18){\notransistor{a_2}{s'_2}}
\put(0,33){\threetransistor{a_1}{s_1}{s_2}{s_3}{s'_1}}
\put(5.7,0){\framebox(10.3,48){}}
\end{picture}}

%% file: essthms.tex
\newtheorem{thm}{Theorem}[section]
\newtheorem{defn}[thm]{Definition}

\newtheorem{constn}[thm]{Construction}

\newtheorem{lotsofremarks}[thm]{Remarks}
\newenvironment{remarks}[1]
	{\begin{lotsofremarks}\label{#1}\end{lotsofremarks}\begin{enumerate}}
	{\end{enumerate}}

\newtheorem{lotsofegs}[thm]{Examples}
\newenvironment{eg}[1]
	{\begin{lotsofegs}\label{#1}\end{lotsofegs}\nopagebreak\begin{enumerate}}
	{\end{enumerate}}

\newtheorem{conjecture}[thm]{Conjecture}


%% file: esshdpics.tex
\newlength{\gwidth}	
\newlength{\gvert}	
\newlength{\gdrop}	
\newlength{\gbaredrop}	
\newlength{\goffset}	
\newlength{\gtemp}	

\newcommand{\present}[1]{%
\makebox[1\gwidth]{%
\rule[-1\gdrop]{0ex}{1\gvert}%
\raisebox{-1\gbaredrop}{#1}}}

\newcommand{\presentl}[1]{%
\makebox[1\gwidth][l]{%
\rule[-1\gdrop]{0ex}{1\gvert}%
\raisebox{-1\gbaredrop}{#1}}}

\newcommand{\presentr}[1]{%
\makebox[1\gwidth][r]{%
\rule[-1\gdrop]{0ex}{1\gvert}%
\raisebox{-1\gbaredrop}{#1}}}

\newcommand{\ginitdims}[2]{
\setlength{\unitlength}{1em}
\setlength{\goffset}{.25\unitlength}
\setlength{\gwidth}{#1\unitlength}
\setlength{\gvert}{#2\unitlength}
\setlength{\gdrop}{.5\gvert}
\addtolength{\gdrop}{-1\goffset}
\setlength{\gbaredrop}{1\gdrop}
\addtolength{\gvert}{.6\unitlength}
\addtolength{\gdrop}{.3\unitlength}}	

\newcommand{\cinitdims}[2]{
\setlength{\unitlength}{1em}
\setlength{\goffset}{.35\unitlength}
\setlength{\gwidth}{#1\unitlength}
\setlength{\gvert}{#2\unitlength}
\setlength{\gdrop}{.5\gvert}
\addtolength{\gdrop}{-1\goffset}
\setlength{\gbaredrop}{1\gdrop}
\addtolength{\gvert}{.6\unitlength}
\addtolength{\gdrop}{.3\unitlength}}	

\newcommand{\gsinitdims}[2]{
\setlength{\unitlength}{0.5em}
\setlength{\goffset}{.25\unitlength}
\setlength{\gwidth}{#1\unitlength}
\setlength{\gvert}{#2\unitlength}
\setlength{\gdrop}{.5\gvert}
\addtolength{\gdrop}{-1\goffset}
\setlength{\gbaredrop}{1\gdrop}
\addtolength{\gvert}{.6\unitlength}
\addtolength{\gdrop}{.3\unitlength}}	

\newcommand{\sidespic}[1]{%
\settowidth{\gtemp}{\ensuremath{#1}}%
\addtolength{\gwidth}{1\gtemp}}

\newcommand{\abovepic}[1]{%
\settoheight{\gtemp}{\ensuremath{#1}}%
\addtolength{\gvert}{1\gtemp}%
\settodepth{\gtemp}{\ensuremath{#1}}%
\addtolength{\gvert}{1\gtemp}}

\newcommand{\belowpic}[1]{%
\settoheight{\gtemp}{\ensuremath{#1}}%
\addtolength{\gvert}{1\gtemp}%
\addtolength{\gdrop}{1\gtemp}%
\settodepth{\gtemp}{\ensuremath{#1}}%
\addtolength{\gvert}{1\gtemp}%
\addtolength{\gdrop}{1\gtemp}}

\newcommand{\cell}[4]{\put(#1,#2){\makebox(0,0)[#3]{\ensuremath{#4}}}}
\mcm{\zmark}{0}{\scriptstyle{\bullet}}

\newcommand{\pregfst}[1]{%
\begin{picture}(0.5,0.2)(-0.5,-0.2)%
\cell{-0.1}{-0.2}{tr}{#1}%
\cell{0}{0}{c}{\zmark}%
\end{picture}}

\mcm{\gfst}{1}{%
\ginitdims{0.5}{0.4}%
\sidespic{#1}%
\belowpic{#1}%
\presentr{\pregfst{#1}}}

\newcommand{\preglst}[1]{%
\begin{picture}(0.5,0.2)(0,-0.2)%
\cell{0.1}{-0.2}{tl}{#1}%
\cell{0.05}{0}{c}{\zmark}%
\end{picture}}

\mcm{\glst}{1}{%
\ginitdims{.5}{.4}%
\sidespic{#1}%
\belowpic{#1}%
\presentl{\preglst{#1}}}

\newcommand{\preglft}[1]{%
\begin{picture}(0,0.2)(0,-0.2)%
\cell{-0.1}{-0.2}{tr}{#1}%
\cell{0.05}{0}{c}{\zmark}%
\end{picture}}

\mcm{\glft}{1}{%
\ginitdims{0}{.4}%
\belowpic{#1}%
\present{\preglft{#1}}}

\newcommand{\pregrgt}[1]{%
\begin{picture}(0,0.2)(0,-0.2)%
\cell{0.1}{-0.2}{tl}{#1}%
\cell{0.05}{0}{c}{\zmark}%
\end{picture}}

\mcm{\grgt}{1}{%
\ginitdims{0}{.4}%
\belowpic{#1}%
\present{\pregrgt{#1}}}

\newcommand{\pregblw}[1]{%
\begin{picture}(0,0.3)(0,-0.3)
\cell{0}{-0.3}{t}{#1}%
\cell{0.05}{0}{c}{\zmark}%
\end{picture}}

\mcm{\gblw}{1}{%
\ginitdims{0}{.6}%
\belowpic{#1}%
\present{\pregblw{#1}}}

\newcommand{\pregfbw}[1]{%
\begin{picture}(0,0.65)(0,-0.65)
\cell{0}{-0.65}{t}{#1}%
\cell{0.05}{0}{c}{\zmark}%
\end{picture}}

\mcm{\gfbw}{1}{%
\ginitdims{0}{1.3}%
\belowpic{#1}%
\present{\pregfbw{#1}}}

\newcommand{\pregzero}[1]{%
\begin{picture}(0.8,0.4)(-0.4,-0.4)
\cell{0}{-0.4}{t}{#1}%
\cell{0}{0}{c}{\zmark}%
\end{picture}}

\mcm{\gzero}{1}{%
\ginitdims{0.8}{.6}%
\belowpic{#1}%
\sidespic{#1}%
\present{\pregzero{#1}}}

\newcommand{\pregone}[1]{%
\begin{picture}(5,0.4)(0,-0.2)%
\cell{2.5}{0.2}{b}{#1}%
\put(0,0){\vector(1,0){5}}%
\end{picture}}

\mcm{\gone}{1}{%
\ginitdims{5}{0.4}%
\abovepic{#1}%
\present{\pregone{#1}}}

\newcommand{\pregtwo}[3]{%
\begin{picture}(5,3.4)(0,-0.2)%
\cell{2.5}{3.2}{b}{#1}%
\cell{2.5}{-.2}{t}{#2}%
\cell{2.7}{1.5}{l}{#3}%
\qbezier(0,1.5)(2.5,4.5)(5,1.5)%
\qbezier(0,1.5)(2.5,-1.5)(5,1.5)%
\put(5,1.5){\vector(1,-1){0}}%
\put(5,1.5){\vector(1,1){0}}%
\put(2.5,2.5){\vector(0,-1){2}}%
\end{picture}}

\mcm{\gtwo}{3}{%
\ginitdims{5}{3.4}%
\abovepic{#1}%
\belowpic{#2}%
\present{\pregtwo{#1}{#2}{#3}}}

\newcommand{\pregthree}[5]{%
\begin{picture}(5,5.4)(0,-1.2)%
\cell{2.5}{4.2}{b}{#1}%
\cell{1.5}{1.7}{b}{#2}%
\cell{2.5}{-1.2}{t}{#3}%
\cell{2.7}{2.75}{l}{#4}%
\cell{2.7}{0.25}{l}{#5}%
\qbezier(0,1.5)(2.5,6.5)(5,1.5)%
\qbezier(0,1.5)(2.5,-3.5)(5,1.5)%
\put(0,1.5){\vector(1,0){5}}%
\put(2.5,3.5){\vector(0,-1){1.5}}%
\put(2.5,1){\vector(0,-1){1.5}}%
\put(5,1.5){\vector(1,-3){0}}%
\put(5,1.5){\vector(1,3){0}}%
\end{picture}}

\mcm{\gthree}{5}{%
\ginitdims{5}{5.4}%
\abovepic{#1}%
\belowpic{#3}%
\present{\pregthree{#1}{#2}{#3}{#4}{#5}}}

\newcommand{\pregfour}[7]{%
\begin{picture}(5,8.4)(0,-2.7)%
\cell{2.5}{5.7}{b}{#1}%
\cell{1.5}{2.8}{b}{#2}%
\cell{1.5}{0.2}{t}{#3}%
\cell{2.5}{-2.7}{t}{#4}%
\cell{2.7}{4.25}{l}{#5}%
\cell{2.7}{1.5}{l}{#6}%
\cell{2.7}{-1.25}{l}{#7}%
\qbezier(0,1.5)(2.5,9.5)(5,1.5)%
\qbezier(0,1.5)(2.5,4)(5,1.5)%
\qbezier(0,1.5)(2.5,-1)(5,1.5)%
\qbezier(0,1.5)(2.5,-6.5)(5,1.5)%
\put(2.5,5.25){\vector(0,-1){2}}%
\put(2.5,2.5){\vector(0,-1){2}}%
\put(2.5,-0.25){\vector(0,-1){2}}%
\put(5,1.5){\vector(1,-4){0}}%
\put(5,1.5){\vector(4,-3){0}}%
\put(5,1.5){\vector(4,3){0}}%
\put(5,1.5){\vector(1,4){0}}%
\end{picture}}

\mcm{\gfour}{7}{%
\ginitdims{5}{8.4}%
\abovepic{#1}%
\belowpic{#4}%
\present{\pregfour{#1}{#2}{#3}{#4}{#5}{#6}{#7}}}

\newcommand{\pregthreecell}[5]{%
\begin{picture}(8,5)(-4,-2.5)%
\cell{0}{2.5}{b}{#1}%
\cell{0}{-2.5}{t}{#2}%
\cell{-1.7}{0}{r}{#3}%
\cell{1.7}{0}{l}{#4}%
\cell{0}{0.2}{b}{#5}%
\qbezier(-4,0)(0,4.2)(4,0)%
\qbezier(-4,0)(0,-4.2)(4,0)%
\qbezier(-0.5,1.8)(-2.5,0)(-0.5,-1.8)%
\qbezier(0.5,1.8)(2.5,0)(0.5,-1.8)%
\put(-1,0){\vector(1,0){2}}%
\put(4,0){\vector(1,-1){0}}%
\put(4,0){\vector(1,1){0}}%
\put(-0.5,-1.8){\vector(1,-1){0}}%
\put(0.5,-1.8){\vector(-1,-1){0}}%
\end{picture}}

\mcm{\gthreecell}{5}{%
\ginitdims{8}{5}%
\abovepic{#1}%
\belowpic{#2}%
\present{\pregthreecell{#1}{#2}{#3}{#4}{#5}}}

\newcommand{\pregthreecellu}{%
\begin{picture}(5,3.4)(-0.5,-0.2)%
\qbezier(-.5,1.5)(2,4.5)(4.5,1.5)%
\qbezier(-.5,1.5)(2,-1.5)(4.5,1.5)%
\qbezier(1.5,2.7)(0.5,1.5)(1.5,0.3)%
\qbezier(2.5,2.7)(3.5,1.5)(2.5,0.3)%
\put(1.3,1.5){\vector(1,0){1.4}}%
\put(4.5,1.5){\vector(1,-1){0}}%
\put(4.5,1.5){\vector(1,1){0}}%
\put(1.5,0.3){\vector(2,-3){0}}%
\put(2.5,0.3){\vector(-2,-3){0}}%
\end{picture}}

\mcm{\gthreecellu}{0}{%
\ginitdims{5}{3.4}%
\present{\pregthreecellu}}

\newcommand{\pregtwocentre}[3]{%
\begin{picture}(5,3.4)(0,-0.2)%
\cell{2.5}{3.2}{b}{#1}%
\cell{2.5}{-.2}{t}{#2}%
\cell{2.5}{1.5}{c}{#3}%
\qbezier(0,1.5)(2.5,4.5)(5,1.5)%
\qbezier(0,1.5)(2.5,-1.5)(5,1.5)%
\put(5,1.5){\vector(1,-1){0}}%
\put(5,1.5){\vector(1,1){0}}%
\put(2.5,2.5){\vector(0,-1){2}}%
\end{picture}}

\mcm{\gtwocentre}{3}{%
\ginitdims{5}{3.4}%
\abovepic{#1}%
\belowpic{#2}%
\present{\pregtwocentre{#1}{#2}{#3}}}

\newcommand{\pregspecialone}[9]{%
\begin{picture}(8,8)(-4,-4)%
\cell{0}{3.9}{b}{#1}%
\cell{-2}{-0.2}{t}{#2}%
\cell{0}{-3.9}{t}{#3}%
\cell{-1.5}{1.1}{r}{#4}%
\cell{0.2}{1.5}{l}{#5}%
\cell{1.5}{1.1}{l}{#6}%
\cell{0.2}{-2}{l}{#7}%
\cell{-0.9}{2.3}{b}{#8}%
\cell{0.9}{2.3}{b}{#9}%
\qbezier(-4,0)(0,8)(4,0)%
\qbezier(-4,0)(0,-8)(4,0)%
\qbezier(-0.5,3.4)(-3.5,2)(-0.5,0.6)%
\qbezier(0.5,3.4)(3.5,2)(0.5,0.6)%
\put(-4,0){\vector(1,0){8}}%
\put(0,3.4){\vector(0,-1){2.8}}%
\put(0,-0.8){\vector(0,-1){2.4}}%
\put(-1.5,2.2){\vector(1,0){1.2}}%
\put(0.3,2.2){\vector(1,0){1.2}}%
\put(4,0){\vector(1,-2){0}}%
\put(4,0){\vector(1,2){0}}%
\put(-0.5,0.6){\vector(2,-1){0}}%
\put(0.5,0.6){\vector(-2,-1){0}}%
\end{picture}}

\mcm{\gspecialone}{9}{%
\ginitdims{8}{8}%
\abovepic{#1}%
\belowpic{#3}%
\present{\pregspecialone{#1}{#2}{#3}{#4}{#5}{#6}{#7}{#8}{#9}}}

\newcommand{\pregspecialtwo}{%
\begin{picture}(5,3.4)(0,-0.2)%
\qbezier(0,1.5)(2.5,4.5)(5,1.5)%
\qbezier(0,1.5)(2.5,-1.5)(5,1.5)%
\qbezier(1.7,2.5)(0,1.5)(1.7,0.5)%
\qbezier(3.3,2.5)(5,1.5)(3.3,0.5)%
\put(5,1.5){\vector(1,-1){0}}%
\put(5,1.5){\vector(1,1){0}}%
\put(1.7,0.5){\vector(3,-2){0}}%
\put(3.3,0.5){\vector(-3,-2){0}}%
\put(2.5,2.5){\vector(0,-1){2}}%
\put(1.2,1.5){\vector(1,0){1}}%
\put(2.8,1.5){\vector(1,0){1}}%
\end{picture}}

\mcm{\gspecialtwo}{0}{%
\ginitdims{5}{3.4}%
\present{\pregspecialtwo}}

\newcommand{\pregspecialthree}{%
\begin{picture}(5,5.4)(0,-1.2)%
\qbezier(0,1.5)(2.5,6.5)(5,1.5)%
\qbezier(0,1.5)(2.5,-3.5)(5,1.5)%
\qbezier(2,3.5)(1,2.75)(2,2)%
\qbezier(3,3.5)(4,2.75)(3,2)%
\qbezier(2,1)(1,0.25)(2,-0.5)%
\qbezier(3,1)(4,0.25)(3,-0.5)%
\put(0,1.5){\vector(1,0){5}}%
\put(1.5,2.75){\vector(1,0){2}}%
\put(1.5,0.25){\vector(1,0){2}}%
\put(5,1.5){\vector(1,-3){0}}%
\put(5,1.5){\vector(1,3){0}}%
\put(2,2){\vector(1,-1){0}}%
\put(3,2){\vector(-1,-1){0}}%
\put(2,-0.5){\vector(1,-1){0}}%
\put(3,-0.5){\vector(-1,-1){0}}%
\end{picture}}

\mcm{\gspecialthree}{0}{%
\ginitdims{5}{5.4}%
\present{\pregspecialthree}}

\newcommand{\pregonew}[1]{%
\begin{picture}(8,0.4)(0,-0.2)%
\cell{4}{0.2}{b}{#1}%
\put(0,0){\vector(1,0){8}}%
\end{picture}}

\mcm{\gonew}{1}{%
\ginitdims{8}{0.4}%
\abovepic{#1}%
\present{\pregonew{#1}}}

\mcm{\gzersu}{0}{%
\gsinitdims{0}{.6}%
\present{\pregblw{}}}

\mcm{\gonesu}{0}{%
\gsinitdims{5}{0.4}%
\present{\pregone{}}}

\mcm{\gtwosu}{0}{%
\gsinitdims{5}{3.4}%
\present{\pregtwo{}{}{}}}

\mcm{\gthreesu}{0}{%
\gsinitdims{5}{5.4}%
\present{\pregthree{}{}{}{}{}}}

\mcm{\gfoursu}{0}{%
\gsinitdims{5}{8.4}%
\present{\pregfour{}{}{}{}{}{}{}}}

\newcommand{\prectwo}[3]%
{\begin{picture}(4.2,3.4)(-0.1,-0.2)%
\cell{2}{3.2}{b}{#1}%
\cell{2}{-0.2}{t}{#2}%
\cell{2.2}{1.5}{l}{#3}%
\qbezier(0,2)(2,4)(4,2)%
\qbezier(0,1)(2,-1)(4,1)%
\put(4,2){\vector(1,-1){0}}%
\put(4,1){\vector(1,1){0}}%
\put(2,2.5){\vector(0,-1){2}}%
\end{picture}}

\mcm{\ctwo}{3}{%
\cinitdims{4.2}{3.4}%
\abovepic{#1}%
\belowpic{#2}%
\present{\prectwo{#1}{#2}{#3}}}

\newcommand{\precthree}[5]{%
\begin{picture}(4.2,5.4)(-0.1,-0.2)%
\cell{2}{5.2}{b}{#1}%
\cell{1}{2.7}{b}{#2}%
\cell{2}{-.2}{t}{#3}%
\cell{2.2}{3.75}{l}{#4}%
\cell{2.2}{1.25}{l}{#5}%
\qbezier(0,3)(2,7)(4,3)%
\qbezier(0,2)(2,-2)(4,2)%
\put(0,2.5){\vector(1,0){4}}%
\put(2,4.5){\vector(0,-1){1.5}}%
\put(2,2){\vector(0,-1){1.5}}%
\put(4,3){\vector(1,-3){0}}%
\put(4,2){\vector(1,3){0}}%
\end{picture}}

\mcm{\cthree}{5}{%
\cinitdims{4.2}{5.4}%
\abovepic{#1}%
\belowpic{#3}%
\present{\precthree{#1}{#2}{#3}{#4}{#5}}}

\newcommand{\prectwoop}[3]%
{\begin{picture}(4.2,3.4)(-0.1,-0.2)%
\cell{2}{3.2}{b}{#1}%
\cell{2}{-0.2}{t}{#2}%
\cell{2.2}{1.5}{l}{#3}%
\qbezier(0,2)(2,4)(4,2)%
\qbezier(0,1)(2,-1)(4,1)%
\put(0,2){\vector(-1,-1){0}}%
\put(0,1){\vector(-1,1){0}}%
\put(2,2.5){\vector(0,-1){2}}%
\end{picture}}

\mcm{\ctwoop}{3}{%
\cinitdims{4.2}{3.4}%
\abovepic{#1}%
\belowpic{#2}%
\present{\prectwoop{#1}{#2}{#3}}}

\newcommand{\prectwopar}[4]{%
\begin{picture}(4.2,3.4)(-0.1,-0.2)%
\cell{2}{3.2}{b}{#1}%
\cell{2}{-0.2}{t}{#2}%
\cell{1.6}{1.5}{r}{#3}%
\cell{2.4}{1.5}{l}{#4}%
\qbezier(0,2)(2,4)(4,2)%
\qbezier(0,1)(2,-1)(4,1)%
\put(4,2){\vector(1,-1){0}}%
\put(4,1){\vector(1,1){0}}%
\put(1.8,2.5){\vector(0,-1){2}}%
\put(2.2,2.5){\vector(0,-1){2}}%
\end{picture}}

\mcm{\ctwopar}{4}{%
\cinitdims{4.2}{3.4}%
\abovepic{#1}%
\belowpic{#2}%
\present{\prectwopar{#1}{#2}{#3}{#4}}}

\newcommand{\precthreein}[5]{%
\begin{picture}(4.2,5.4)(-0.1,-0.2)%
\cell{2}{5.2}{b}{#1}%
\cell{1}{2.7}{b}{#2}%
\cell{2}{-.2}{t}{#3}%
\cell{2.2}{3.75}{l}{#4}%
\cell{2.2}{1.25}{l}{#5}%
\qbezier(0,3)(2,7)(4,3)%
\qbezier(0,2)(2,-2)(4,2)%
\put(0,2.5){\vector(1,0){4}}%
\put(2,4.5){\vector(0,-1){1.5}}%
\put(2,0.5){\vector(0,1){1.5}}%
\put(4,3){\vector(1,-3){0}}%
\put(4,2){\vector(1,3){0}}%
\end{picture}}

\mcm{\cthreein}{5}{%
\cinitdims{4.2}{5.4}%
\abovepic{#1}%
\belowpic{#3}%
\present{\precthreein{#1}{#2}{#3}{#4}{#5}}}

\newcommand{\precthreecell}[5]{%
\begin{picture}(8.2,5)(-4.1,-2.5)%
\cell{0}{2.5}{b}{#1}%
\cell{0}{-2.5}{t}{#2}%
\cell{-1.7}{0}{r}{#3}%
\cell{1.7}{0}{l}{#4}%
\cell{0}{0.2}{b}{#5}%
\qbezier(-4,0.5)(0,4)(4,0.5)%
\qbezier(-4,-0.5)(0,-4)(4,-0.5)%
\qbezier(-0.5,2)(-2.5,0)(-0.5,-2)%
\qbezier(0.5,2)(2.5,0)(0.5,-2)%
\put(-1,0){\vector(1,0){2}}%
\put(4,0.5){\vector(1,-1){0}}%
\put(4,-0.5){\vector(1,1){0}}%
\put(-0.5,-2){\vector(1,-1){0}}%
\put(0.5,-2){\vector(-1,-1){0}}%
\end{picture}}

\mcm{\cthreecell}{5}{%
\cinitdims{8.2}{5}%
\abovepic{#1}%
\belowpic{#2}%
\present{\precthreecell{#1}{#2}{#3}{#4}{#5}}}

\newcommand{\precthreecellpar}[6]{%
\begin{picture}(8.2,5)(-4.1,-2.5)%
\cell{0}{2.5}{b}{#1}%
\cell{0}{-2.5}{t}{#2}%
\cell{-1.7}{0}{r}{#3}%
\cell{1.7}{0}{l}{#4}%
\cell{0}{0.4}{b}{#5}%
\cell{0}{-0.4}{t}{#6}%
\qbezier(-4,0.5)(0,4)(4,0.5)%
\qbezier(-4,-0.5)(0,-4)(4,-0.5)%
\qbezier(-0.5,2)(-2.5,0)(-0.5,-2)%
\qbezier(0.5,2)(2.5,0)(0.5,-2)%
\put(-1,0.2){\vector(1,0){2}}%
\put(-1,-0.2){\vector(1,0){2}}%
\put(4,0.5){\vector(1,-1){0}}%
\put(4,-0.5){\vector(1,1){0}}%
\put(-0.5,-2){\vector(1,-1){0}}%
\put(0.5,-2){\vector(-1,-1){0}}%
\end{picture}}

\mcm{\cthreecellpar}{6}{%
\cinitdims{8.2}{5}%
\abovepic{#1}%
\belowpic{#2}%
\present{\precthreecellpar{#1}{#2}{#3}{#4}{#5}{#6}}}

\newcommand{\prectwocellv}[5]{%
\begin{picture}(3.4,4.2)(0.8,0.9)%
\cell{2.5}{5.1}{b}{#1}%
\cell{2.5}{0.9}{t}{#2}%
\cell{0.8}{3}{r}{#3}%
\cell{4.2}{3}{l}{#4}%
\cell{2.5}{3.2}{b}{#5}%
\qbezier(2,5)(0,3)(2,1)%
\qbezier(3,5)(5,3)(3,1)%
\put(2,1){\vector(1,-1){0}}%
\put(3,1){\vector(-1,-1){0}}%
\put(1.5,3){\vector(1,0){2}}%
\end{picture}}

\mcm{\ctwocellv}{5}{%
\cinitdims{3.4}{4.2}%
\abovepic{#1}%
\belowpic{#2}%
\sidespic{#3}%
\sidespic{#4}%
\present{\prectwocellv{#1}{#2}{#3}{#4}{#5}}}

\newcommand{\precthreecellv}[7]{%
\begin{picture}(5,8.2)(0.5,-1.6)%
\cell{3}{6.6}{b}{#1}%
\cell{3}{-1.6}{t}{#2}%
\cell{0.5}{2.5}{r}{#3}%
\cell{5.5}{2.5}{l}{#4}%
\cell{3}{4.2}{b}{#5}%
\cell{3}{0.8}{t}{#6}%
\cell{3.2}{2.5}{l}{#7}%
\qbezier(3.5,6.5)(7,2.5)(3.5,-1.5)%
\qbezier(2.5,6.5)(-1,2.5)(2.5,-1.5)%
\put(2.5,-1.5){\vector(1,-1){0}}%
\put(3.5,-1.5){\vector(-1,-1){0}}%
\qbezier(1,3)(3,5)(5,3)%
\qbezier(1,2)(3,0)(5,2)%
\put(5,3){\vector(1,-1){0}}%
\put(5,2){\vector(1,1){0}}%
\put(3,3.5){\vector(0,-1){2}}%
\end{picture}}

\mcm{\cthreecellv}{7}{%
\cinitdims{5}{8.2}%
\abovepic{#1}%
\belowpic{#2}%
\sidespic{#3}%
\sidespic{#4}%
\present{\precthreecellv{#1}{#2}{#3}{#4}{#5}{#6}{#7}}}

\newcommand{\pretopez}[2]{%
\begin{picture}(2.6,2.3)(-1.3,-2.2)%
\cell{0}{-2.2}{t}{#1}%
\cell{0}{-1.2}{c}{#2}%
\qbezier(0,0)(-2,-2)(0,-2)%
\qbezier(0,0)(2,-2)(0,-2)%
\put(0,0){\vector(-1,1){0}}%
\end{picture}}

\mcm{\topez}{2}{%
\ginitdims{2.6}{2.3}%
\belowpic{#1}%
\present{\pretopez{#1}{#2}}}

\newcommand{\pretopea}[3]{%
\begin{picture}(4,1.9)(-2,-0,2)%
\cell{0}{1.7}{b}{#1}%
\cell{0}{-0.2}{t}{#2}%
\cell{0}{0.7}{c}{#3}%
\qbezier(-2,0)(0,3)(2,0)%
\put(-2,0){\vector(1,0){4}}%
\put(2,0){\vector(2,-3){0}}%
\end{picture}}

\mcm{\topea}{3}{%
\ginitdims{4}{1.9}%
\abovepic{#1}%
\belowpic{#2}%
\present{\pretopea{#1}{#2}{#3}}}

\newcommand{\pretopeb}[4]{%
\begin{picture}(4,2.2)(-2,-0.2)%
\cell{-1.1}{1}{br}{#1}%
\cell{1.1}{1}{bl}{#2}%
\cell{0}{-0.2}{t}{#3}%
\cell{0}{0.8}{c}{#4}%
\put(-2,0){\vector(1,1){2}}%
\put(0,2){\vector(1,-1){2}}%
\put(-2,0){\vector(1,0){4}}%
\end{picture}}

\mcm{\topeb}{4}{%
\ginitdims{4}{2.2}%
\belowpic{#3}%
\present{\pretopeb{#1}{#2}{#3}{#4}}}

\newcommand{\pretopec}[5]{%
\begin{picture}(4,2.2)(-2,-0.2)%
\cell{-1.8}{1}{br}{#1}%
\cell{0}{2.2}{b}{#2}%
\cell{1.8}{1}{bl}{#3}%
\cell{0}{-0.2}{t}{#4}%
\cell{0}{0.8}{c}{#5}%
\put(-2,0){\vector(1,2){1}}%
\put(-1,2){\vector(1,0){2}}%
\put(1,2){\vector(1,-2){1}}%
\put(-2,0){\vector(1,0){4}}%
\end{picture}}

\mcm{\topec}{5}{%
\ginitdims{4}{2.2}%
\sidespic{#1}%
\abovepic{#2}%
\sidespic{#3}%
\belowpic{#4}%
\present{\pretopec{#1}{#2}{#3}{#4}{#5}}}

\newcommand{\pretoped}[6]{%
\begin{picture}(4,2.5)(-2,-0.2)%
\cell{-2}{0.6}{br}{#1}%
\cell{-0.7}{2.2}{br}{#2}%
\cell{0.7}{2.2}{bl}{#3}%
\cell{2}{0.6}{bl}{#4}%
\cell{0}{-0.2}{t}{#5}%
\cell{0}{0.8}{c}{#6}%
\put(-2,0){\vector(1,3){0.5}}%
\put(-1.5,1.5){\vector(3,2){1.5}}%
\put(0,2.5){\vector(3,-2){1.5}}%
\put(1.5,1.5){\vector(1,-3){0.5}}%
\put(-2,0){\vector(1,0){4}}%
\end{picture}}

\mcm{\toped}{6}{%
\ginitdims{4}{2.5}%
\sidespic{#1}%
\abovepic{#2}%
\abovepic{#3}%
\sidespic{#4}%
\belowpic{#5}%
\present{\pretoped{#1}{#2}{#3}{#4}{#5}{#6}}}

\newcommand{\pretopeq}[5]{%
\begin{picture}(4,2.5)(-2,-0.2)%
\cell{-2}{0.6}{br}{#1}%
\cell{-1}{2.2}{br}{#2}%
\cell{2}{0.6}{bl}{#3}%
\cell{0}{-0.2}{t}{#4}%
\cell{0}{0.8}{c}{#5}%
\put(-2,0){\vector(1,3){0.5}}%
\put(-1.5,1.5){\vector(1,1){1}}%
\cell{0.9}{2.3}{c}{\ddots}
\put(1.5,1.5){\vector(1,-3){0.5}}%
\put(-2,0){\vector(1,0){4}}%
\end{picture}}

\mcm{\topeq}{5}{%
\ginitdims{4}{2.5}%
\sidespic{#1}%
\abovepic{#2}%
\sidespic{#3}%
\belowpic{#4}%
\present{\pretopeq{#1}{#2}{#3}{#4}{#5}}}

\newcommand{\pretopebase}[1]{%
\begin{picture}(4,0.4)(0,-0.2)%
\cell{2}{0.2}{b}{#1}%
\put(0,0){\vector(1,0){4}}%
\end{picture}}

\mcm{\topebase}{1}{%
\ginitdims{4}{0.4}%
\abovepic{#1}%
\present{\pretopebase{#1}}}

\newcommand{\pretopezs}[2]{%
\begin{picture}(2.6,2.3)(-1.3,-2.2)%
\cell{0}{-2.2}{t}{#1}%
\cell{0}{-1.2}{c}{#2}%
\qbezier(0,0)(-2,-2)(0,-2)%
\qbezier(0,0)(2,-2)(0,-2)%
\end{picture}}

\mcm{\topezs}{2}{%
\ginitdims{2.6}{2.3}%
\belowpic{#1}%
\present{\pretopezs{#1}{#2}}}

\newcommand{\pretopeas}[3]{%
\begin{picture}(4,1.9)(-2,-0,2)%
\cell{0}{1.7}{b}{#1}%
\cell{0}{-0.2}{t}{#2}%
\cell{0}{0.7}{c}{#3}%
\qbezier(-2,0)(0,3)(2,0)%
\put(-2,0){\line(1,0){4}}%
\end{picture}}

\mcm{\topeas}{3}{%
\ginitdims{4}{1.9}%
\abovepic{#1}%
\belowpic{#2}%
\present{\pretopeas{#1}{#2}{#3}}}

\newcommand{\pretopebs}[4]{%
\begin{picture}(4,2.2)(-2,-0.2)%
\cell{-1.1}{1}{br}{#1}%
\cell{1.1}{1}{bl}{#2}%
\cell{0}{-0.2}{t}{#3}%
\cell{0}{0.8}{c}{#4}%
\put(-2,0){\line(1,1){2}}%
\put(0,2){\line(1,-1){2}}%
\put(-2,0){\line(1,0){4}}%
\end{picture}}

\mcm{\topebs}{4}{%
\ginitdims{4}{2.2}%
\belowpic{#3}%
\present{\pretopebs{#1}{#2}{#3}{#4}}}

\newcommand{\pretopecs}[5]{%
\begin{picture}(4,2.2)(-2,-0.2)%
\cell{-1.8}{1}{br}{#1}%
\cell{0}{2.2}{b}{#2}%
\cell{1.8}{1}{bl}{#3}%
\cell{0}{-0.2}{t}{#4}%
\cell{0}{0.8}{c}{#5}%
\put(-2,0){\line(1,2){1}}%
\put(-1,2){\line(1,0){2}}%
\put(1,2){\line(1,-2){1}}%
\put(-2,0){\line(1,0){4}}%
\end{picture}}

\mcm{\topecs}{5}{%
\ginitdims{4}{2.2}%
\sidespic{#1}%
\abovepic{#2}%
\sidespic{#3}%
\belowpic{#4}%
\present{\pretopecs{#1}{#2}{#3}{#4}{#5}}}

\newcommand{\pretopeds}[6]{%
\begin{picture}(4,2.5)(-2,-0.2)%
\cell{-2}{0.6}{br}{#1}%
\cell{-0.7}{2.2}{br}{#2}%
\cell{0.7}{2.2}{bl}{#3}%
\cell{2}{0.6}{bl}{#4}%
\cell{0}{-0.2}{t}{#5}%
\cell{0}{0.8}{c}{#6}%
\put(-2,0){\line(1,3){0.5}}%
\put(-1.5,1.5){\line(3,2){1.5}}%
\put(0,2.5){\line(3,-2){1.5}}%
\put(1.5,1.5){\line(1,-3){0.5}}%
\put(-2,0){\line(1,0){4}}%
\end{picture}}

\mcm{\topeds}{6}{%
\ginitdims{4}{2.5}%
\sidespic{#1}%
\abovepic{#2}%
\abovepic{#3}%
\sidespic{#4}%
\belowpic{#5}%
\present{\pretopeds{#1}{#2}{#3}{#4}{#5}{#6}}}

\newcommand{\pretopeqs}[5]{%
\begin{picture}(4,2.5)(-2,-0.2)%
\cell{-2}{0.6}{br}{#1}%
\cell{-1}{2.2}{br}{#2}%
\cell{2}{0.6}{bl}{#3}%
\cell{0}{-0.2}{t}{#4}%
\cell{0}{0.8}{c}{#5}%
\put(-2,0){\line(1,3){0.5}}%
\put(-1.5,1.5){\line(1,1){1}}%
\cell{0.9}{2.3}{c}{\ddots}
\put(1.5,1.5){\line(1,-3){0.5}}%
\put(-2,0){\line(1,0){4}}%
\end{picture}}

\mcm{\topeqs}{5}{%
\ginitdims{4}{2.5}%
\sidespic{#1}%
\abovepic{#2}%
\sidespic{#3}%
\belowpic{#4}%
\present{\pretopeqs{#1}{#2}{#3}{#4}{#5}}}

\newcommand{\pretopebases}[1]{%
\begin{picture}(4,0.4)(0,-0.2)%
\cell{2}{0.2}{b}{#1}%
\put(0,0){\line(1,0){4}}%
\end{picture}}

\mcm{\topebases}{1}{%
\ginitdims{4}{0.4}%
\abovepic{#1}%
\present{\pretopebases{#1}}}

%% file: esstitle.tex
\title{Structures in Higher-Dimensional Category Theory}
\author{Tom Leinster\\ \\
	\normalsize{Department of Pure Mathematics, University of
	Cambridge}\\ 
	\normalsize{Email: leinster@dpmms.cam.ac.uk}\\
	\normalsize{Web: http://www.dpmms.cam.ac.uk/$\sim$leinster}}
\date{13 January, 1998\\ {\small Revised June 98, October 98}\\
	\vspace{8mm}
	\normalsize
	\textbf{Abstract}\\
	\vspace{3mm} \raggedright \setlength{\rightskip}{0pt}
This is an exposition of some of the constructions which have arisen in
higher-dimensional category theory. We start with a review of the general
theory of operads and multicategories (cf.\ \cite{Bur}, \cite{Her},
\cite{GOM}). Using this we give a new definition of $n$-category (a variation
on Batanin's); we also give an informal definition in pictures.  Next we
discuss Gray-categories and their place in coherence problems. Finally, we
present various constructions relevant to the opetopic definitions of
$n$-category.\newline\\
\vspace{-2mm}
New material includes our definition of lax $n$-category; a suggestion for a
definition of lax cubical $n$-category; a characterization of small
Gray-categories as the small substructures of \tooc; a conjecture on
coherence theorems in higher dimensions; a construction of the category of
trees and, more generally, of $n$-pasting diagrams; and an analogue of the
Baez-Dolan slicing process in the general theory of operads.\newline\\}

\maketitle

%% file: preface.tex
\section*{Preface to the arXiv version}

This paper was written in early 1998.  At the time I did not post it on the
electronic archive because the files, containing as they do many
memory-intensive graphics, were considered too large; and later I did not post
it because there were various shortcomings that I wanted to make good before
committing it to permanent storage.  But technology has progressed, and I
have become wise enough to know that I will never get round to revising this
paper in its present form: so here it is.

The main aspects I would like to have changed are as follows.  The choice of
the term `lax $n$-category' was probably misguided, and I would now
substitute `weak $n$-category'.  Chapters~\ref{ch:om} and~\ref{ch:glob} are
mostly superseded by my thesis (\texttt{math.CT/0011106}), where the ideas
are explained more precisely and, I think, more clearly.
Chapters~\ref{ch:gray} and~\ref{ch:opap}, however, do not yet appear in any
form elsewhere.  In Chapter~\ref{ch:gray} I should have explained more
clearly the non-standard usage of \tooc, which denotes 2-categories,
homomorphisms, strong transformations and modifications, and in
section~\ref{sec:hd-coh-conj} I would have liked to add some further points
and subtract at least one dubious statement.  Finally, the subject has
progressed significantly since the time of writing and there are many new
references to be added, in particular to the work of Cheng relating the
opetopes of Chapter~\ref{ch:opap} to the opetopes of~\cite{BD} and the
multitopes of~\cite{HMP}.  The bibliography of my thesis or (better) of my
survey paper \texttt{math.CT/0107188} should provide a decent substitute.

\subsection*{}

\begin{flushright}
September 2001
\end{flushright}

%% file: essintro.tex
\chapter*{Introduction}
\ucontents{chapter}{Introduction}

The subject of this essay is higher-dimensional category theory; the theme is
to make sense of the diverse structures which have arisen in its pursuit. For
instance: various different notions of `operad' have been developed, and we
provide a general theory covering many of them. Or: 
\index{Gray-category}%
Gray-categories were
introduced to address questions of coherence for tricategories, and here we
give two characterizations of them (both different from the original
definition) in an attempt to see what the pattern is for coherence in higher
dimensions.

The description of mathematical structures in primitive terms constitutes
much of what is here. Typically, we give informal arguments using pictures
(for example, the 
\index{Batanin|(}%
Batanin operads of Chapter~\ref{ch:glob}). Where proofs are necessary---and
there are few big theorems or startling results---they tend to be bland,
seemingly formal and low in content. We usually omit these. 

To date there appear to have been three main approaches to the problem of
defining `lax 
\index{n-category@$n$-category!various definitions of}%
$n$-category'. (See also the introduction to 
\index{Simpson}%
\cite{SimCMS} for some notes on the evolution of the subject.) Zouhair
\index{Tamsamani}%
Tamsamani has made a definition which generalizes the correspondence between
a category and its 
\index{nerve}%
nerve $\Delta^{\op}\go\Set$ (\cite{Tam}); we say nothing
about this, but for some speculation on page~\pageref{p:tam-rmk}. 
\index{Baez!-Dolan}%
Baez and Dolan presented their definition in terms of
\index{opetopes!and $n$-categories}%
opetopes, and this
inspired a related description by 
\index{Hermida!-Makkai-Power}%
Hermida, Makkai and Power (\cite{BD}, \cite{HMP}). Finally, Batanin defined
lax $n$- and $\omega$-categories using globular structures (\cite{Bat}). On
the second and third approaches we have much more to say.

The contents of this essay are as follows.  A preliminary chapter reviews
some basic bicategory theory, which we shall call upon later.
Chapter~\ref{ch:om} presents a general theory of operads and
multicategories. The language alone allows us to give a very concise
definition of a Batanin operad; in Chapter~\ref{ch:glob} we explain what a
Batanin operad is in elementary language, and how lax $n$-categories may then
be defined. We also include a parallel (but truncated) development for
cubical categories, the higher-dimensional versions of double categories.
Chapter~\ref{ch:gray} is on Gray-categories, important in coherence for
tricategories. Some ideas from Chapter~\ref{ch:glob} are employed here.
Chapter~\ref{ch:opap}, however, calls only upon Chapter~\ref{ch:om}: using
some of the results there, various opetopic structures and methods are
illuminated. There is no attempt to describe the opetopic approach as a
whole; the theory of Chapter~\ref{ch:om} seems less well suited to this
style than to Batanin's.

\subsection*{Terminology}

We distinguish between 
\index{strict!vs lax@\emph{vs.} lax}%
\index{lax vs strict@lax \emph{vs.} strict}%
\emph{strict}\/ $n$-categories at one end of the scale, and
\emph{lax}\/ $n$-categories at the other. Some authors have used
\emph{weak}\/ instead of lax; `lax' seems more evocative and a better
opposite to `strict'.
\index{Baez}%
Baez has argued convincingly that of the strict and lax
versions, the lax are more commonly-occurring and more natural as an idea, so
this ought to be the default in the nomenclature. Thus he uses the term
`$n$-category' in contrast to `strict $n$-category'. However, the terminology
for dimension 2 is quite well-established, where bicategories are the lax
version and 2-categories the strict. Our policy is always to say `lax' or
`strict' except in dimensions 2 and 3, and there use \emph{2-category}\/ and
\emph{3-category}\/ for the strict versions and \emph{bicategory} and
\emph{tricategory} for the lax.

For us, an \emph{operad}\/ is a multicategory with just one object; this
conflicts with \cite{BD}, as explained after Definition~\ref{defn:multicat}.

\subsection*{Related Work}

\emph{Preliminaries on Bicategories.}\/ The basic definitions are taken from 
\index{B\'{e}nabou}%
\cite{Ben} and
\index{Gray}%
\cite{Gray}, and the outline of the coherence theorem from 
\index{Street}%
\cite{StFB} and 
\index{Gordon-Power-Street}%
\cite{GPS}. The same material is covered in more detail in
\cite{BB}; see also 
\index{Lack}%
\cite{Lack} for another
summary. I have not seen the ideas on `bias' elsewhere; however, 
\index{Tamsamani}%
Tamsamani appears to establish that his definition of
$n$-category is `equivalent to' the usual one in the case $n=2$ (\cite{Tam}),
and this must require some consideration of the bias issue.

\emph{\ref{ch:om}: Operads and Multicategories.}\/ The material here was first
published as \cite{GOM}, of which this chapter is an abbreviated version. At
that time the ideas were new to the author; however, the definition of
\Cartpr-multicategory was also given by 
\index{Burroni}%
Burroni in 1971 and 
\index{Hermida}%
Hermida in 1997
(\cite{Bur}, \cite{Her}). Some of the other ideas in the chapter also appear
in one or both of these sources. The condition that a monad be cartesian is
close to a condition in 
\index{Carboni-Johnstone}%
\cite{CJ}, results from which on familial
representability are used here and appear to bear on the construction of the
free strict $\omega$-category functor in~\ref{sec:trees}. It
seems that 
\index{Kelly}%
Kelly's theory of clubs (\cite{KelAA}, \cite{KelMV}, \cite{KelCDC}) also has
common ground with this chapter.

\emph{\ref{ch:glob}: The Globular Approach.}\/ This is a variant of Batanin's
definition of lax $n$-category, which appears in \cite{Bat} and is summarized
in 
\index{Street}%
\cite{StRMB}. The operads used in this chapter are the same as Batanin's, but
our notions of contractibility differ. The final section, on cubical
structures, is to my knowledge original.

\emph{\ref{ch:gray}: Gray-categories.}\/ These were defined in \cite{GPS}
following
\index{Gray}%
\cite{Gray}. We use an equivalent definition, given in
\cite{Bat}
\index{Batanin|)}%
. Sections \ref{sec:tooc-gray} and \ref{sec:gray-tri} are based on
these sources, with what appears to be a new emphasis (on which processes
are canonical). Sections \ref{sec:char-of-gray} and \ref{sec:hd-coh-conj}, on
Cayley representation, seem to be new.

\emph{\ref{ch:opap}: The Opetopic Approach.}\/ Opetopes seem to have been
defined first in 
\index{Baez!-Dolan}%
\cite{BD}; they are also explained in 
\index{Baez}%
\cite{Baez} and used
in the 
\index{Hermida!-Makkai-Power}%
\cite{HMP} approach to lax $n$-categories. My understanding of the
latter is based on 
\index{Hyland}%
\cite{Hy}. Various ideas in this chapter also appear in 
\index{Hermida}%
\cite{Her}. Categories of trees are employed in 
\index{Borcherds}%
\cite{Bor},
\index{Snydal}%
\cite{Sny}, 
\index{Kontsevich-Manin}%
\cite{KMGWC}, \cite{KMQCP} and 
\index{Soibelman}%
\cite{Soi}.

\emph{Other.}\/ The 
\index{Tamsamani}%
Tamsamani approach is laid out in \cite{Tam} and explored
further in 
\index{Simpson}%
Simpson's papers \cite{SimCMS} and \cite{SimLn}. Another set of questions
about lax structures is suggested by the 
\index{relaxed multilinear category}%
relaxed multilinear categories of
\index{Borcherds}%
\index{Snydal}%
Borcherds (\cite{Bor}, \cite{Sny}) and the very similar 
\index{pseudo-monoidal
category}%
pseudo-monoidal categories of
\index{Soibelman}%
Soibelman (\cite{Soi}); see also \ref{eg:multicats}(\ref{eg:tree-mti}),
(\ref{eg:soib-op}) and \ref{sec:cat-trees}.

\subsection*{Acknowledgements}

This work was supported by a PhD grant from EPSRC. The document was prepared
in \LaTeX, using Paul 
\index{Taylor}%
Taylor's diagrams package---essentially, for those diagrams without
curves. It was originally written for the J. T. Knight essay competition at
Cambridge University. 

I would like to thank Peter 
\index{Johnstone}%
Johnstone, Craig 
\index{Snydal}%
Snydal, and especially Martin
\index{Hyland}%
Hyland, for the help, direction and encouragement they have continued to give
me. John
\index{Power}%
Power and Ross
\index{Street}%
Street have imparted valuable advice on tricategories, and I have had useful
discussions with Robin
\index{Cockett}%
Cockett on structured categories.

In previous versions of this paper, the second chapter was billed as an
account of 
\index{Batanin!misrepresentation}%
Batanin's definition of lax $\omega$-category. This is not the case. I thank
Michael Batanin for pointing this out to me, and apologise for the previous
misrepresentation of his work. The introduction to Chapter~\ref{ch:glob}
explains more fully how our two definitions differ.

%% file: essprelims.tex
\chapter*{Preliminaries on Bicategories}
\ucontents{chapter}{Preliminaries on Bicategories}
\index{bicategory|(}

Here we review the basic properties of bicategories and state our
terminology. The final section, on 
\index{bias}%
`biased' \emph{vs.} `unbiased' bicategories, looks
at what happens when composition of arities 2 and 0 is replaced by
composition of arbitrary arities.

\cite{BB} is a more detailed summary of most of this material.

\subsection*{Basic terminology}

We will typically denote 0-cells of a bicategory \bee\ by $A$, $B$, \ldots,
1-cells by $f$, $g$, \ldots\ and 2-cells by $\alpha$, $\beta$, \ldots, e.g.\
$A\ctwo{f}{g}{\alpha}B$. The
\index{composition!vertical}%
`vertical' composite of 2-cells
$\cdot\cthree{}{}{}{\alpha}{\beta}\cdot$ is written $\beta\of\alpha$ or
$\beta\alpha$, and the 
\index{composition!horizontal}%
`horizontal' composite of 2-cells
\[
\cdot\ctwo{}{}{\alpha}\cdot\ctwo{}{}{\alpha'}\cdot
\]
is written
$\alpha'*\alpha$. 

A 
\index{morphism of bicategories|emph}%
\emph{morphism}\/ $\bee\go\beep$ consists of a map $F$ of the underlying
graphs together with 
\index{coherence!data}%
`coherence' 2-cells $Fg\of Ff \go F(g\of f)$ and $1\go
F1$ satisfying some axioms. If these 2-cells are all invertible then $F$ is
called a 
\index{homomorphism!of bicategories|emph}%
\emph{homomorphism}; if they are identities (so that $Fg\of Ff =
F(g\of f)$ and $F1=1$) then $F$ is called a 
\index{strict!homomorphism of bicategories|emph}%
\emph{strict homomorphism}.

A 
\index{transformation of bicategories|emph}%
\emph{transformation}\/ \bee\ctwo{F}{G}{\sigma}\beep\ between morphisms
consists of 1-cells \gobyc{\sigma_{A}}{FA}{GA} and `coherence' 2-cells
\begin{ntdiag}
FA		&	&\rTo^{Ff}	&	&FB		\\
		&	&		&\	&		\\
\dTo<{\sigma_A}	&	&\ruTo>{\sigma_f}&	&\dTo>{\sigma_B}	\\
		&\	&		&	&		\\
GA		&	&\rTo_{Gf}	&	&GB		\\
\end{ntdiag}
satisfying axioms. If the $\sigma_f$ are all invertible [identities] then
$\sigma$ is a 
\index{strong transformation of bicategories|emph}%
\index{strict!transformation of bicategories|emph}%
\emph{strong [strict] transformation}.

A 
\index{modification|emph}%
\emph{modification}\/ $F\ctwo{\sigma}{\twid{\sigma}}{}G$ between
transformations consists of 2-cells $FA\ctwo{\sigma_A}{\twid{\sigma}_A}{}GA$
satisfying axioms.

We will chiefly be interested in the versions where the coherence cells are
isomorphisms: that is, homomorphisms and strong transformations (and
modifications). Plain morphisms and transformations will not be used at all.

\subsection*{Duality}

Given a bicategory \bee, we may form a 
\index{dual bicategory|emph}%
\index{opposite bicategory|emph}%
dual bicategory $\bee^\op$ by
reversing the 1-cells but not the 2-cells. Thus if \bee\ has a 2-cell
$A\ctwo{f}{g}{\alpha}B$ then $\bee^\op$ has a 2-cell 
$A\ctwoop{f}{g}{\alpha}B$.

\subsection*{Internal features}

We can make certain definitions in an arbitrary bicategory \bee\ by
generalizing from the case $\bee=\Cat$. An 
\index{equivalence!internal|emph}%
\emph{(internal) equivalence}\/ in
\bee\ consists of a pair of 1-cells $A\oppair{f}{g}B$
together with invertible 2-cells $1\go g\of f$ and $f\of g \go 1$. We also
say that $f$ is an equivalence.

A 
\index{monad!in bicategory|emph}%
\emph{monad}\/ in a bicategory consists of a 0-cell $A$, a 1-cell $A
\goby{t} A$, and 2-cells
\[
A\cthreein{1}{t}{t\of t}{\eta}{\mu}A,
\]
such that the diagrams
\[\left.
\begin{diagram}
t\of 1	&\rTo^{t\eta}	&t\of t		&\lTo^{\eta t}	&1\of t	\\
   	&\rdTo<{\diso}	&\dTo<{\mu}	&\ldTo>{\diso}	&	\\
   	&     		&t  		&     		&   	\\
\end{diagram}
\right.
\ \ \ 
\left.
\begin{diagram}
   		&     		&t\of(t\of t)	&\rTo^{t\mu}	&t\of t\\
   		&\ldLine<{\diso}&   		&     		&   \\
(t\of t)\of t	&     		&   		& 		&\dTo>{\mu}\\
\dTo<{\mu t}	&     		&   		&     		&   \\
t\of t		&     		&\rTo_{\mu}	&     		&t  \\
\end{diagram}
\right.\]
commute.

\index{functor!bicategory|emph}%
\subsection*{Functor bicategories}

Given bicategories \bee\ and \beep, there is a bicategory
\ftrcat{\bee}{\beep} with 0-cells homomorphisms, 1-cells strong
transformations and 2-cells modifications. This is a 2-category if \beep\ is.
In particular, \Pshf\ is a 2-category for any \bee.

\index{biequivalence|emph}%
\subsection*{Biequivalence}		\label{p:bieq}

Let \bee\ and \beep\ be bicategories. A \emph{biequivalence} from \bee\ to
\beep\ consists of a pair of homomorphisms \bee\oppair{F}{G}\beep\
together with an equivalence $1 \go G \of F$ inside the bicategory 
\ftrbi{\bee}{\bee} and an equivalence $F \of G \go 1$ inside
\ftrbi{\beep}{\beep}. We also say that $F$ is a biequivalence and that \bee\
is biequivalent to \beep. Just as for equivalence of plain categories, there
is an alternative criterion for biequivalence: namely, that a homomorphism
\gobyc{F}{\bee}{\beep} is a biequivalence if and only if $F$ is 
\index{local equivalence|emph}%
locally an equivalence and is 
\index{surjection-up-to-equivalence|emph}%
surjective-up-to-equivalence on objects. The
former condition means that each functor $\homset{\bee}{A}{B}\go
\homset{\beep}{FA}{FB}$ is an equivalence; the latter
that if $B'$ is any 0-cell of \beep\ then there is some 0-cell $B$ of \bee\
such that $FB$ is (internally) equivalent to $B'$.

\index{coherence!for bicategories|(}%
\index{Yoneda embedding!for bicategories|(}%
\subsection*{Coherence}

We can construct representables for a bicategory \bee: for each 0-cell $A$, a
homomorphism \gobyc{\homset{\bee}{\dashbk}{A}}{\bee^\op}{\Cat}, and similarly
1- and 2-cells. Thus we get a Yoneda homomorphism \gobyc{Y}{\bee}{\Pshf}. It
is straightforward to calculate that $Y$ is locally an equivalence, by
showing that locally it is full, faithful and essentially surjective on
objects. From this it follows that $Y$ provides a biequivalence from \bee\ to
its full image, a 2-category. Hence every bicategory is biequivalent to a
2-category.  We call this `the coherence theorem for bicategories'; it
implies that in a suitable sense, every diagram of coherence 2-cells in a
bicategory commutes.%
\index{coherence!for bicategories|)}%
\index{Yoneda embedding!for bicategories|)}%

\index{bias|(}%
\subsection*{Bias}		\label{sec:bias}

The traditional definition of a bicategory is `biased' towards binary and
nullary compositions, in that only these are given explicit mention. For
instance, there is no specified ternary composite of 1-cells,
$\triple{h}{g}{f} \goesto hgf$, only the derived ones like $h(gf)$ and
$((h1)g)(f1)$. It is necessary to be biased in order to achieve a 
\index{finite axiomatization}%
finite axiomatization. However, such finite axiomatizations become
progressively more complex in higher dimensions (see
\index{Gordon-Power-Street}%
\cite{GPS}, for instance), so the
prevailing approach to defining higher-dimensional categories is an
`unbiased' one, treating all arities even-handedly (usually via 
\index{operad}%
operads).

We would therefore like to define `unbiased bicategory' and see how this
notion compares to the traditional one.

First, let us define an \emph{unbiased category}. The data consists of a
collection $\cat{C}_0$ of objects, a collection \homset{\cat{C}}{A}{B} for
each $A,B\elt\cat{C}_0$, and then for each sequence
$\mb{A}=\tuple{A_0}{A_1}{A_n}$ of objects ($n\geq0$), a function
\[
\gobyc{c_{\mb{A}}}{\homset{\cat{C}}{A_{n-1}}{A_n} \times\cdots\times
\homset{\cat{C}}{A_0}{A_1}}{\homset{\cat{C}}{A_0}{A_n}}.
\]
One then defines the 
\index{term}%
\index{one-term@1-term|emph}%
\emph{1-terms}\/ to be the class of functions generated
from the $c_{\mb{A}}$'s by (binary and nullary) products and
composition. The axioms are that if 
\[
\homset{\cat{C}}{A_{n-1}}{A_n} \times\cdots\times \homset{\cat{C}}{A_0}{A_1}
\parpair{s}{t}
\homset{\cat{C}}{A_0}{A_n}
\]
are two 1-terms, then $s=t$. Unbiased functors and natural transformations
are defined similarly, to make a 2-category \fcat{U\hyph Cat}. Clearly an
unbiased category is uniquely determined by its binary and nullary
compositions, so the forgetful 2-functor (=strict homomorphism) $\fcat{U\hyph
Cat}\go\Cat$ is an isomorphism.

Next we can define an \emph{unbiased bicategory}, in the same style as an
unbiased category, and \emph{unbiased homomorphisms}, \ldots. Any unbiased
bicategory has an underlying bicategory; conversely, we may non-canonically
choose a process by which any bicategory extends to an unbiased bicategory.
The crux is that any unbiased bicategory is isomorphic to one coming from a
traditional bicategory, in the sense of there being an invertible
(non-strict) unbiased homomorphism between the two. Note that this is two
levels better than we might have asked for: it's an isomorphism, not just an
(unbiased) equivalence or biequivalence. With suitable definitions in place,
this translates to a statement about the similarity of the tricategories
\Bicat\ and \fcat{U\hyph Bicat}.

\index{bias|)}
\index{bicategory|)}

%% file: essom.tex
\index{operad|(}
\index{multicategory|(}
\chapter{Operads and Multicategories}	\label{ch:om}

In this chapter we introduce the language of operads and multicategories to
be used in the rest of the essay. The simplest kind of operad---a 
\index{operad!plain|emph}%
\emph{plain
operad}--- consists of a sequence $C(0)$, $C(1)$, \ldots of sets together
with an `identity' element of $C(1)$ and `composition' functions
\[
C(n_{1})\times\cdots\times C(n_{k})\times C(k)
\go C(n_{1}+\cdots +n_{k}),
\]
obeying associativity and identity laws. (In the original definition,
\index{May}%
\cite{May}, the $C(n)$'s were not just sets but spaces with 
\index{symmetric group action}%
symmetric group action.) The simplest kind of multicategory consists of a
collection $C_{0}$ of objects, and arrows like
\[
s_{1}, \ldots, s_{n} \goby{a} s
\]
($s_{1}, \ldots, s_{n}, s \elt C_{0}$), together with composition functions
and identity elements obeying associativity and unit laws. These will be
called 
\index{multicategory!plain|emph}%
\emph{plain multicategories}; a plain operad is therefore a one-object
plain multicategory.

The general idea now is that there's nothing special about \emph{sequences}
of objects: the domain of an arrow might form another shape instead, such as
a tree of objects or just a single object (as in a normal category). Indeed,
the objects do not even need to form a set. Maybe a graph or a category would
do just as well. Together, what these generalizations amount to is the
replacement of the
\index{monoid!free}%
free-monoid monad on \Set\ with some other monad on some
other category.

This generalization is put into practice as follows. The graph structure of a
plain multicategory is a diagram
\begin{slopeydiag}
	&	&C_1	&	&	\\
	&\ldTo<{\dom}	&	&\rdTo>{\cod}	&	\\
C_{0}^{*}&	&	&	&C_0	\\
\end{slopeydiag}
in \Set, where \stbk\ is the free-monoid monad. Now, just as a (small) 
category can be described as a diagram 
\begin{slopeydiag}
		&	&D_1		&	&	\\
		&\ldTo	&		&\rdTo	&	\\
D_0		&	&		&	&D_0	\\
\end{slopeydiag}
in \Sets\ together with identity and composition functions
\[
D_{0}\go D_{1} \mr{,\ \ } 
D_{1}\times_{D_{0}}D_{1}\go D_{1}
\]
satisfying some axioms, so we may describe the multicategory structure on
\gph{C_1}{C_0} 
by manipulation of certain diagrams in \Sets. In general, we take a category
\ess\ and a monad \ust\ on \ess\ satisfying some simple conditions, and
define `\Cartpr-multicategory'. Thus a category is a
\pr{\Set}{\id}-multicategory.

Section~\ref{sec:cart-mnds} describes the simple conditions on \ess\ and
\ust\ required in order that everything that follows will work. Many examples
are given. Section~\ref{sec:multicats} explains what an \Cartpr-multicategory
is and how the examples relate to existing notions of multicategory. In
particular, a concise definition of Batanin operads is given. Most of these
existing notions carry with them the concept of an \emph{algebra} for an
operad/multicategory; \ref{sec:algs} defines algebras in our general
setting. Section~\ref{sec:struc} is on `\Cartpr-structured categories', which
are to \Cartpr-multicategories as strict monoidal categories are to plain
multicategories; in~\ref{sec:free-multi} we sketch the construction of the
free multicategory on a graph. The last two sections are included because of
their impact on the opetopic approach to $n$-categories: they enable, for
instance, a compact construction of the opetopes (\ref{sec:opetopes}). They
will not be used in Chapters~\ref{ch:glob} and~\ref{ch:gray}.

This chapter is an abbreviated version of~\cite{GOM}.

\section{Cartesian Monads}	\label{sec:cart-mnds}

In this section we introduce the conditions required of a monad \Mnd\ on a
category \ess, in order that we may (in~\ref{sec:multicats}) define
the notion of an \Cartpr-multicategory. The conditions are that the category
and the monad are both cartesian, as defined now.

\begin{defn}
\index{cartesian!category|emph}%
A category is called \emph{cartesian} if it has all finite limits.
\end{defn}

\begin{defn}		\label{defn-cart}
\index{cartesian!monad|emph}%
A monad \Mnd\ on a category \cat{S} is called \emph{cartesian} if
\begin{enumerate}
\item $\eta$ and $\mu$ are 
\index{cartesian!natural transformation|emph}%
cartesian natural transformations, i.e.\
for any $X \goby{f} Y$ in \cat{S} the naturality squares
\[\left.
\begin{diagram}
X	&\rTo^{\eta_{X}}&X^{*}		\\
\dTo<{f}&		&\dTo>{f^{*}}	\\
Y	&\rTo^{\eta_{Y}}&Y^{*}		\\
\end{diagram}
\right.
\mi{and}
\left.
\begin{diagram}
X^{**}		&\rTo^{\mu_{X}}	&X^{*}		\\
\dTo<{f^{**}}	&		&\dTo>{f^{*}}	\\
Y^{**}		&\rTo^{\mu_{Y}}	&Y^{*}		\\
\end{diagram}
\right.\]
are pullbacks, and \label{nts-cart}
\item \stbk\ preserves pullbacks. \label{pb-pres}
\end{enumerate}
\end{defn}

\begin{eg}{eg:cart-mnds}

\item The identity monad on any category is clearly cartesian.

\item \label{free-monoid-is-cart}
\index{monoid!free}%
Let $\ess = \Sets$ and let \ust\ be the monoid monad, i.e. the
monad arising from the adjunction
\begin{diagram}
\fcat{Monoid}	&\pile{\rTo \\ \ \ \top \\ \lTo}	&\Sets.
\end{diagram}
Certainly \ess\ is cartesian. It is easy to calculate that \ust, too, is
cartesian (\cite[1.4(ii)]{GOM}).

\item \label{eg:comm-not-cart}
\index{monoid!commutative}%
A non-example. Let $\ess=\Sets$ and let \Mnd\ be the free 
commutative monoid
monad. This is not cartesian: e.g.\ the naturality square for $\mu$ at $2\go
1$ is not a pullback. 

\item \label{eg:sr-theories}
\index{finitary!algebraic theory}%
\index{strongly regular|(emph}%
Let $\ess = \Sets$. Any finitary algebraic theory gives a monad
on \ess; which are cartesian? Without answering this question
completely, we indicate a certain class of theories which do give
cartesian monads. An equation (made up of variables and finitary
operators) is said to be \emph{strongly regular} if the same variables
appear in the same order, without repetition, on each side. Thus
\[\begin{array}{ccc}
(x.y).z=x.(y.z) & \mr{and} & (x\uparrow y)\uparrow z = x\uparrow(y.z),\\
\end{array}\]
but not
\[\begin{array}{cccc}
x+(y+(-y))=x, & x.y=y.x & \mr{or} & (x.x).y=x.(x.y),\\
\end{array}\]
qualify. A theory is called \emph{strongly regular} if it can be presented by
operators and strongly regular equations. In
Example~(\ref{free-monoid-is-cart}), the only property of the theory of
monoids that we actually needed was its strong regularity: for in general,
the monad yielded by any strongly regular theory is cartesian.

This last result, and the notion of strong regularity, are due to
\index{Carboni-Johnstone}%
Carboni and Johnstone. They show in \cite[Proposition 3.2 via Theorem
2.6]{CJ} that a theory is strongly regular iff $\eta$ and $\mu$ are cartesian
natural transformations and $\stbk$ preserves wide pullbacks. A 
\index{wide pullback|emph}%
\index{pullback!wide|emph}%
\emph{wide pullback} is by definition a limit of shape
\begin{diagram}
\cdot	&\cdot	&\cdot		&	&\cdots	&		&\cdot	\\
	&\rdTo(4,2)&\rdTo(3,2)	&\rdTo(2,2)&	&\ldTo(2,2)	&	\\
	&	&		&	&\cdot	&		&,	\\
\end{diagram}
where the top row is a set of any size (perhaps infinite). When the set is of
size 2 this is an ordinary pullback, so the monad from a strongly regular
theory is indeed cartesian. (Examples~(\ref{eg:exceptions-mnd}) and
(\ref{eg:tree-mnd}) can also be found in \cite{CJ}.)%
\index{strongly regular|)emph}

\item \label{eg:exceptions-mnd}
Let $\ess = \Sets$, and let + denote binary coproduct: then the endofunctor
$\dashbk +1$ on \ess\ has a natural monad structure. This monad is cartesian,
corresponding to the theory of 
\index{pointed sets}%
pointed sets.

\item \label{eg:tree-mnd}
\index{tree!Baez-Dolan-style!monad}%
Let $\ess = \Sets$, and consider the finitary algebraic theory
on $\ess$ generated by one $n$-ary operation for each $n\elt\nat$, and
no equations. This theory is strongly regular, so the induced monad
\Mnd\ on $\ess$ is cartesian.

If $X$ is any set then $X^*$ can be described inductively by:
\begin{itemize}
	\item if $x\elt X$ then $x\elt X^*$
	\item if $\range{t_1}{t_n} \elt X^*$ then
$\abftuple{t_1}{t_n}\elt X^*$.
\end{itemize}
We can draw any element of $X^*$ as a tree with leaves labelled by
elements of $X$:
\begin{itemize}
	\item $x\elt X$ is drawn as \nl{x}
	\item if \range{t_1}{t_n} are drawn as \range{T_1}{T_n} then
\abftuple{t_1}{t_n} is drawn as
$\left.
\begin{tree}
\nodelt{T_1}&	&\nodelt{T_2}&	&	&\cdots	& & &\nodelt{T_n}\\
	&\rt{4}	&	&\rt{2}	&	&	& &\lt{4} &	\\
	&	&	&	&\node	&	& & &		\\
\end{tree}
\right.$
, or if $n=0$, as $\left.\EmptyOne\right.$.
\end{itemize}
Thus the element $\atuplebts{\atuplebts{x_{1},x_{2},\atuplebts{}}, x_{3},
\atuplebts{x_{4},x_{5}}}$ of $X^*$ is drawn as
\[\left.\Oak{x_1}{x_2}{x_3}{x_4}{x_5}\right. .\]

The unit $X\go X^*$ is $x\goesto\nl{x}$, and multiplication $X^{**}
\go X^{*}$ takes an $X^*$-labelled tree (e.g.\ 
\[\left.\Pear{t_1}{t_2}{}\right.,\]
with
\[t_{1}=\left.\Orange{x_1}{x_2}{}\right.
\mr{\ and\ }
t_{2}=\left.\Apple{x_3}{x_4}{}\right. ) \]
and gives an $X$-labelled tree by substituting at the leaves (here,
\[\left.\MixedFruit{x_1}{x_2}{x_3}{x_4}{}\right. ).\]

\item
\index{monoidal category!free}%
On the category \Cat\ of small categories and functors, there is the
free strict 
\index{monoidal category!monad}%
monoidal category monad. Both \Cat\ and the monad are
cartesian.

\item		\label{eg:glob-mnd}
\index{globular set|emph}%
A \emph{globular set} is a diagram
\begin{diagram}
\cdots &\pile{\rTo \\ \rTo} &X_{n+1} &\pile{\rTo^s \\ \rTo_t} &X_n
&\pile{\rTo^s \\ \rTo_t} &\cdots &\pile{\rTo \\ \rTo} &X_1
&\pile{\rTo^s \\ \rTo_t} &X_0 \\
\end{diagram}
in \Sets\ satisfying the `globularity equations' $ss=st$ and
\gobyc{ts=tt}{X_{n+1}}{X_{n-1}}.  The
underlying graph of a strict 
\index{omega-category@$\omega$-category}%
$\omega$-category is a globular set: $X_n$ is
the set of $n$-cells, and $s$ and $t$ are the source and target
functions. One can construct the free strict $\omega$-category monad on the
category of globular sets and show that it is cartesian. Moreover, the
category of globular sets is cartesian, being a presheaf category (i.e.\ of
the form \ftrcat{\scat{G}^{\op}}{\Set}). This
example is explained further in Chapter~\ref{ch:glob}.

\end{eg}

\pagebreak

\section{Multicategories}	\label{sec:multicats}

We now describe what an \Cartpr-multicategory is, where \ust\ is a cartesian
monad on a cartesian category \ess. As mentioned in the introduction to this
chapter, this is a generalization of the (well-known) description of a small
category as a 
\index{monad!in bicategory}%
monad object in the bicategory of spans.

We will use the phrase `\Cartpr\ is 
\index{cartesian}%
cartesian' to mean that \ess\ is a
cartesian category and \Mnd\ is a cartesian monad on \ess.

\begin{constn} \label{constn:bicat} \end{constn}
Let \Cartpr\ be cartesian. We construct a bicategory \cat{B} from
\Cartpr, which in the case $\ust=\id$ is the bicategory of spans in
\ess. Hermida calls \cat{B} the 
\index{Kleisli bicategory of spans}%
\index{spans!Kleisli bicategory of}%
`Kleisli bicategory of spans' in 
\index{Hermida}%
\cite{Her};
the formal similarity between the definition of \cat{B} and the usual
construction of a Kleisli category is evident.

\begin{description}
\item[0-cell:] Object $S$ of \ess.
\item[1-cell $R \go S$:] Diagram
\begin{slopeydiag}
	&	&A	&	&	\\
	&\ldTo	&	&\rdTo	&	\\
R^{*}	&	&	&	&S	\\
\end{slopeydiag}
in \ess.
\item[2-cell $A \go A'$:] Commutative diagram
\begin{slopeydiag}
	&	&A	&	&	\\
	&\ldTo	&	&\rdTo	&	\\
R^{*}	&	&\dTo	&	&S	\\
	&\luTo	&	&\ruTo	&	\\
	&	&A'	&	&	\\
\end{slopeydiag}
in \ess.
\item[1-cell composition:] To define this we need to choose particular
\index{pullback!choice of}%
\index{non-canonical choice}%
pullbacks in \ess, and in everything that follows we assume this has
been done. Take 
\[\left.
\begin{slopeydiag}
	&	&A	&	&	\\
	&\ldTo<{d}&	&\rdTo>{c}&	\\
R^{*}	&	&	&	&S	\\
\end{slopeydiag}
\right.
\mr{\ and\ }
\left.
\begin{slopeydiag}
	&	&B	&	&	\\
	&\ldTo<{q}&	&\rdTo>{p}&	\\
S^{*}	&	&	&	&T	\\
\end{slopeydiag}
\right.
;
\]
then their composite is given by the diagram
\begin{slopeydiag}
   &       &   &       &   &       &B\of A\Spbk&  &   &       &   \\
   &       &   &       &   &\ldTo  &      &\rdTo  &   &       &   \\
   &       &   &       &A^*&       &      &       &B  &       &   \\
   &       &   &\ldTo<{d^*}&&\rdTo>{c^*}& &\ldTo<{q}& &\rdTo>{p}& \\
   &       &R^{**}&    &   &       &S^*   &       &   &       &T  \\
   &\ldTo<{\mu_R}&&    &   &       &      &       &   &       &   \\
R^*&       &   &       &   &       &      &       &   &       &   \\
\end{slopeydiag}
where the right-angle mark in the top square indicates that the square
is a pullback.

\item[1-cell identities:] The identity on $S$ is
\begin{slopeydiag}
	&	&S	&	&	&\\
	&\ldTo<{\eta_S}&&\rdTo>{1}&	&\\
S^{*}	&	&	&	&S	&.\\
\end{slopeydiag}

\item[2-cell identities and compositions:] Identities and vertical
composition are as in \ess. Horizontal composition is given in an
obvious way.

\end{description}

Because the choice of pullbacks is arbitrary, 1-cell composition does
not obey strict associative and identity laws. That it obeys them up
to invertible 2-cells is a consequence of the fact that \Mnd\ is
cartesian. \done

\begin{defn}	\label{defn:multicat}
Let \Cartpr\ be cartesian. Then an 
\index{multicategory|emph}%
\emph{\Cartpr-multicategory} is a monad in the associated bicategory \cat{B} of
Construction~\ref{constn:bicat}.
\end{defn}

An \Cartpr-multicategory therefore consists of a diagram
\spaan{C_1}{C_{0}^*}{C_0}{d}{c} in \ess\ and maps $C_{0} \goby{\ids} C_{1}$,
$C_{1}\of C_{1} 
\goby{\comp} C_{1}$ satisfying associative and identity laws. Think of $C_0$
as `objects', $C_1$ as `arrows', $d$ as `domain' and $c$ as `codomain'.  Such a
multicategory will be called an 
\index{multicategory!on object@\emph{on} object|emph}%
\index{on|emph}%
\Cartpr-multicategory \emph{on $C_0$}, or if
$C_{0}=1$ an
\index{operad|emph}%
\emph{\Cartpr-operad}. 
\index{Baez!-Dolan}%
(Baez and Dolan, in \cite{BD}, use `operad' or 
\index{operad!typed}%
`typed operad' for the same kind of purpose as we use `multicategory', and
`untyped operad' where we use `operad'.)

It is inherent that everything is small: when $\ess=\Sets$, for instance, the
objects and arrows form sets, not classes. For plain multicategories, at
least, there seems to be no practical difficulty in using 
\index{multicategory!large}%
\index{smallness}%
large versions too.

\begin{defn}
Let \Cartpr\ be cartesian.\nopagebreak
\begin{enumerate}
\item \label{defn:graph}
An 
\index{graph|emph}%
\emph{\Cartpr-graph} (on an object $C_0$) is a diagram \gph{C_1}{C_0} in
\ess. A \emph{map of \Cartpr-graphs}
\[\left.
\begin{slopeydiag}
	&	&C_1	&	&	\\
	&\ldTo	&	&\rdTo	&	\\	
C_{0}^{*}	&	&	&	&C_{0}	\\
\end{slopeydiag}
\right.
\go
\left.
\begin{slopeydiag}
		&	&\twid{C_1}&	&	\\	
		&\ldTo	&	&\rdTo	&	\\
\twid{C_0}^{*}	&	&	&	&\twid{C_0}\\
\end{slopeydiag}
\right.\]
is a pair \pr{C_{0} \goby{f_0} \twid{C_0}}{C_{1}\goby{f_1} \twid{C_1}} of maps
in
\ess\ such that
\begin{slopeydiag}
	&	&C_1	&	&	\\
	&\ldTo	&\dTo>{f_1}&\rdTo	&	\\
C_{0}^{*}	&	&	&	&C_0	\\
\dTo<{f_{0}^*}&	&\twid{C_1}&	&\dTo>{f_0}\\
	&\ldTo	&	&\rdTo	&	\\
\twid{C_0}^{*}&	&	&	&\twid{C_0}\\
\end{slopeydiag}
commutes.

\item A 
\index{multicategory!map of}%
\emph{map of \Cartpr-multicategories} $C \go \twid{C}$ (with graphs
as in~(\ref{defn:graph})) is a map $f$ of their graphs such that the
diagrams 
\[\left.
\begin{diagram}
C_0		&\rTo^{f_0}		&\twid{C_0}		\\	
\dTo<{\ids}	&			&\dTo>{\twid{\ids}}	\\
C_1		&\rTo^{f_1}		&\twid{C_1}		\\
\end{diagram}
\right.
\ \ \ 
\left.
\begin{diagram}
&C_{1}\of C_1	&\rTo^{f_{1}*f_1}	&\twid{C_1}\of\twid{C_1}	\\
&\dTo<{\comp}	&			&\dTo>{\twid{\comp}}	\\
&C_1		&\rTo^{f_1}		&\twid{C_1}		\\	
\end{diagram}
\right.\]
commute.

\end{enumerate}
\end{defn}

\begin{remarks}{rmks:maps}
\item Fix $S\elt\ess$. Then we may consider the category of
\Cartpr-graphs on $S$, whose morphisms $f = \pr{S \goby{f_0} S}
{C_{1} \goby{f_1} \twid{C_1}}$ all have $f_{0}=1$. This is just the slice
category $\frac{\ess}{S^{*}\times S}$. It is also the full sub-bicategory of
\cat{B} whose only object is $S$, and is therefore a monoidal
category. The category of \Cartpr-multicategories on $S$ is then the
category $\Mon(\frac{\ess}{S^{*}\times S})$ of monoids in 
$\frac{\ess}{S^{*}\times S}$.

\item A 
\index{pullback!choice of}%
\index{non-canonical choice}%
choice of pullbacks in \ess\ was made; changing that choice
gives an isomorphic category of \Cartpr-multicategories. 

\item
If \pr{\ess'}{\ush} is also cartesian then a 
\index{monad!(op)functor}%
monad functor
$\Cartpr\go\pr{\ess'}{\ush}$ gives a functor
$\Multicatof{\Cartpr}\go\Multicatof{\pr{\ess'}{\ush}}$, and the same is true
of monad opfunctors. See 
\index{Street}%
\cite{StFTM} for the terminology and \cite[4.4]{GOM}
for more details.
\end{remarks}

\begin{eg}{eg:multicats}

\item
Let $\Cartpr=\pr{\Sets}{\id}$. Then \cat{B} is the
bicategory of 
\index{spans!bicategory of}%
spans, and a monad in \cat{B} is just a (small)
category. Thus categories are \pr{\Sets}{\id}-multicategories.
Functors are maps of such. More generally, if \ess\ is any cartesian category
then \pr{\ess}{\id}-multicategories are 
\index{internal category}%
internal categories in \ess.

\item
\index{monoid!free}%
Let $\Cartpr=\pr{\Sets}{\mr{free\ monoid}}$. Specifying an \Cartpr-graph
\graph{C_1}{C_0}{d}{c} is equivalent to specifying a set
\mtihom{C}{s_1}{s_n}{s} for each \range{s_1}{s_n,s}\elt$C_0$ ($n\geq
0$); if $a\elt\mtihom{C}{s_1}{s_n}{s}$ then we write
\[
\range{s_1}{s_n} \goby{a} s
\]
or
\[
\ctransistor{a}{s_1}{s_2}{s_n}{s}
\]
or
\[
\begin{opetope}
	&	&	&\cnr	&\ldots	&	&	&	\\
	&\cnr	&\ruEdge(2,1)<{s_2}&&	&	&\cnr	&	\\
\ruEdge(1,2)<{s_1}&&	&	&\Downarrow a&	&	&\rdEdge(1,2)>{s_n}\\
\cnr	&	&	&\rEdge_{s}&	&	&	&\cnr	\\
\end{opetope}
\ \ \ \ .
\]

In the associated bicategory, the identity 1-cell
\graph{C_0}{C_0}{\eta_{C_0}}{1} on $C_0$ has
\[
\mtihom{C_0}{s_1}{s_n}{s}=
\left\{
\begin{array}{ll}
1		&\mr{if\ }n=1\mr{\ and\ }s_{1}=s	\\
\emptyset	&\mr{otherwise.}
\end{array}
\right.
\]
The composite 1-cell $C_{1}\of C_1$ is
\[
\{\pr{\bftuple{a_1}{a_n}}{a} \such da=\bftuple{ca_1}{ca_n}\},
\]
i.e.\ is the set of diagrams
\begin{equation} \label{diag:comptrans}
\comptrans{a_1}{a_2}{a_n}{a}
\end{equation}
with the evident domain and codomain functions.

We then have a function \ids\ assigning to each $s\elt C_0$ a member of
\mtiendo{C}{s}, and a function \comp\ composing diagrams
like~(\ref{diag:comptrans}). These are required to obey associative and
identity laws. Thus a \pr{\Sets}{\mr{free\ monoid}}-multicategory is just a
\index{multicategory!plain}%
plain multicategory and a \pr{\Sets}{\mr{free\ monoid}}-operad is a 
\index{operad!plain}%
plain operad.

\item
Let $\Cartpr=\pr{\Set}{\dashbk +1}$. It is not hard to see
(\cite[2.6(iv)]{GOM}) that an \Cartpr-multicategory is a (small) category
$C$ together with a functor $C\go\Set$. To put it another way, an
\Cartpr-multicategory is a 
\index{opfibration!discrete}%
discrete opfibration (between small categories);
in fact, the category of \Cartpr-multicategories is the category of discrete
opfibrations.

\item		\label{eg:tree-mti}
\index{tree!Baez-Dolan-style!monad}%
Let $\Cartpr = \pr{\Sets}{\mr{tree\ monad}}$, as
in~\ref{eg:cart-mnds}(\ref{eg:tree-mnd}). An \Cartpr-multicategory
consists of a set $C_0$ of objects, and sets like
\[
C\left(\Pear{s_1}{s_2}{s}\right)
\]
($s_{1}, s_{2}, s \elt C_0$), together with a unit element of each
$C(\nlal{s}{s})$ and composition 
\pagebreak
functions like
\begin{eqnarray*}
\left\{	C\left( \Orange{r_1}{r_2}{s_1}\right)
	\times
	C\left( \Apple{r_3}{r_4}{s_2}\right)
\right\}
\times
C\left( \Pear{s_1}{s_2}{s}\right)\\
\go
C\left( \MixedFruit{r_1}{r_2}{r_3}{r_4}{s}\right)
\end{eqnarray*}
($r_{1}, r_{2}, r_{3}, r_{4} \elt C_0$). These are to satisfy associativity
and identity laws.

When $C_{0}=1$, so that we're considering \Cartpr-operads, the graph
structure is comprised of sets like $C\left(\Pear{}{}{}\right)$.

The \Cartpr-multicategories are a simpler version of 
\index{Soibelman}%
Soibelman's 
\index{pseudo-monoidal category}%
pseudo-monoidal categories (\cite{Soi}); they omit the aspect of
\index{tree!Baez-Dolan-style!map of}%
maps between trees. A similar relation is borne to
\index{Borcherds}%
Borcherds' 
\index{relaxed multilinear category}%
relaxed multilinear categories (\cite{Bor}, 
\index{Snydal}%
\cite{Sny}).

\item		\label{eg:soib-op}
When $\ess=\fcat{Cat}$ and \ust\ is the free strict 
\index{monoidal category!monad}%
monoidal category monad, an \Cartpr-operad is what 
\index{Soibelman}%
Soibelman calls a 
\index{strict!monoidal 2-operad}%
strict monoidal 2-operad in
\cite[2.1]{Soi}. Such a structure might also be thought of as a plain operad
\index{enrichment!of operad in Cat@of operad in \Cat}%
enriched in \Cat.

\item		\label{eg:glob-ops}
\index{globular set}%
\index{omega-category@$\omega$-category}%
\index{operad!Batanin}%
Let $\Cartpr=\pr{\mr{Globular\ sets}}{\mr{free\ strict\
}\omega\hyph\mr{category}}$. An \Cartpr-operad is exactly what Batanin
calls an operad (see Chapter~\ref{ch:glob}).

\end{eg}

\section{Algebras}	\label{sec:algs}
\index{algebra!for multicategory or operad|(}

Most of the existing notions of operad (e.g.\ 
\index{May}%
\cite{May}, 
\index{Baez!-Dolan}%
\cite{BD},
\index{operad!Batanin}%
\cite{Bat}) carry with them the idea of an \emph{algebra} for an operad. We
now define a category of algebras for any \Cartpr-multicategory.

In the case $\Cartpr=\pr{\Set}{\id}$, where we are dealing with a normal
category $C$, the category of algebras will be
\ftrcat{C}{\Sets}. Now \ftrcat{C}{\Sets} is equivalent to the
category of 
\index{opfibration!discrete|(}%
discrete opfibrations over $C$, where by definition, a
discrete opfibration over $C$ is a functor $D\goby{f}C$
such that for each arrow in $C$, every lift to $D$ of the source of
the arrow extends uniquely to a lift of the whole arrow. This means precisely
that the left-hand half of the diagram
\begin{slopeydiag}
	&		&D_1		&		&	\\
	&\ldTo<{d}	&\dTo>{f_1}	&\rdTo>{c}	&	\\
D_0	&		&		&		&D_0	\\
\dTo<{f_0}&		&C_1		&		&\dTo>{f_0}\\
	&\ldTo<{d}	&		&\rdTo>{c}	&	\\
C_0	&		&		&		&C_0	\\
\end{slopeydiag}
is a pullback square. So by analogy, we will say that a map
$D\goby{f}C$ of \Cartpr-multicategories is a 
\index{opfibration!discrete!of multicategories|emph}%
\emph{discrete opfibration} if
\begin{diagram}
D_{0}^{*}	&\lTo^{d}	&D_1		\\
\dTo<{f_{0}^{*}}&		&\dTo>{f_1}	\\
C_{0}^{*}	&\lTo^{d}	&C_1		\\
\end{diagram}
is a pullback square, and define an 
\index{algebra!for multicategory or operad|emph}%
\emph{algebra} for $C$ to be a discrete
opfibration over $C$. A \emph{map} $(D\goby{f}C)$ to
$(\twid{D} \goby{\twid{f}} C)$ of algebras for $C$
consists of a map $D\goby{g}\twid{D}$ of multicategories such
that $\twid{f}g = f$. We thereby obtain the category of algebras $\Alg(C)$
for any \Cartpr-multicategory $C$. 

There is an alternative definition of algebra, which gives a category of
algebras equivalent to the original one. It is longer to state but seems to
be easier to use in practice. The initial observation is that if $C$ is a
normal category then the forgetful functor $\ftrcat{C}{\Set} \go \Set^{C_0}$
is monadic, the monad $T$ on $\Set^{C_0}$ being given by
\[
(TX)s=\coprod_{s' \goby{a} s \mr{\ in\ } C} Xs'
\]
($X\elt\Set^{C_0}$, $s\elt C_{0}$). Transferring $T$ across the
equivalence $\Set^{C_0} \eqv \Set/C_{0}$, we get a monad $T'$ on $\Set/C_{0}$
such that $\ftrcat{C}{\Set}$ is equivalent to the category of $T'$-algebras.
It turns out (\cite[3.2]{GOM}) that if \slob{X}{p}{C_0} is an object of
$\Set/C_{0}$ then $T'\bktdslob{X}{p}{C_0}$ is the right-hand diagonal of the
diagram
\[\left.
\begin{slopeydiag}
	&	&\cdot\Spbk&	&	&	&	\\
	&\ldTo	&	&\rdTo	&	&	&	\\
X	&	&	&	&C_{1}&	&	\\
	&\rdTo<{p}&	&\ldTo>{d}&	&\rdTo>{c}&	\\
	&	&C_{0}	&	&	&	&C_{0}\\
\end{slopeydiag}
\right.
.
\]
We are now ready to generalize to an arbitrary cartesian \Cartpr.

\begin{constn} \label{constn:ind-mnd} \end{constn}
Let \Cartpr\ be cartesian: then any \Cartpr-multicategory $C$
gives rise to a monad on $\ess/C_0$. The alternative definition of the
category of algebras for $C$ is as the category of algebras for
this monad.

The functor part of the monad on \Slice will be called \blbk; in what
follows, we'll write $\bktdslob{X}{p}{C_0}^{\blob} =
\slob{X_\blob}{p_\blob}{C_0}$. The details omitted below are given in
\cite[3.3]{GOM}.

\begin{itemize}
\item \slob{X_\blob}{p_\blob}{C_0} is the composite down the right-hand
diagonal of
\[\left.
\begin{slopeydiag}
	&	&\cdot\Spbk&	&	&	&	\\
	&\ldTo	&	&\rdTo	&	&	&	\\
X^*	&	&	&	&C_1	&	&	\\
	&\rdTo<{p^*}&	&\ldTo>{d}&	&\rdTo>{c}&	\\
	&	&C_{0}^*&	&	&	&C_0	\\
\end{slopeydiag}
\right.
.
\]
The functor \blbk\ is defined on morphisms in an obvious way.

\item The unit at \slob{X}{p}{C_0} is given by
\[\left.
\begin{slopeydiag}
	&			&X		&	&	\\
	&\ldTo(2,5)<{\eta_X}	&		&\rdTo>{p}&	\\
	&			&\dGet~{\unit_p}	&	&C_0	\\
	&			&X_{\blob}\Spbk	&	&\dTo>{\ids}\\
	&\ldTo			&		&\rdTo	&	\\
X^{*}	&			&		&	&C_1	\\
	&\rdTo<{p^*}		&		&\ldTo>{d}&	\\
	&			&C_{0}^{*}		&	&	\\
\end{slopeydiag}
\right.
\ \ \ \ \ \ .
\]
The multiplication is defined in a similar, if slightly more complicated,
fashion.
\end{itemize}

It is now straightforward to check that \Imnd\ forms a monad on \Slice, and
that the category of algebras for \Imnd\ is equivalent to the category of
discrete opfibrations over $C$. We indicate roughly how this isomorphism
works. Given an algebra $\bktdslob{X}{p}{C_0}^{\blob} \goby{h}
\bktdslob{X}{p}{C_0}$ for \blbk, there is a diagram
\begin{slopeydiag}
	&	&X_\blob	&	&	\\
	&\ldTo	&\dTo		&\rdTo>{h}&	\\
X^{*}	&	&		&	&X	\\
\dTo<{p^*}&	&C_1		&	&\dTo>{p}\\
	&\ldTo<{d}&		&\rdTo>{c}&	\\
C_{0}^{*}&	&		&	&C_0\\
\end{slopeydiag}
the left-hand half of which is a pullback square,
and a multicategory structure on \gph{X^\blob}{X} such that this diagram
defines a map of multicategories---and therefore a discrete opfibration. 

In the case $\Cartpr=\pr{\Set}{\id}$, where \gph{C_1}{C_0} is a normal
category \spn{C_1}{C_0}{C_0}, the object \bktdslob{X}{p}{C_0} of $\Set/C_0$
gives a map $C_{0} \go \Set$ by $x \goesto p^{-1}\{x\}$, and the algebra
structure $h$ extends this to a functor $C\go\Set$. The category
\spn{X_\blob}{X}{X} is then the
\index{opfibration!Grothendieck}%
Grothendieck opfibration of the functor.

Finally, with the \pr{\Sets}{\id} case of plain categories in mind, we would
expect a map $C \go C'$ of multicategories to 
\index{Alg@\Alg!functoriality of}%
yield a functor $\Alg(C) \og
\Alg(C')$.  This is indeed the case, as may easily be verified using either
of the two definitions.%
\index{opfibration!discrete|)}%

\begin{eg}{eg:algs}

\item 
When $\Cartpr=\pr{\Sets}{\id}$, $\Alg(C)$ is equivalent to
$\ftrcat{C}{\Sets}$ (or isomorphic, if we use the second definition of \Alg).

\item
\index{monoid!free}%
When $\Cartpr=\pr{\Sets}{\mr{free\ monoid}}$, so that an
\Cartpr-multicategory is a multicategory of the familiar kind, we already
have an idea of what an algebra for $C$ should be: a 
\index{multifunctor}%
`multifunctor $C\go\Sets$'. That is, an algebra for $C$ should consist of:
	\begin{itemize}
	\item for each $s\elt C_0$, a set $X(s)$
	\item for each $\range{s_1}{s_n}\goby{a}s$ in $C$, a function
	$X(s_{1})\times\cdots\times X(s_{n}) \go X(s)$, in a way compatible
	with identities and composition.  
	\end{itemize}
In fact, this is the same as either definition of algebra given above. We
work with the second one: algebras for $C$ are algebras for the induced monad
\ubl\ on \Slice. If \slob{X}{p}{C_0} is an object of \Slice\ and we put
$X(s)=p^{-1}\{s\}$ then
\begin{eqnarray*}
X_{\blob}	&=&	\{\pr{\bftuple{x_1}{x_n}}{f} \such
			df=\bftuple{px_1}{px_n}\}\\
		&=&	\{X(s_{1})\times\cdots\times X(s_{n}) \times
			\mtihom{C}{s_1}{s_n}{s} \such
			\range{s_1}{s_{n},s}\elt C_0\},
\end{eqnarray*}
and an algebra structure on \slob{X}{p}{C_0} therefore consists of a function
\[
X(s_{1})\times\cdots\times X(s_{n}) \go X(s)
\]
for each member of \mtihom{C}{s_1}{s_n}{s}, subject to certain laws.

\item
\index{tree!Baez-Dolan-style!monad}%
Let \Cartpr\ be the tree monad on \Sets; for simplicity, let us just consider
\emph{operads} $C$ for \Cartpr---thus the object-set $C_0$ is 1. An algebra for
$C$ consists of a set $X$ together with a function $X_{\blob}\go X$
satisfying some axioms. One can calculate that an element of $X_{\blob}$
consists of an $X$-labelling of the leaves of a tree $T$, together with a
member of $C(T)$.  An $X$-labelling of an $n$-leafed tree $T$ is just a
member of $X^n$, so one can view the algebra structure $X_{\blob}\go X$ on
$X$ as: for each number $n$, $n$-leafed tree $T$, and element of $C(T)$, a
function $X^{n}\go{X}$.  These functions are required to be compatible with
gluing of trees in an evident way.

\item			\label{eg:glob-algs}
\index{globular set}%
\index{omega-category@$\omega$-category}%
For $\Cartpr=\pr{\mr{Globular\ sets}}{\mr{free\ strict\ }
\omega\hyph\mr{category}}$, 
\index{Batanin}%
Batanin constructs a certain operad $K$, the
`universal contractible operad'. He then defines a
\index{omega-category@$\omega$-category!lax}%
lax $\omega$-category to be an algebra for $K$. 
\index{Batanin!misrepresentation}%
We follow exactly this strategy for defining lax $\omega$-category in
Chapter~\ref{ch:glob}, except that our operad $K$ is different from his.

\item	\label{eg:terminal}
\index{terminal multicategory}%
\index{multicategory!terminal}%
\index{operad!terminal}%
The graph \graph{1^*}{1}{1}{!} is terminal amongst all \Cartpr-graphs. It
carries a unique multicategory structure, since a terminal object in a
monoidal category always carries a unique monoid structure. It then becomes
the terminal \Cartpr-multicategory. The induced monad on $\ess/1$ is just
\Mnd, and so an algebra for the terminal multicategory is just an algebra for
\ust. (E.g.\ an algebra for the terminal plain operad is a monoid.)
This observation can aid recognition of when a theory of operads or
multicategories fits into our scheme. For instance, if we were to read
\index{Batanin}%
Batanin's paper and learn that, in his terminology, an algebra for the
terminal operad is a strict $\omega$-category (\cite[p.\ 51, example
3]{Bat}), then we might suspect that his operads were \Cartpr-operads for the
free strict $\omega$-category monad \ust\ on some suitable category \ess---as
indeed they are.

\end{eg}

\index{algebra!for multicategory or operad|)}

\section{Structured Categories}	\label{sec:struc}
\index{structured category|(}%

The observation from which this section takes off is that any strict 
\index{monoidal category!and multicategories}%
monoidal category has an underlying multicategory. (All monoidal categories
and maps between them will be strict in this section. For the time being,
`multicategory' means plain multicategory.)  Explicitly, if
\pr{\scat{C}}{\otimes} is a monoidal category, then the underlying
multicategory $C$ has the same object-set as \scat{C} and has homsets defined
by
\[
\mtihom{\Hom_{C}}{s_1}{s_n}{s} = 
\homset{\Hom_{\scat{C}}}{s_{1}\otimes\cdots\otimes s_{n}}{s}
\]
for objects \range{s_1}{s_n,s}. Composition and identities in $C$ are easily
defined.

There is a converse process: given any multicategory $C$,
there is a 
\index{monoidal category!free}%
`free' monoidal category \scat{C} on it.  Informally, an
object/arrow of \scat{C} is a sequence of objects/arrows of $C$.  Thus the
objects of \scat{C} are of the form \bftuple{s_1}{s_n} ($s_i \elt C_0$), and a
typical arrow $\tuplebts{s_1,s_2,s_3,s_4,s_5} \go
\tuplebts{s'_1,s'_2,s'_3}$ is a sequence \tuplebts{a_1,a_2,a_3} of elements
of $A$ with domains and codomains as illustrated:
\begin{equation}	\label{diag:arrows-in-mon-cat} 
\freemoncatpic
\end{equation}
The tensor in \scat{C} is just juxtaposition.

For example, the terminal multicategory \fcat{1} has one object and, for each
$n\elt\nat$, one arrow of the form
\[
n \left\{\rule{0ex}{5.3ex}\right.\left.
\ctransistor{}{}{}{}{}
\right.
\ ;
\]
figure~\ref{diag:arrows-in-mon-cat} (above) indicates that the `free'
monoidal category on the multicategory \fcat{1} is 
\index{simplicial category}%
$\Delta$, the category of finite ordinals (including 0), with addition as
$\otimes$.

The name `free' is justified: that is, there is an adjunction
\begin{diagram}
\fcat{Monoidal\ Categories}	\\
\uTo\dashv\dTo			\\
\fcat{Multicategories}		\\
\end{diagram}
where the two functors are those described above. Moreover, this adjunction
is monadic. (But the forgetful functor does \emph{not} provide a full
embedding of \fcat{Monoidal\ Categories} into \fcat{Multicategories}. It is
faithful, but not full: there is a multicategory map $\fcat{1} \go \Delta$
sending the unique object of \fcat{1} to the object 1 of $\Delta$, and this
map does not preserve the monoidal structure.)

Naturally, we would like to generalize from $\Cartpr=\pr{\Sets}{\mr{free\
monoid}}$ to any cartesian \Cartpr. To do this, we need a notion of
`\Cartpr-structured category', which in the case \pr{\Sets}{\mr{free\
monoid}} just means monoidal category. A monoidal category is a category
object in \fcat{Monoid}, so it is reasonable to define an
\index{structured category|emph}%
\emph{\Cartpr-structured category} to be an
\pr{\ess^\stbk}{\id}-multicategory---that is, a category object in the
category $\ess^\stbk$ of algebras for the monad \stbk\ on \ess.

It is now possible to describe a monadic adjunction
\index{structured category!free}%
\begin{diagram}
\Strucof{\Cartpr}			\\
\uTo<{F}\dashv\dTo>{U}			\\
\Multicatof{\Cartpr}			\\
\end{diagram}
generalizing that above. The effect of the functors $U$ and $F$ on objects is
as outlined now. Given an \Cartpr-structured category $D$, with algebraic
structure $D_{0}^{*} \goby{\otimes} D_{0}$ and $D_{1}^{*} \goby{\otimes}
D_{1}$, the graph \gph{C_1}{D_0} of $UD$ is given by
\[
\begin{slopeydiag}
	&		&C_{1}\Spbk&	&	&	&	\\
	&\ldTo		&	&\rdTo	&	&	&	\\
D_{0}^{*}&		&	&	&D_{1}	&	&	\\
	&\rdTo<{\otimes}&	&\ldTo	&	&\rdTo	&	\\
	&		&D_0	&	&	&	&D_0	\\
\end{slopeydiag}
\ .
\]
Given an \Cartpr-multicategory $C$, the category $FC$ has
graph
\begin{slopeydiag}
	&	&	&	&C_{1}^{*}&	&	&	&	\\
	&	&	&\ldTo	&	&\rdTo(4,4)&	&	&	\\
	&	&C_{0}^{**}&	&	&	&	&	&	\\
	&\ldTo<{\mu_{C_0}}&&	&	&	&	&	&	\\
C_{0}^*	&	&	&	&	&	&	&	&C_{0}^*\\
\end{slopeydiag}
and the algebraic structures $C_{1}^{**} \goby{\otimes} C_{1}^{*}$, $C_{0}^{**}
\goby{\otimes} C_{0}^{*}$ are components of $\mu$.

\index{structured category|)}%

\section{The Free Multicategory Monad}	\label{sec:free-multi}
\index{multicategory!free|(}%
\index{operad!free|(}%

Let \Cartpr\ be cartesian. Subject to certain further conditions on \ess\ and
\ust, the following are true:
\begin{itemize}
	\item the forgetful functor
		\[
		\Multicatof{\Cartpr} \go \Graphof{\Cartpr}
		\]
	has a left adjoint, the `free \Cartpr-multicategory functor'
	\item the adjunction is monadic
	\item the monad on \Graphof{\Cartpr} is cartesian (and the category
	\Graphof{\Cartpr} is cartesian, as always)
\end{itemize}

The precise nature of the conditions and the free multicategory
construction is not important for our purposes. However, as a gesture
towards supporting the assertions above, here is a sketch of a construction
of the free functor. Given an \Cartpr-graph $G$, one defines a sequence
$(\gph{A^{(n)}}{G_0})_{n\elt\nat}$ of \Cartpr-graphs by $A^{(0)}=G_1$ and
$A^{(n+1)}=G_{0}+G_{1}\of A^{(n)}$, where $+$ is binary coproduct and \of\ is
composition in the bicategory of spans for \Cartpr. The free multicategory
\gph{A}{G_0} on $G$ is given as a colimit of the $A^{(n)}$'s. (In the case
$\Cartpr=\pr{\Set}{\id}$ and $G_{0}=1$ we are forming the 
\index{monoid!free}%
free monoid on a set $G$, and the formula $A^{(n+1)}=1+G\of A^{(n)}$
expresses the fact that a word on a set is either the empty word or an
element of the set followed by a word.)  To make this work it is necessary to
assume that in \ess\ certain colimits exist and interact with pullbacks in a
suitable way, and that the functor \stbk\ preserves certain
\index{filtered colimit}%
filtered colimits.

As may be apparent from the foregoing sketch of the construction, taking the
free multicategory on a graph does not change the objects-object $S$
($=G_0$). It follows that all the statements in the first paragraph of this
section hold if we replace the forgetful functor
\[
\Multicatof{\Cartpr} \go \Graphof{\Cartpr}
\]
by the forgetful functor 
\[
\Cartpr\mbox{-Multicategories on }S \go \Cartpr\mbox{-Graphs on }S,
\]
for any $S\elt\ess$. We therefore have the monad `free \Cartpr-multicategory
on $S$' on the cartesian category $\frac{\ess}{S^{*}\times S} =
\Cartpr\mbox{-Graphs on }S$. 

All we will need to know about the conditions on \ess\ and \ust\ is:
\begin{itemize}
\item the conditions on \ess\ are satisfied if \ess\ is \Set, or any presheaf
category
\item the conditions on \ust\ are satisfied if \ust\ is the identity monad,
or in fact just as long as the functor \stbk\ is 
\index{finitary!functor}%
finitary
\item if \Cartpr\ satisfies the conditions then so too does
\[
\pr{\Cartpr\mbox{-Graphs on }S}{\mbox{free }\Cartpr\mbox{-multicategory on }S},
\]
for each $S\elt\ess$.
\end{itemize}
For example, these conditions hold for any 
\index{finitary!algebraic theory}%
cartesian finitary algebraic theory on
\Set.

\index{operad|)}
\index{multicategory|)}
\index{multicategory!free|)}
\index{operad!free|)}

%% file: essglob.tex
\chapter{The Globular Approach}		\label{ch:glob}
\index{n-category@$n$-category!lax|(}
\index{Batanin|(}

Perhaps the most direct approach to defining `lax $n$-category' is the kind
suggested by Batanin in \cite{Bat}. In this chapter we present a definition
of lax $n$-category which is a variation on Batanin's. An informal version
can be given using no technical terms at all (\ref{sec:infout}).  The
definition also has a very simple structure when described in the language of
general operads, as already seen in \ref{eg:cart-mnds}(\ref{eg:glob-mnd}),
\ref{eg:multicats}(\ref{eg:glob-ops}) and \ref{eg:algs}(\ref{eg:glob-algs}).

\index{Batanin!misrepresentation}%
Our definition of lax $n$-category is very close to Batanin's, although not
(as asserted previously) the same. The operads he uses \emph{are} the same as
the \Cartpr-operads here (for a certain choice of \Cartpr), and part of the
purpose of this chapter is to explain in elementary language and pictures
what these \Cartpr-operads are, so that the knowledgeable reader may
understand that the two kinds of operad coincide. However, the other main
concept used in the definition, contractibility, is defined differently in
our two accounts.

The final section of the chapter is on another theme. It is a brief
explanation of how we might go about defining (strict and lax) 
\index{cubical!category}%
$n$-tuple categories, which are to $n$-categories as 
\index{double category}%
double categories are to 2-categories. In other words, it is a cubical rather
than a globular approach, and the development is exactly
analogous. Conceivably there is some connection here with
\index{nerve}%
\index{Tamsamani}%
\label{p:tam-rmk}%
Tamsamani's $n$-categories. A category is a finite-limit-preserving functor
$\Delta^{\op}\go\Set$. In \cite{Tam}, a lax $n$-category is a functor
$(\Delta^{n})^{\op}\go\Set$ satisfying three axioms. The first is a condition
called `truncatability'. The second is a generalization of the
`finite-limit-preserving' condition. The third is a 
\index{degeneracy!of cubical set}%
degeneracy condition, analogous to the fact that a 2-category is a double
category with a certain degeneracy in its graph. The present author's
ignorance prevents further discussion.

Section~\ref{sec:infout} is an informal definition of lax
$n$-category;~\ref{sec:formout} describes the formal structure of the
definition. Section~\ref{sec:trees} concerns globular diagrams like
\[
\gfst{}\gtwo{}{}{}\gblw{}\gone{}\gblw{}\gfour{}{}{}{}{}{}{}\glst{}
\]
and how they may be manipulated: these are the cells of the free strict
$\omega$-category $\wun^*$ on the terminal globular set. One way to define
$X^*$ for general $X$ is to define $\wun^*$ and then each fibre in the map
$X^{*}\goby{!^*}\wun^*$, and we discuss this here. In
\ref{sec:ops-and-algs}, operads and their algebras are explained in
elementary terms. We then reach the definition of lax $n$- and
$\omega$-categories in~\ref{sec:con-and-n-cat}, via the concept of a
contraction on an operad (which is not part of the general theory of
operads). Finally, section~\ref{sec:cub-ap} is on cubical structures, as
explained above.

The chief reference for this chapter is \cite{Bat}, for the material on
operads and a different view of contractibility.
\index{Street}%
Street has also produced an account of Batanin's work, \cite{StRMB}.

\section{Informal Outline}	\label{sec:infout}

As in the traditional conception, the graph structure of a lax $n$-category
consists of 0-cells \gzero{A}, 1-cells \gfst{A}\gone{f}\glst{B}, 2-cells
\gfst{A}\gtwo{f}{g}{\alpha}\glst{B}, \ldots. There are then various ways of
\index{composition!ways of composing|(}%
composing these cells; just how many ways and how they interact depends on
whether we are dealing with strict or lax $n$-categories, or something in
between. In a 
\index{n-category@$n$-category!strict}%
strict $n$-category, there will be precisely one way of
composing a diagram like
\begin{equation}
\gfst{A}%
\gthree{f}{f'}{f''}{\alpha}{\alpha'}%
\grgt{B}%
\gone{g}%
\glft{C}%
\gtwo{h}{h'}{\gamma}%
\glst{D}
\label{pic:typical-labelled-2-glob}
\end{equation}
to obtain a 2-cell: that is, any two different ways of doing it (e.g.\
compose $\alpha'$ with $\alpha$, and $\gamma$ with $g$, then the two of these
together) give exactly the same resulting 2-cell. In a lax $n$-category
there will be many ways, but the resulting 2-cells will all be equivalent in a
suitably weak sense.

Our method of describing what ways of composing are available in a lax
$n$-category depends on one simple principle, the 
\index{contraction!principle|(}%
`contraction principle'.
Take, for example, the diagram~(\ref{pic:typical-labelled-2-glob})
above. 
\label{p:star}
Suppose we have already constructed two ways of composing a generic diagram
\[
\gfst{}%
\gone{p}%
\gblw{}%
\gone{q}%
\gblw{}%
\gone{r}%
\glst{}
\]
of 1-cells, namely $(rq)p$ and $r(qp)$. Then we deduce that there is a way of
composing diagram~(\ref{pic:typical-labelled-2-glob}) to get a 2-cell of
the form \gfst{A}\gtwo{(hg)f}{h'(gf'')}{}\glst{D}. As another example of the
principle, this time in one higher dimension, take a diagram
\begin{equation}	\label{pic:typical-labelled-3-glob}
\gfst{A}%
\gspecialone{f}{f'}{f''}{\alpha}{\alpha'}{\alpha''}{\beta}{x}{y}%
\grgt{B}%
\gonew{g}%
\glft{C}%
\gthreecell{h}{h'}{\gamma}{\gamma'}{z}%
\glst{D}.
\end{equation}
Suppose we have constructed two ways of composing a diagram of the shape
of~(\ref{pic:typical-labelled-2-glob}) to a 2-cell, each of which invokes the
same way of composing the 1-cells
\[
\gfst{}%
\gone{}%
\gblw{}%
\gone{}%
\gblw{}%
\gone{}%
\glst{}
\]
along the top and bottom. Say, for instance, that the first way of
composing~(\ref{pic:typical-labelled-2-glob}) results in a 2-cell
\gfst{A}\gtwo{(hg)f}{h'(gf'')}{\delta}\glst{D} and the second way in a 2-cell
\gfst{A}\gtwo{(hg)f}{h'(gf'')}{\delta'}\glst{D}. Then the contraction
principle says that there is a way of
composing~(\ref{pic:typical-labelled-3-glob}) to get a 3-cell of the form
\[
\gfst{A}\gthreecell{(hg)f}{h'(gf'')}{\delta}{\delta'}{}\glst{D}.
\]

In general, we may state the contraction principle as follows. Suppose we are
given an $(n+1)$-dimensional diagram and two ways of composing the
$n$-dimensional diagram at its source/target, such that these two ways match
on the $(n-1)$-dimensional source and target. Then there's a way of composing
the $(n+1)$-dimensional diagram, inducing the first way on its source and the
second way on its target. (In our first example, we implicitly used the
fact that the two ways of composing
\[
\gfst{}%
\gone{p}%
\gblw{}%
\gone{q}%
\gblw{}%
\gone{r}%
\glst{},
\]
$(rq)p$ and $r(qp)$, do the
same thing to the bounding 0-cells: nothing at all.)

Let us see how the contraction principle generates the ways of composing in a
lax $\omega$-category. First of all, there's a way of composing a
(`generic') 0-cell \gzero{A} to get a 0-cell: do nothing. Now, given a diagram
\[
\gfst{A_0}%
\gone{f_1}%
\gblw{A_1}%
\gone{f_2}%
\diagspace%
\cdots%
\diagspace%
\gone{f_n}%
\glst{A_n}
\]
of 1-cells, note that there are two (identical) ways of composing the 0-cell
diagram \gzero{}, both of which are the `do nothing' way just referred to. So
by the contraction principle, there arises a way of composing
\[
\gfst{A_0}%
\gone{f_1}%
\gblw{A_1}%
\gone{f_2}%
\diagspace%
\cdots%
\diagspace%
\gone{f_n}%
\glst{A_n}
\]
to a 1-cell \gfst{A_0}\gone{}\glst{A_n}; we call this 1-cell
\gfst{A_0}\gone{f_{n}\ldots f_{2}f_{1}}\glst{A_n}. Putting together these
ways of composing 1-cells, we can take a diagram like
\[
\gfst{}%
\gone{f_1}%
\gblw{}%
\gone{f_2}%
\gblw{}%
\gone{f_3}%
\gblw{}%
\gone{f_4}%
\gblw{}%
\gone{f_5}%
\glst{}
\]
and do a composition like
\begin{eqnarray*}
&
\gfst{}%
\gone{f_1}%
\gblw{}%
\gone{f_2}%
\gblw{}%
\gone{f_3}%
\gblw{}%
\gone{f_4}%
\gblw{}%
\gone{f_5}%
\glst{}\\
\goesto	&
\gfst{}%
\gone{f_{3}f_{2}f_{1}}%
\gblw{}%
\gone{1}%
\gblw{}%
\gone{f_{5}f_{4}}%
\glst{}\\
\goesto	&
\gfst{}%
\gone{(f_{5}f_{4})1(f_{3}f_{2}f_{1})}%
\glst{}
\end{eqnarray*}
where to get `1' we've taken $n=0$ in the above.
Next, we may take a 2-dimensional diagram
like~(\ref{pic:typical-labelled-2-glob}) and find a way of composing it by
using its 1-dimensional boundary, as on page~\pageref{p:star}. For another
example, take the operation `do nothing' on \gfst{}\gone{}\glst{} to get a
way of composing \gfst{}\gthree{}{}{}{}{}\glst{}, and similarly the operation
\[
\gfst{}%
\gone{f}%
\gblw{}%
\gone{g}%
\glst{}
\diagspace\goesto\diagspace
\gfst{}%
\gone{gf}%
\glst{}
\]
on 
\[
\gfst{}\gone{}\gblw{}\gone{}\glst{}
\]
to get a way of composing
\[
\gfst{}\gtwo{}{}{}\gblw{}\gtwo{}{}{}\glst{};
\]
putting these together gives two
different ways of composing 
\[
\gfst{}\gthree{}{}{}{}{}\gfbw{}\gthree{}{}{}{}{}\glst{}, 
\]
familiar from the interchange law. There is also the all-at-once way of
composing this shape, given by applying the contraction principle to
\[
\gfst{}%
\gone{f}%
\gblw{}%
\gone{g}%
\glst{}
\diagspace\goesto\diagspace
\gfst{}%
\gone{gf}%
\glst{}.
\]

At this point the reader might be wondering what has happened to the
\index{coherence!axioms}
coherence conditions that exist in, say, bicategories. The answer lies in how
we realise
\index{n-category@$n$-category!as degenerate $\omega$-category}%
\index{finite dimensionality}%
lax $n$-categories as lax $\omega$-categories: namely, a lax
$n$-category is a lax $\omega$-category all of whose $m$-cells are identities
when $m>n$. 
\index{interchange law}%
Take, for instance, lax 2-categories and the two traditional ways
of composing a diagram 
\[
\gfst{A}%
\gthree{f}{f'}{f''}{\alpha}{\alpha'}%
\gfbw{B}%
\gthree{g}{g'}{g''}{\beta}{\beta'}%
\glst{C}
\]
to give a 2-cell, i.e.\
\gfst{A}\gtwo{gf}{g''f''}{\gamma_1}\glst{C} 
and  
\gfst{A}\gtwo{gf}{g''f''}{\gamma_2}\glst{C} 
where $\gamma_{1} = (\beta'\of\beta)*(\alpha'\of\alpha)$ and $\gamma_{2} =
(\beta'*\alpha')\of(\beta*\alpha)$. As these two ways restrict to the same
way on the 1-dimensional source and target, there is a 3-cell of the form
\[
\gfst{A}\gthreecell{gf}{g''f''}{\gamma_1}{\gamma_2}{}\glst{C}: 
\]
for apply the contraction principle to the
\index{degeneracy!of globular diagram}%
degenerate 3-cell diagram
\[
\gfst{}%
\gthree{}{}{}{}{}%
\gfbw{}%
\gthree{}{}{}{}{}%
\glst{}.
\] 
But in a (lax) 2-category the only 3-cells are identities, so $\gamma_{1} =
\gamma_{2}$, as expected. In this manner we can see that the different ways of
composing a diagram of 2-cells in a lax 2-category depend only on what
compositions they induce on the diagrams of 1-cells at its source and target,
just as in a traditional bicategory; and, indeed, that bicategories and our
lax 2-categories are the same but for inessential detail (the 
\index{bias}%
`bias' issue mentioned on page~\pageref{sec:bias}).%
\index{contraction!principle|)}%
\index{composition!ways of composing|)}%

\section{Formal Outline}	\label{sec:formout}

As indicated previously, the language of general
operads allows us to state a definition of lax $\omega$-category
quite easily. A certain category \ess\ and monad \ust\ are defined, and a
certain \Cartpr-operad 
\index{K@$K$ (operad)}%
$K$; lax $\omega$-categories are then just $K$-algebras. The point where our
definition 
\index{Batanin!misrepresentation}%
departs from Batanin's is the choice of the operad $K$. For us,
$K$ will be the structure formed by the `ways of composing' discussed in the
previous section.

The category \ess\ is that of 
\index{globular set}%
globular sets, as defined in~\ref{eg:cart-mnds}(\ref{eg:glob-mnd});
thus the underlying graph of a lax $\omega$-category will be a globular set.
The monad \ust\ is the 
\index{omega-category@$\omega$-category!free|(}%
free strict $\omega$-category monad. In outline, a
\index{omega-category@$\omega$-category!strict|emph}%
\emph{strict $\omega$-category} consists of a globular set $X$ and for each
$n\geq k$ a binary composition function
\[
X(n) \times_{X(k)} X(n) \go X(n)
\]
and an identity function
\[
X(k) \go X(n),
\]
satisfying the appropriate source-target relations, and such that the
compositions and identities obey strictly all the possible associative,
identity and interchange laws. The forgetful functor from strict
$\omega$-categories to globular sets has a left adjoint, described,
essentially, in \cite{Bat} (and see also~\ref{sec:trees} below). It is
straightforward to calculate that the adjunction is monadic and that the
monad \stbk\ on the (cartesian) category of globular sets is cartesian. We
may therefore speak of \Cartpr-multicategories and operads. 

Our major task now is to define the operad $K$. \emph{En route}, I will
attempt to explain what is going on in elementary terms.

\section{Trees}		\label{sec:trees}
\index{tree!Batanin-style|(}

It is instructive to contemplate the globular set $\wun^*$, where
\[
\wun = (\cdots \pile{\rTo \\ \rTo} 1 \pile{\rTo \\ \rTo} \cdots
\pile{\rTo \\ \rTo} 1).	\\
\] 
The free strict $\omega$-category functor takes a globular set $X$ and
creates formally all possible composites in it, to make $X^*$. Thus a typical
element of $\wun^{*}(2)$ looks like 
\begin{equation}
\gfst{}%
\gfour{}{}{}{}{}{}{}%
\grgt{}%
\gone{}%
\glft{}%
\gtwo{}{}{}%
\glst{},
\label{pic:four-one-three-glob}
\end{equation}
where each $k$-cell drawn represents the unique member of $\wun(k)$. Note
that because of 
\index{composition!nullary}%
identities (which we think of throughout as nullary
composites), this diagram might be thought of as representing an element of
$\wun^{*}(n)$ for any $n\geq 2$. Let us call a globular diagram representing
an element of $\wun^{*}(n)$ an
\index{glob|emph}
\emph{$n$-glob}. (The collections of $m$-globs and of $n$-globs are
considered disjoint, when $m\neq n$.) This
2-glob~(\ref{pic:four-one-three-glob}) has a source and a target, both of
which are the 1-glob
\[
\gfst{}%
\gone{}%
\gblw{}%
\gone{}%
\gblw{}%
\gone{}%
\glst{}.
\]
(Since all cells in \wun\ have the same source and target---are
\index{endomorphism}%
`endomorphisms'---it is inevitable that the same should be true in $\wun^*$.)

Having described $\wun^*$ as a globular set, we next turn to its strict
$\omega$-category structure: in other words, how globs may be composed.

\index{composition!of globs|(}%
Typical binary compositions are illustrated by
\[
\epsfig{file=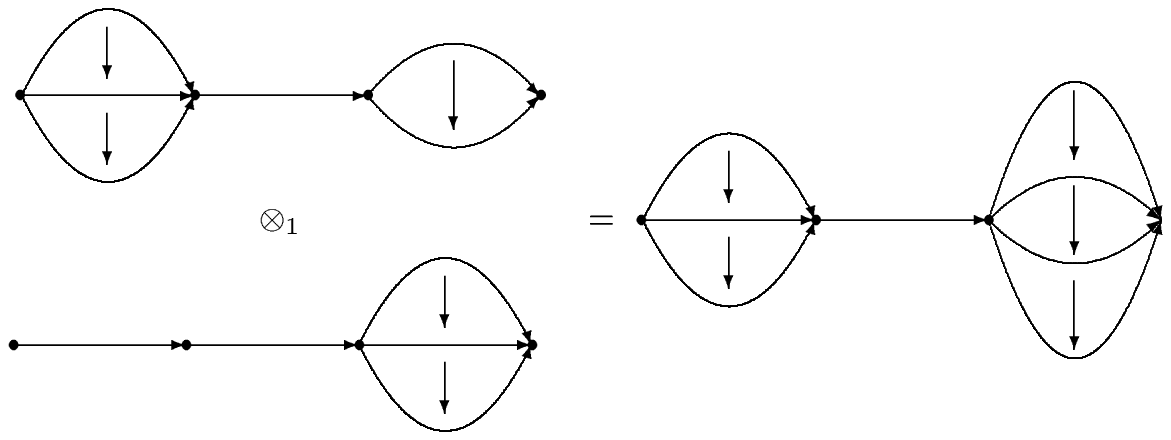}\label{pic:bin-comps}
\]
and
\[
\epsfig{file=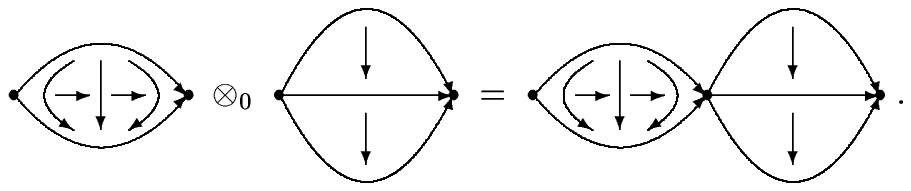}
\]
These compositions are possible because the sources/targets match
appropriately: e.g.\ in the first calculation, where we are gluing along
1-cells (indicated by $\otimes_1$), the 1-dimensional parts of the two
arguments are the same. A typical 
\index{composition!nullary}%
nullary composition---identity---is
\[
\begin{array}{ccc}	
\gfst{}\gone{}\gblw{}\gone{}\gblw{}\gone{}\glst{}	&
\goesto	&
\gfst{}\gone{}\gblw{}\gone{}\gblw{}\gone{}\glst{}	\\
\elt\wun^{*}(1)	&	&\elt\wun^{*}(2).
\end{array}
\]

Now in keeping with the 
\index{operadic philosophy}%
operadic philosophy, we do not wish to restrict our
attention merely to binary and nullary compositions, but rather to treat all
shapes of composition even-handedly. We may think of the first binary
composition above as 
\index{composition!indexing shape of}%
indexed%
\label{p:indexing}
by $\gfst{}\gthree{}{}{}{}{}\glst{} \elt \wun^{*}(2)$, because we were
composing one 2-cell with another by joining along their bounding 1-cells.
The composition can be represented as
\begin{equation}	\label{pic:first-bin-comp}
\epsfig{file=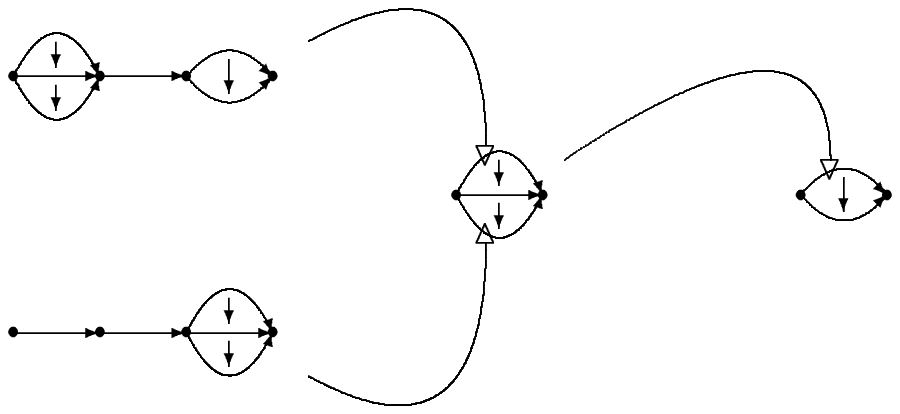}.
\end{equation}
In general, the ways of composing globs are indexed by globs themselves. For
instance,
\begin{equation}	\label{pic:gen-comps}
\epsfig{file=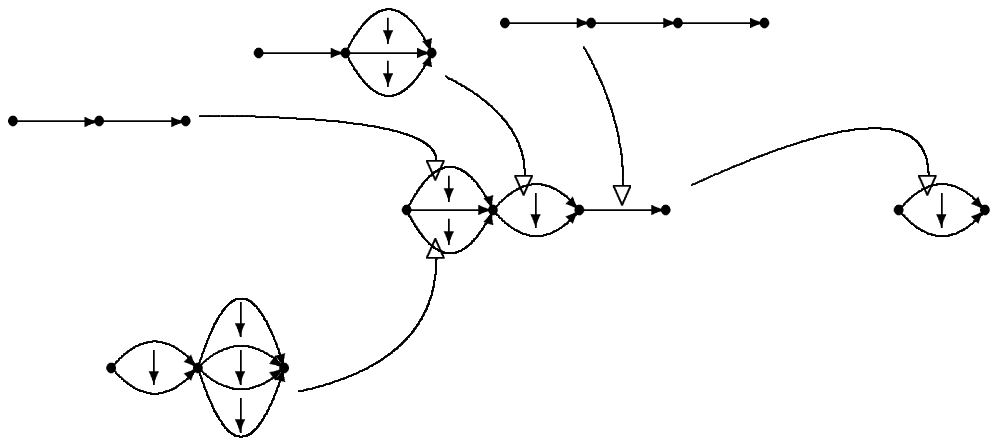}
\end{equation}
represents the composition
\[
\epsfig{file=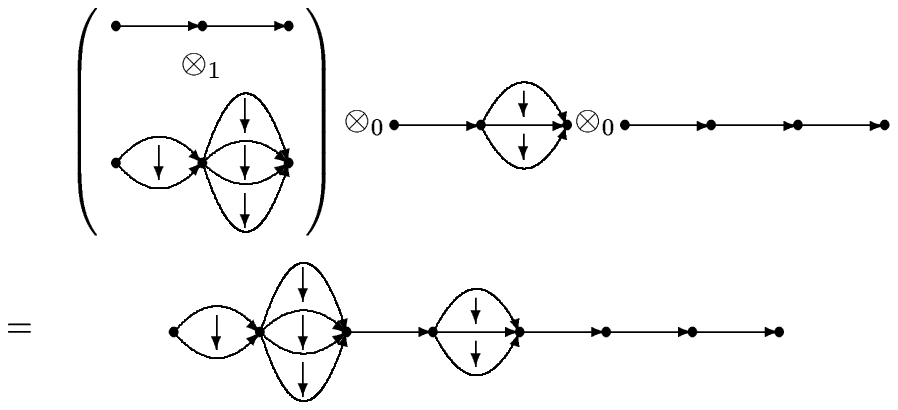}.
\]
\index{composition!of globs|)}%

We have now described the strict $\omega$-category $\wun^*$. The next
observation is that given a globular set $X$ and a glob $\tau$, there arises
a set 
\index{X@$X^\tau$!for globular set|emph}%
$X^\tau$. Formally, this is the fibre over $\tau$ in the map $X^{*}
\goby{!^*} \wun^{*}$. Informally, it's a 
\index{glob!labelled}%
labelling of $\tau$ in $X$. For
example, if $\tau$ is the glob of diagram~(\ref{pic:four-one-three-glob})
then $X^\tau$ has elements \tuplebts{A,B,C,D,
f,f',f'',f''',g,h,h', \alpha,\alpha',\alpha'',\beta} where $\range{A}{D}\elt
X(0)$, $\range{f}{h'}\elt X(1)$, $\range{\alpha}{\beta}\elt X(2)$, and
$s(\alpha)=f$, $t(\alpha)=f'$, etc., as in 
the picture
\[
\epsfig{file=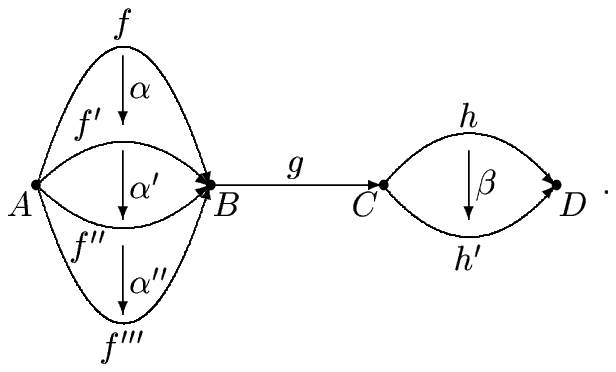}
\]

There is an alternative way to represent elements of $\wun^{*}(n)$, much
exploited by Batanin: as trees. (These trees differ from those of
Example~\ref{eg:cart-mnds}(\ref{eg:tree-mnd}) and Chapter~\ref{ch:opap}, in
both minor details and intent.) For instance, we translate the glob
\[
\gfst{}%
\gfour{}{}{}{}{}{}{}%
\grgt{}%
\gone{}%
\glft{}%
\gtwo{}{}{}%
\glst{}
\]
into the tree $\treedc$.  The thinking here is that the glob is 3 1-cells
long, so we start the tree as $\treec$; then the first column is 3 2-cells
high, the second 0, and the third 1, so the tree becomes $\treedc$; finally,
there are no 3-cells so the tree stops there. Formally, we may define an
\index{tree!Batanin-style|emph}%
\emph{$n$-stage tree} ($n\elt\nat$) to be a diagram
\[
\tau(n)\go\tau(n-1)\go\cdots\go\tau(1)\go\tau(0)=1
\]
in the category 
\index{simplicial category}%
$\Delta$ of all finite ordinals. The tree $\treedc \elt
\wun^{*}(2)$ corresponds to a certain diagram $4\go 3\go 1$ in $\Delta$, for
example; note that if $\tau$ is an $n$-stage tree with $\tau(n)=0$ then the
height of the picture of $\tau$ will be less than $n$. The 
\index{tree!Batanin-style!source and target of}%
source/target $\bdry\tau$ of an $n$-stage tree $\tau$ is the $(n-1)$-stage
tree obtained by removing all the nodes at height $n$, or formally,
truncating
\[
\tau(n)\go\tau(n-1)\go\cdots\go\tau(1)\go\tau(0)
\]
to
\[
\tau(n-1)\go\cdots\go\tau(1)\go\tau(0).
\]

It can be proved that the elements of $\wun^{*}(n)$ correspond one-to-one with
the $n$-stage trees, and that $\bdry$ provides the source and target (see
\cite{Bat} and the following paragraphs). Henceforth we assume this, write
\index{Tr@\Tr\ (Batanin trees)}%
\Tr\ for $\wun^{*}$, and prefer the word `tree' to `glob'. 

(An equivalent definition of tree uses the fact that an $(n+1)$-stage
tree is just a sequence of $n$-stage trees. Thus we may define $\Tr(0)=1$,
$\Tr(n+1)$ to be the free monoid%
\index{monoid!free}
on $\Tr(n)$, and
$\Tr(n+1)\goby{\bdry}\Tr(n)$ to be the free monoid functor applied to
$\Tr(n)\goby{\bdry}\Tr(n-1)$.) 

As a demonstration of the tree notation, the binary compositions of
page~\pageref{pic:bin-comps} look like
\begin{eqnarray*}
\treedcbang \otimes_{1} \treeccbang	&=	&\treegcbang	\\
\treebaa \otimes_{0} \treeba	&=	&\treebcb	
\end{eqnarray*}
and the 
\index{composition!nullary}%
nullary one looks like
\begin{eqnarray*}
\treec		&\goesto	&\treec		\\
\elt\Tr(1)	&		&\elt\Tr(2).	
\end{eqnarray*}
The 
\index{composition!indexing shape of}%
indexing shape of the first binary composition is $\treeba\elt\Tr(2)$,
and the whole composition is represented as
\[
\epsfig{file=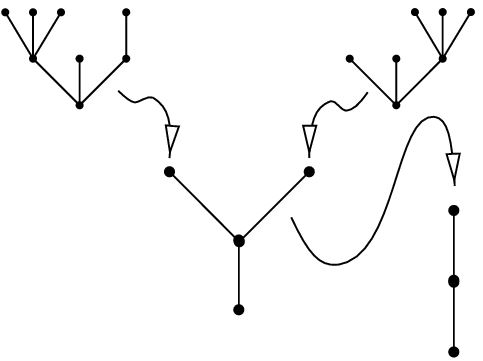}
\ .
\]
The more general, tree/glob-indexed composition on
page~\pageref{pic:gen-comps} translates to
\[
\epsfig{file=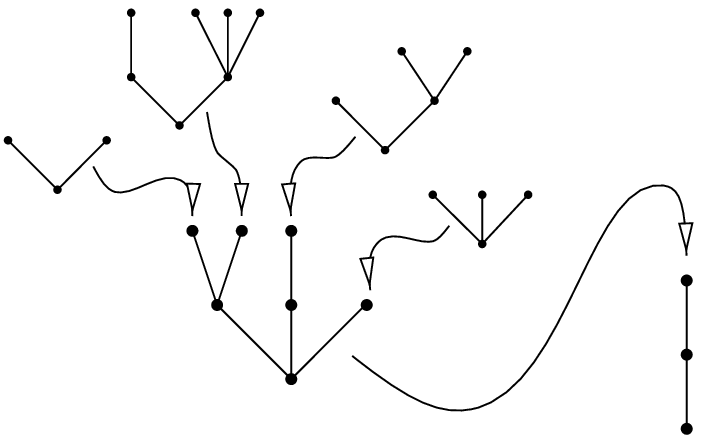}
\ ,
\]
representing the composition
\[
\left(\treeb\otimes_{1}\treedb\right) \otimes_{0}\treebb\otimes_{0}\treec
= \treefg.
\]

The foregoing considerations get us most of the way to a rigorous definition
of the 
\index{omega-category@$\omega$-category!free|emph}%
free strict $\omega$-category monad \stbk, as follows.
\begin{itemize}
\item
A tree is defined as a certain kind of diagram in $\Delta$ (as above).
\item
For each tree $\tau$ is defined a globular set 
\index{tau@$\hat{\tau}$!globular set}%
$\hat{\tau}$. This seems to be most easily done with the the second,
free-monoid, definition of tree, but we omit the details. (Or see \cite{Bat},
where $\hat{\tau}$ is called $T^*$.) Pictorially, we
take a tree such as $\treeec$, and draw a corresponding globular diagram
\[
\epsfig{file=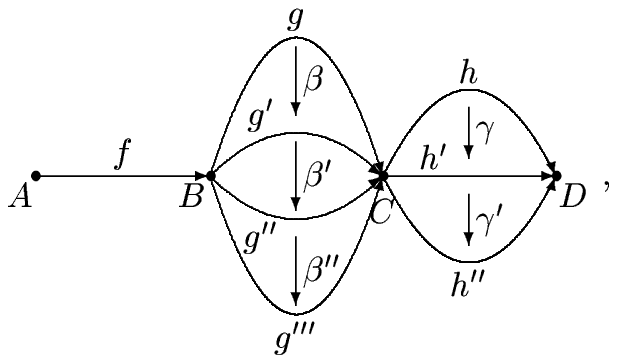}
\]
where $A, f, \beta, \ldots$ are formal symbols. This in turn
yields a 
globular set
\[
\hat{\tau}=(\cdots\emptyset
\pile{\rTo \\ \rTo}	\{\beta,\beta',\beta'',\gamma,\gamma'\}
\pile{\rTo \\ \rTo}	\{f,g,g',g'',g''',h,h',h''\}	
\pile{\rTo \\ \rTo}	\{A,B,C,D\}).	
\]

\item
Let $X$ be a globular set: then 
\index{X@$X^\tau$!for globular set}%
$X^\tau$ is defined as the set of maps
$\hat{\tau} \go X$ of globular sets.
\item
Finally, define $X^{*}(n)=\coprod_{\tau\elt\Tr(n)}X^{\tau}$.
\end{itemize}
With the rest of the details in place, it can then be shown that \stbk\
deserves the name `free strict $\omega$-category monad', and in particular
that \Tr\ is the free strict $\omega$-category on \wun.%
\index{omega-category@$\omega$-category!free|)}%
\index{tree!Batanin-style|)}

\section{Operads and Algebras}	\label{sec:ops-and-algs}
\index{operad!Batanin|(}%
\index{algebra!for Batanin operad|(}%

It has already been stated that if \ess\ is the category of globular sets and
\ust\ the free strict $\omega$-category monad, then \Cartpr-operads and their
algebras are what Batanin calls operads and algebras.%
\footnote{\label{fn:mgcs}%
In fact, \cite{Bat} is presented in the wider context of 
\index{monoidal globular category}%
monoidal globular categories, but this need not concern us here. The Batanin
operads referred to here are sometimes given the fuller name of
`$\omega$-operads in
\textit{Span}' in \cite{Bat}.}
We now explain these in elementary terms.

A 
\index{collection|emph}%
\emph{collection} is a graph 
\begin{slopeydiag}
		&		&C	&		&	\\
		&\ldTo<{d}	&	&\rdTo>{!}	&	\\
\Tr=\wun^{*}	&		&	&		&\wun	\\
\end{slopeydiag}
in \ess. In other words, a collection consists of a set $C(\tau)$ for each
tree $\tau$, together with a pair of functions
$
\begin{diagram}
C(\tau)	&\pile{\rTo^{s} \\ \rTo_{t}}	&C(\bdry\tau),	\\
\end{diagram}
$
satisfying the usual globularity relations $ss=st$ and $ts=tt$. An operad
structure on the collection $C$ consists of identities and compositions with
suitable properties. The identities consist of an element of $C(\upsilon_n)$
for each $n$, where 
\index{upsilon@$\upsilon_{n}$}%
\index{tree!Batanin-style!straight}%
$\upsilon_{n} \elt \Tr(n)$ is the tree
\[
\left.\treeaaq\ \right\} n \mr{\ edges}
\]
(corresponding to the glob which is a single $n$-cell). For 
\index{composition!in Batanin operad}%
composition in $C$,
consider a diagram
\begin{equation}		\label{pic:compn-in-operad}
\epsfig{file=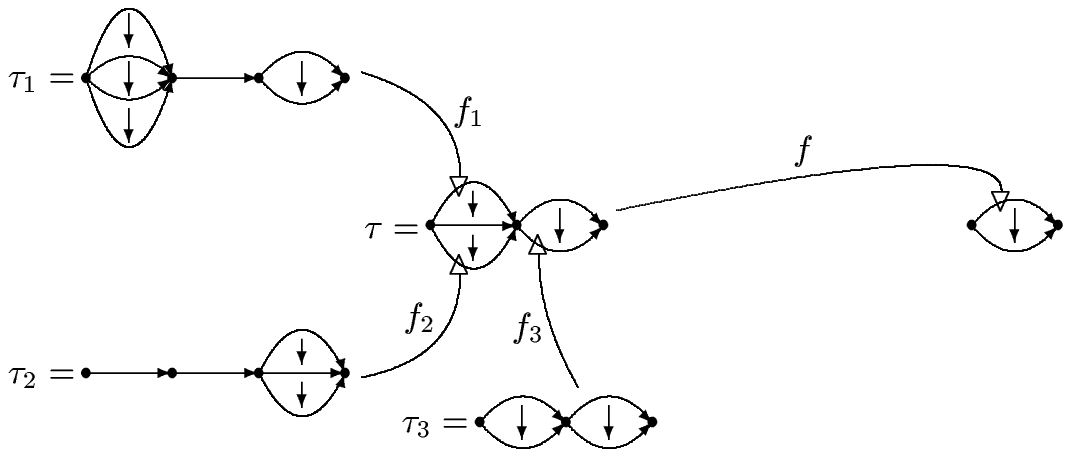}.
\end{equation}
This is meant to indicate that $f_{1}\elt C(\tau_{1})$, $f_{2}\elt
C(\tau_{2})$, $f_{3}\elt C(\tau_{3})$, $f\elt C(\tau)$, and that $f_{1}$,
$f_{2}$, $f_3$ match 
suitably on their sources and targets (e.g.\ $t(f_{1}) =
s(f_{2})$). Then composition should produce from this data a diagram
\[
\epsfig{file=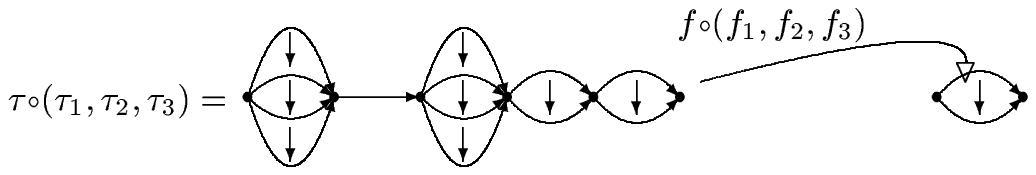}\ ,
\]
i.e.\ an element of $C(\tau\of(\tau_{1},\tau_{2},\tau_{3}))$.

According to the general theory, composition is a map $C\of C \go C$ over
\Tr, where $C\of C$ is given by the diagram
\[
\begin{slopeydiag}
   &       &   &       &   &       &C\of C\Spbk&  &   \\
   &       &   &       &   &\ldTo<{\mr{pr}}&&\rdTo&   \\
   &       &   &       &C^*&       &      &       &C  \\
   &       &   &\ldTo<{d^*}&&\rdTo<{!^*}& &\ldTo>{d}& \\
   &       &1^{**}&    &   &       &\Tr   &       &   \\
   &\ldTo<{\mu_1}&&    &   &       &      &       &   \\
\Tr&       &   &       &   &       &      &       &   \\
\end{slopeydiag}
\ \ \ \ .
\]
To see what a typical element of $C\of C$ is, consider $\tau$ above, an
element of $\Tr(2)$. In the map $C^{*}\goby{!^*}\Tr$, a typical element of
the fibre over $\tau$ is the left-hand half of
diagram~(\ref{pic:compn-in-operad}). In the map $C\goby{d}\Tr$, a typical
element of the fibre over $\tau$ is the right-hand half
of~(\ref{pic:compn-in-operad}). The map $\gobyc{\mu_{1}\of d^{*} \of
\mr{pr}}{C\of C}{\Tr}$ sends~(\ref{pic:compn-in-operad}) to the tree
$\tau\of(\tau_{1},\tau_{2},\tau_{3})$. So, as hoped for above, composition
takes data such as~(\ref{pic:compn-in-operad}) and produces an element of
$C(\tau\of(\tau_{1},\tau_{2},\tau_{3}))$. This composition is moreover
required to be associative and to obey unit laws with respect to the
identities. For example, if we have a diagram
\begin{diagram}[height=1.5em]
\cdot	&	&	&	&	&	&	\\
	&\rdTo>{f_{11}}&&	&	&	&	\\
\cdot	&\rTo^{f_{12}}&\cdot&	&	&	&	\\
	&\ruTo>{f_{13}}&&\rdTo(2,3)>{f_1}&&	&	\\
\cdot	&	&	&	&	&	&	\\
	&	&	&	&\cdot	&\rTo^{f}&\cdot	\\
	&	&	&\ruTo>{f_2}&	&	&	\\
\cdot	&\rTo_{1}&\cdot	&	&	&	&	\\
\end{diagram}
of the kind~(\ref{pic:compn-in-operad}), then
\[
f\of(f_{1}\of(f_{11},f_{12},f_{13}),f_{2}) =
(f\of(f_{1},f_{2})) \of (f_{11},f_{12},f_{13},1).
\]

We have now seen that an operad consists of a set $C(\tau)$ for each $\tau$,
with source and target functions, and compositions between the $C(\tau)$'s
according to the compositions for trees. Think of an element of $C(\tau)$ as
a way of composing a diagram shaped like the glob $\tau$; then it makes
sense that an algebra for $C$ should be a globular set $X$ with for each
$f\elt C(\tau)$ a function $X^{\tau} \goby{\bar{f}} X(n)$, such that
composition of these functions $\bar{f}$ commutes with composition in the
operad. So for instance, if 
\[
f \elt C\left(
\gfst{}\gfour{}{}{}{}{}{}{}\grgt{}\gone{}\glst{}
\right)
\]
and 
\[
\gfst{A}%
\gfour{p}{p'}{p''}{p'''}{\alpha}{\alpha'}{\alpha''}%
\grgt{B}%
\gone{q}%
\glst{C}
\]
is a picture in $X$, then $\bar{f}$ composes this picture
to give a 2-cell of $X$. This is indeed what the general theory says an
algebra is. For an algebra structure on $X$ is a suitable map $X_{\blob}
\goby{h} X$, where $X_{\blob}$ is the pullback
\[
\begin{slopeydiag}
	&		&X_{\blob}\Spbk	&		&	\\
	&\ldTo		&		&\rdTo		&	\\
X^*	&		&		&		&C	\\
	&\rdTo<{!^*}	&		&\ldTo>{d}	&	\\
	&		&\Tr		&		&	\\
\end{slopeydiag}
\ \ \ \ ;
\]
this means that 
$X_{\blob}(n)=\coprod_{\tau\elt\Tr(n)}(X^{\tau}\times C(\tau))$ 
and that $h$ is a sequence of functions 
$\left(\coprod_{\tau\elt\Tr(n)}(X^{\tau}\times C(\tau))
\go X(n) \right)_{n\elt\nat}$ 
obeying the usual laws.

We have now discussed what operads and their algebras look like, and it is
time to come to the main point of the chapter.%
\index{operad!Batanin|)}%
\index{algebra!for Batanin operad|)}%

\section{Contractions and Lax $n$-Categories}	\label{sec:con-and-n-cat}

This section formalizes the 
\index{contraction!principle}%
`contraction principle' of~\ref{sec:infout},
allowing us to define the operad $K$ whose algebras are the lax
$\omega$-categories.

If $C$ is a collection, $\sigma$ a tree, and $f, f' \elt C(\sigma)$, let
us say \pr{f}{f'} is a 
\index{matching pair|emph}%
\emph{matching pair} for $\sigma$ if $s(f)=s(f')\elt
C(\bdry\sigma)$ and $t(f)=t(f')\elt C(\bdry\sigma)$. Here it is intended that
if $\sigma$ is the 0-stage tree then any pair $f,f'\elt C(\sigma)$ is
matching.
\begin{defn}	\label{p:defn-contr}
Let \slob{C}{}{\Tr} be a collection. A 
\index{contraction|emph}%
\emph{contraction} $\psi$ on $C$ is a function assigning to each
$\tau\elt\Tr(n)$ ($n\geq 1$) and matching pair
\pr{f}{f'} for $\bdry\tau$, an element $\psi\pr{f}{f'} \elt C(\tau)$ such that
$s(\psi\pr{f}{f'})=f$ and $t(\psi\pr{f}{f'})=f'$.
\end{defn}

We use the category of 
\index{operad!with contraction@-with-contraction}%
\emph{operads-with-contraction}, in which an object is
an operad with a specified contraction (on the underlying collection), and
the maps are the operad maps preserving the specified contractions.

It should be apparent that the operad for lax $\omega$-categories, described
informally in~\ref{sec:infout}, carries a contraction. Indeed, the presence of
the contraction was all that we used to generate the operad: thus we define
\index{K@$K$ (operad)|emph}%
$K$ to be the initial operad-with-contraction. It is of course necessary to
verify that such an initial object exists. Here it is hoped that the reader
will be content with the heuristic argument laid out in~\ref{sec:infout},
where starting from the empty collection it was suggested how the operad and
contraction structures could be added in freely.

Let us say a globular set $X$ is 
\index{globular set!finite-dimensional}%
\index{finite dimensionality}%
\index{n-dimensional@$n$-dimensional}%
\emph{$n$-dimensional} if for all $m\geq n$,
\[
\gobyc{s=t}{X(m+1)}{X(m)}
\]
and this map is an isomorphism. 
\begin{defn}
\index{omega-category@$\omega$-category!lax|emph}%
\index{n-category@$n$-category!lax|emph}%
A \emph{lax $\omega$-category} is an algebra for $K$. A \emph{lax
$n$-category} is a lax $\omega$-category whose underlying globular set is
$n$-dimensional.
\end{defn}
We remark that despite being able to make this definition, there is no
immediately obvious way to define 
\index{lax vs strict@lax \emph{vs.} strict}%
\index{strict!vs lax@\emph{vs.} lax}%
\index{functor!lax \emph{vs.} strict}%
lax functors, transformations, \ldots. In
particular, the category of $K$-algebras consists of lax $\omega$-categories
and strict $\omega$-functors.

Our final observation is that the terminal operad \slob{\Tr}{1}{\Tr} has a
unique contraction on it, and is then the terminal operad-with-contraction.
As always, the algebras for the terminal operad are the algebras for the monad
\ust\ on \ess\ (see~\ref{eg:algs}(\ref{eg:terminal})). So the algebras for
the terminal operad-with-contraction are the strict $\omega$-categories, and
the algebras for the initial operad-with-contraction are the lax
$\omega$-categories. One might therefore hope that in some sense the category
of operads-with-contraction indexes the different possible strengths for
theories of $\omega$-categories.

\section{The Cubical Approach}		\label{sec:cub-ap}
\index{cubical!category|(}

2-categories, strict and lax, are not the only 2-dimensional generalization
of categories: there is also the notion of a double category. This section
presents the beginnings of an approach to $n$-tuple categories. It runs in
very close parallel with our approach to $n$-categories, hence its
inclusion in this chapter.%
\index{Batanin|)}%
\index{n-category@$n$-category!lax|)}

A 
\index{double category|emph}%
\emph{double category} may be defined as a category object in
\fcat{Cat}. More descriptively, the graph structure consists of collections
of
\begin{itemize}
	\item 0-cells $A$
	\item horizontal 1-cells $f$
	\item vertical 1-cells $p$
	\item 2-cells $\alpha$
\end{itemize}
and various source and target functions, as illustrated by the picture
\[
\begin{diagram}
A_1		&\rTo^{f_1}		&A_2		\\
\dTo<{p_1}	&\alpha			&\dTo>{p_2}	\\
A_3		&\rTo_{f_2}		&A_4		\\
\end{diagram}
\ .
\]
The category structure consists of identities and composition functions for
2-cells and both kinds of 1-cell, obeying strict associativity, identity and
interchange laws; see 
\index{Kelly!-Street}%
\cite{KS} for more details. It is worth keeping in mind
that a 2-category is a double category in which all the vertical 1-cells are
identities, and similar 
\index{degeneracy!of cubical set}%
degeneracy properties should hold in higher dimensions.

More generally, let us define 
\index{cubical!set|emph}%
\emph{$n$-cubical set} for any $n\elt\nat$; the
intention is that a 2-cubical set will be the underlying graph of a double
category. So, let 
\index{Cube@\fcat{Cube_n}|emph}%
\fcat{Cube_n} be the category with
\begin{description}
	\item[objects:] subsets $D$ of $\{ \range{0,1}{n-1} \}$ 
	\item[maps $D \go D'$:] the inclusion $D\sub D'$, together with 
	a function $D'\without D \go \{0,1\}$ (in which context 0 should be
	read as `source' and 1 as `target') 
	\item[composition:] place functions side-by-side.
\end{description}
Then we define an $n$-cubical set to be a functor
$\fcat{Cube_n}^{\op}\go\Set$. For instance, we may think of a 2-cubical set
$X$ as:
\begin{itemize}
	\item $X\emptyset = \{ \mbox{0-cells} \}$
	\item $X\{0\} = \{ \mbox{horizontal 1-cells} \}$
	\item $X\{1\} = \{ \mbox{vertical 1-cells} \}$
	\item $X\{0,1\} = \{ \mbox{2-cells} \}$
\end{itemize}
and, for instance, the map $\{1\} \go \{0,1\}$ given by
\[
\{0,1\} \without \{1\} = \{0\} \goby{0} \{0,1\}
\]
sends $\alpha\elt X \{0,1\}$ to $p_{1} \elt X\{1\}$, in the diagram above.

The category \ess\ of $n$-cubical sets is a presheaf category, therefore
cartesian. A \emph{strict $n$-tuple category} is an $n$-cubical set together
with various compositions and identities, as for double categories, all
obeying strict laws. The free strict $n$-tuple category $X^*$ on an
$n$-cubical set $X$ consists of all formal composites of cells in $X$. In
particular, cells of $\wun^*$ look like cuboids. For example, a typical
element of $\wun^{*}\{0,1,2\}$ may be depicted as
\[
\epsfig{file=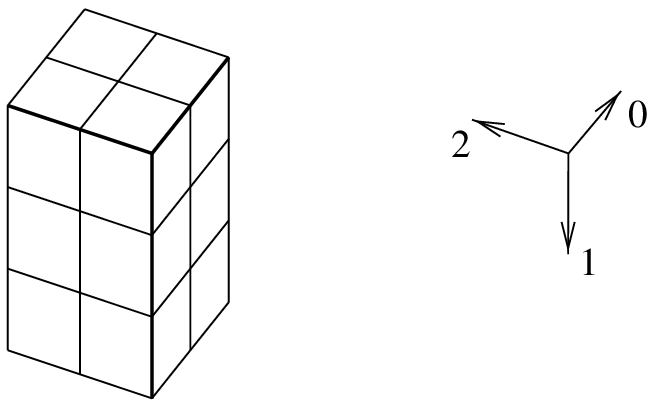}
\]
(where the device on the right is just a `compass' to show which element of
$\{0,1,2\}$ points in which direction). Equally, this cuboid could be used to
illustrate an element of $\wun^{*}(D)$ for any $D\supseteq
\{0,1,2\}$. Composition in $\wun^{*}$ is indicated by
\begin{equation}
\epsfig{file=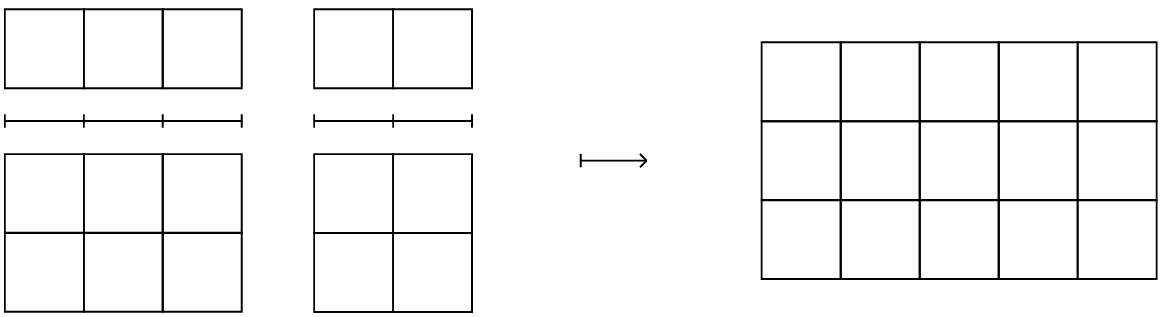,height=60pt}\ ,		\label{pic:cube-comp}
\end{equation}
the 
\index{composition!indexing shape of}%
indexing shape of this composition being
\begin{equation}
\epsfig{file=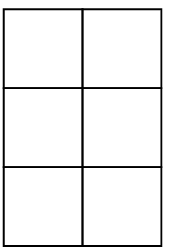,height=47pt}		\label{pic:two-by-three}
\end{equation}
(cf.\ the discussion of indexing on page~\pageref{p:indexing}). Given an
$n$-cubical set $X$ and a cuboid $\tau\elt\wun^{*}(D)$, there arises a set
\index{X@$X^\tau$!for cubical set|emph}%
$X^\tau$. Formally, this is the fibre over $\tau$ in the map $X^{*}
\goby{!^*} \wun^*$; informally, $X^\tau$ consists of ways of labelling the
cuboid $\tau$ in $X$.

In parallel to the representation of a glob as a diagram
$\tau(n)\go\cdots\go\tau(0)$ in $\Delta$, a cuboid can be represented as its
sequence of edge-lengths. So, define a 
\index{cuboid}%
\emph{$D$-cuboid} ($D\sub\{\range{0,1}{n-1}\}$) to be a map $D\go\nat$. To
each cuboid $\tau$ is associated an $n$-cubical set 
\index{tau@$\hat{\tau}$!cubical set}%
$\hat{\tau}$: e.g.\ if
$n=2$ and $\tau$ is the cuboid~(\ref{pic:two-by-three}) above then
$\hat{\tau}$ is a certain 2-cubical set with 6 2-cells, 8 horizontal 1-cells,
9 vertical 1-cells and 12 0-cells. Then define $X^\tau$ to be the set of maps
$\hat{\tau}\go X$, for any $n$-cubical set $X$, and define
$X^{*}(D)=\coprod_{(D\mr{-cuboids\ }\tau)}X^{\tau}$. This goes most of the
way towards a rigorous definition of the functor \stbk, the monad structure
of which can then be described; the resulting monad \ust\ on \ess\ is
cartesian. (The manipulations involved are rather easier here than in the
globular setting, because we have cuboids instead of trees.)

\index{operad!cubical}%
\index{cubical!operad}%
\index{operad!Batanin}%
\Cartpr-operads can be understood in much the same way as Batanin's operads.
A Batanin operad associates to each tree a set, and has composition functions
corresponding to the composition of trees; a cubical operad associates to
each cuboid a set, and has composition functions corresponding to the
composition of cuboids. If $C$ is an operad and $\tau$ a cuboid then an
element of $C(\tau)$ can be thought of as a way of composing a diagram of
shape $\tau$; thus an algebra for $C$ consists of a cubical set $X$ together
with suitably compatible functors $X^{\tau} \goby{\bar{f}} X(D)$ ($f\elt
C(\tau)$, $\tau \elt \wun^{*}(D)$). (Note that the restriction to
\emph{finite} sets $\{\range{0,1}{n-1}\}$ has been quite unnecessary, and was
only done for the sake of simplicity.)

The final step of the analogy would be to develop an appropriate notion of
\index{contraction!cubical version}%
contraction in the cubical setting, in order to define lax $n$-tuple
categories. This has not, as far as I know, been carried out.

\index{cubical!category|)}

%% file: essgray.tex
\chapter{Gray-Categories}		\label{ch:gray}
\index{Gray-category|(}

Gray-categories are significant because of their place in 
\index{coherence!for tricategories|(}%
coherence for tricategories, itself perhaps the most challenging feature of
the higher-dimensional landscape. Whereas any bicategory is biequivalent to a
2-category, it is not true that any tricategory is triequivalent to a
3-category.  Instead, as revealed in
\index{Gordon-Power-Street}%
\cite{GPS}, any tricategory is triequivalent to a Gray-category.

In this chapter we take the view that 
\index{Gray-category!viewpoint on}%
Gray-categories are structures in their own right. It is only by making a
\index{non-canonical choice|(}%
non-canonical choice that we can embed the class of Gray-categories in the
class of tricategories, and so speak of Gray-categories as a particular kind
of tricategory, as in the coherence theorem just stated. An important case of
this is that 2-categories (and homomorphisms, \ldots) naturally form a
Gray-category, not a tricategory; this boils down to the fact that there is
no canonical way to define the horizontal composite of strong
transformations. Similarly, bicategories do not naturally form a tricategory,
only a 
\index{near-Gray-category}%
`near-Gray-category'. This raises questions over the conventional
wisdom that lax 
\index{n+1-category@$(n+1)$-category of $n$-categories}%
$n$-categories should form a lax $(n+1)$-category.

The necessity of this non-canonical choice seems intimately related to the
coherence issue. 
\index{monoidal category!braided}%
\index{braided monoidal category}%
Braided monoidal categories (BMC's; see 
\index{Joyal-Street}%
\cite{JS}) are a
case in point. Any BMC `is' (after some non-canonical choices) a tricategory
with one 0-cell and one 1-cell. Equivalences of BMC's then correspond to
triequivalences of their tricategories. Moreover, the property of a braiding
being a 
\index{monoidal category!symmetric}%
\index{symmetric monoidal category}%
symmetry is preserved under BMC-equivalence, and the braiding is a
symmetry if the tricategory is a 3-category. It follows that any BMC whose
braiding is not a symmetry provides an example of a tricategory not
triequivalent to a 3-category. (It doesn't quite follow immediately as the
3-category needn't be a BMC; see 
\index{Gordon-Power-Street}%
\cite[chapter 8]{GPS}.) Now braidings on a given
monoidal category come in pairs
(\gobyc{\gamma_{AB},\gamma_{BA}^{-1}}{A\otimes B}{B\otimes A}, so to speak);
to choose one particular braiding over its partner is non-canonical. The
exception is when the braiding is equal to its partner, which happens exactly
when the braiding is a symmetry. It is hard not to wonder if this
non-canonical choice is related to the non-canonical choice involved (above)
for Gray-categories, particularly given the nature of the latter choice.
However, these connections remain mysterious, so we leave the matter there.%
\index{coherence!for tricategories|)}%
\index{non-canonical choice|)}%

\index{Cayley representation}%
The Cayley representation functor for categories has an obvious analogue in
bicategories and tricategories. We examine its properties, and thus obtain a
\index{Gray-category!Cayley representation of}%
\index{Cayley representation}%
representation theorem for Gray-categories: any small Gray-category `is' a
substructure of \tooc. This means, of course, that any manipulation of
diagrams in a general Gray-category might as well be done in \tooc\ (in the
same sense as is true for categories and \Set); however, once one knows the
definition of a Gray-category, it is probably less distracting to work in the
general context. A more interesting aspect of our representation theorem is
its relation to the coherence theorem. This leads us to make a conjecture on
coherence theorems for higher-dimensional categories. 

Section~\ref{sec:ssq} is on sesquicategories, which provide our route into
defining a Gray-category (\ref{sec:defn-gray}). We exhibit \tooc\ as a
Gray-category in~\ref{sec:tooc-gray}, and any Gray-category as a tricategory
(non-canonically) in~\ref{sec:gray-tri}. In~\ref{sec:homms-trieqs} we pause to
gather some definitions, before moving on to the Cayley characterization of
Gray-categories (\ref{sec:char-of-gray}) and a conjecture on coherence in
higher dimensions (\ref{sec:hd-coh-conj}).

\section{Sesquicategories}	\label{sec:ssq}
\index{sesquicategory|(}%

In short, a 
\index{sesquicategory|emph}%
\emph{sesquicategory} consists of a category \scat{C} and a
factorization
\begin{diagram}
\scat{C}^{\op} \times \scat{C}	&\rGet		&\Cat		\\
				&\rdTo<{\Hom}	&\dTo>{\ob}	\\
				&		&\Set.		\\
\end{diagram}
Thus for each pair of objects $A,B\elt \scat{C}_0$ there is a category
\homset{\scat{C}}{A}{B} whose objects are the morphisms from $A$ to $B$. We
may think of the objects of \scat{C} as 0-cells, the morphisms in \scat{C} as
1-cells, and the morphisms in the categories \homset{\scat{C}}{A}{B} as
2-cells. For instance, any 2-category has an underlying sesquicategory. 

\begin{figure}
\epsfig{file=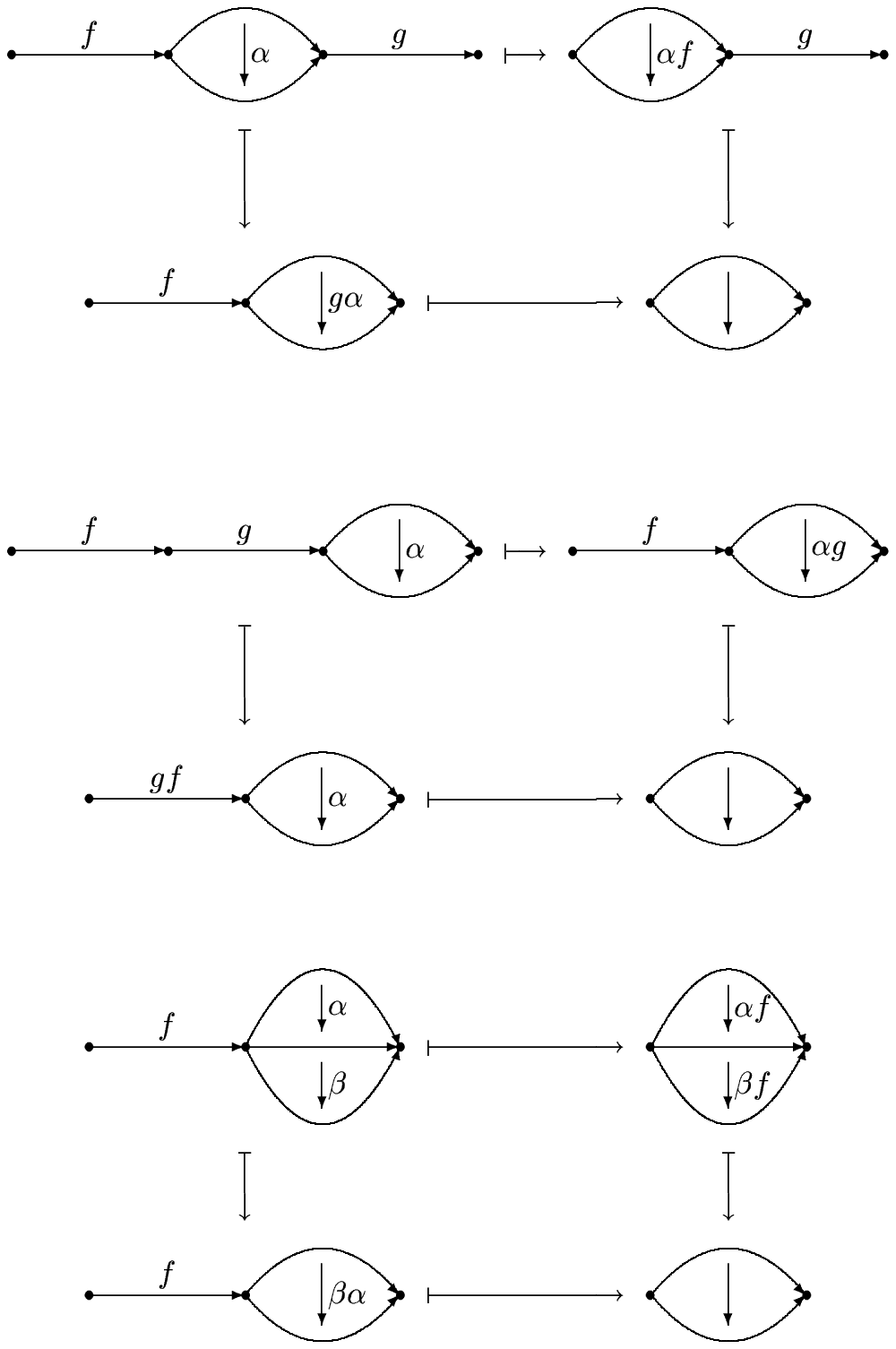}
\caption{Some of the axioms for a sesquicategory. The first diagram, for
instance, says that $(g\alpha)f = g(\alpha f)$.}
\label{fig:ssq-axioms}
\end{figure}

The definition of a sesquicategory can be recast as follows. The graph
structure of a sesquicategory \scat{C} is the same as that of a bicategory:
that is, a globular set
\begin{diagram}
\scat{C}_{2} &\pile{\rTo \\ \rTo} &\scat{C}_{1} &\pile{\rTo \\ \rTo} &
\scat{C}_{0} \\
\end{diagram}
truncated at dimension 2. For 1-cells there are 
\index{composition!in sesquicategory}%
binary compositions and identities,
\[
\gfst{}\gone{}\gblw{}\gone{}\glst{}
\diagspace\goesto\diagspace
\gfst{}\gone{}\glst{}
\]
and			
\[
\gzero{}				
\diagspace\goesto\diagspace
\gfst{}\gone{}\glst{},	
\]
and for 2-cells there are vertical binary compositions and identities,
\[
\gfst{}\gthree{}{}{}{}{}\glst{}		
\diagspace\goesto\diagspace
\gfst{}\gtwo{}{}{}\glst{}
\]
and
\[
\gfst{}\gone{}\glst{}			
\diagspace\goesto\diagspace
\gfst{}\gtwo{}{}{}\glst{}.
\]
Finally, there are compositions as illustrated by
\index{whiskering|(}%
\[
\gfst{}\gone{f}\gblw{}\gtwo{}{}{\alpha}\glst{}	
\diagspace\goesto\diagspace
\gfst{}\gtwo{}{}{\alpha f}\glst{}
\]
and
\[
\gfst{}\gtwo{}{}{\alpha}\gblw{}\gone{f}\glst{}	
\diagspace\goesto\diagspace
\gfst{}\gtwo{}{}{f\alpha}\glst{}.
\]
(This last kind of composition arises because if $B\go B'$ in \scat{C} then
the factorization $\scat{C}^{\op} \times \scat{C} \go \Cat$ gives a functor
$\homset{\scat{C}}{A}{B} \go \homset{\scat{C}}{A}{B'}$.) The compositions and
identities are subject to various axioms. 1-cell composition and vertical
2-cell composition are to obey associative and identity laws---in other
words, the 0- and 1-cells must form a category, as must the 1- and 2-cells.
The remaining axioms are indicated in Figure~\ref{fig:ssq-axioms}.
There, the last
two diagrams are binary and left-handed in character; we should really have
drawn six more to cover the nullary and/or right-handed versions too.

The crucial difference between a sesquicategory and a 2-category is that in a
sesquicategory there is no specified 
\index{composition!horizontal}%
\index{non-canonical choice}%
horizontal composition of 2-cells. We
can derive one, according to the picture
\begin{equation}	\label{pic:derived-horiz}
\epsfig{file=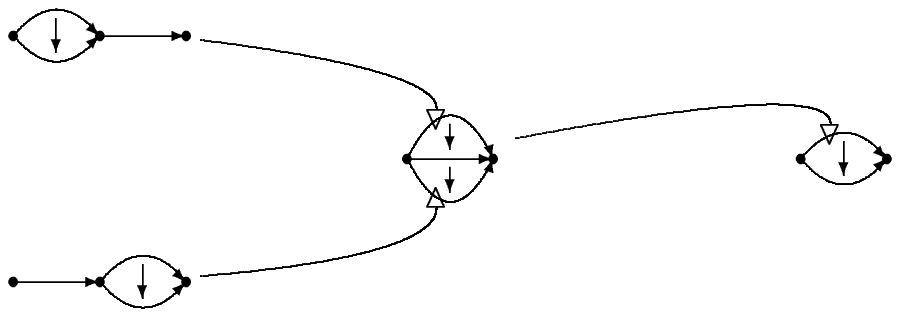}
\end{equation}
(cf.\ page~\pageref{pic:first-bin-comp}, diagrams~(\ref{pic:first-bin-comp})
and~(\ref{pic:gen-comps})), but equally we could have used the symmetrically
opposite version, and these two horizontal compositions are not the same. So
if 
\[
\gfst{}%
\gtwo{f}{f'}{\alpha}%
\gfbw{}%
\gtwo{g}{g'}{\beta}%
\glst{}
\]
is a diagram in a sesquicategory then there are 2-cells
\[
\gfst{}\gtwocentre{gf}{g'f'}{\beta f'\of g\alpha}\glst{} 
\diagspace\mr{and}\diagspace 
\gfst{}\gtwocentre{gf}{g'f'}{g'\alpha\of\beta f}\glst{},
\]
but in general they will be different.
We return to this point repeatedly.

Our elementary description of sesquicategories allows us to realise them as
algebras for a certain 
\index{operad!Batanin!for sesquicategories}%
Batanin operad. In this chapter%
\label{p:fin-dim-attitude}
we will be concerned with 
\index{finite dimensionality}%
finite-dimensional structures, so use $n$-globular sets
\begin{diagram}
X(n) &\pile{\rTo \\ \rTo} &X(n-1) &\pile{\rTo \\ \rTo} &\cdots&\pile{\rTo \\
\rTo}& 
X(0) \\
\end{diagram}
and $n$-operads (that is, \Cartpr-operads where \ust\ is the free strict
$n$-category monad) rather than the $\omega$-versions favoured in
Chapter~\ref{ch:glob}.  The 2-operad 
\index{Ssq@\Ssq|emph}%
\Ssq, whose algebras are
sesquicategories, is generated by a single element of $\Ssq(\tau)$ for the
tree $\tau$ taking any one of the values
\begin{eqnarray*}
\gfst{}\gone{}\gblw{}\gone{}\glst{}, \diagspace
\gzero{}
&\mbox{(1-stage trees),}	\\
\gfst{}\gthree{}{}{}{}{}\glst{}, \diagspace
\gfst{}\gone{}\glst{}, 		&\\
\gfst{}\gone{}\glft{}\gtwo{}{}{}\glst{}, \diagspace
\gfst{}\gtwo{}{}{}\grgt{}\gone{}\glst{}
&\mbox{(2-stage trees),}
\end{eqnarray*}
subject to the relations described above. That it is possible to specify an
operad by 
\index{generators and relations}%
\index{operad!Batanin!by generators and relations}%
generators and relations is established in \cite{Bat}; given that
(by~\ref{sec:free-multi}) the category of operads is monadic over the
category of collections, itself a presheaf category, this does not seem
surprising.%
\index{whiskering|)}%
\index{sesquicategory|)}%
\index{composition!in sesquicategory}%

\section{The Definition of a Gray-Category}	\label{sec:defn-gray}

The original definition 
\index{Gordon-Power-Street}%
(\cite{GPS}) of a Gray-category was as a category 
\index{enrichment!in Gray@in \textbf{Gray}}%
enriched in \textbf{Gray}, where \textbf{Gray} is the category of
2-categories with a certain symmetric monoidal structure. Here we ignore this
definition, favouring the equivalent one of
\index{Batanin!definition of Gray-category}%
\cite[p.\ 59]{Bat}, and therefore
write `Gray-category' rather than `\textbf{Gray}-category'.

\index{Gray-category|emph}%
\emph{Gray-categories} are defined to be algebras for the 3-operad 
\index{Gy@\Gy|emph}%
\Gy. This in turn is defined to be the 3-operad got by taking the 2-operad
\Ssq\ and contracting from dimension 2 to dimension 3. In other words, for
$\tau\elt\Tr(n)$,
\[
\Gy(\tau)=
\left\{
\begin{array}{ll}
\Ssq(\tau),		&\mbox{\ $n\elt\{0,1,2\}$},\\
\Ssq(\bdry\tau)\times \Ssq(\bdry\tau),		&\mbox{\ $n=3$}.
\end{array}
\right.
\]
(Note that if $\tau$ is a 3-stage tree then $\Ssq(\bdry^{2}\tau)$ has just
one element; thus any two elements of $\Ssq(\bdry\tau)$ automatically match
at their source and target.) The operad structure on \Gy\ is uniquely
determined by saying that the 2-operad it restricts to is \Ssq.

The underlying graph of a Gray-category is therefore a 3-globular set. In
dimensions 0--2 there are the same compositions as in \Ssq, obeying the same
laws. A way of 
\index{composition!in Gray-category}%
composing a 3-cell diagram such as
\[
\epsfig{file=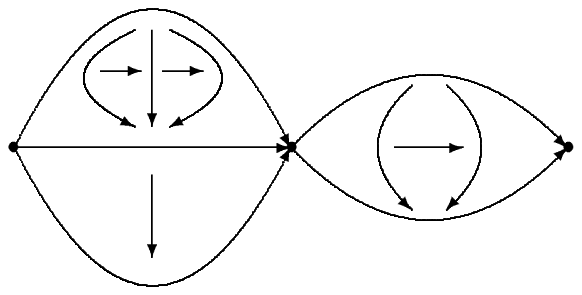}
\]
to a single 3-cell is given exactly as a way of composing up its
2-dimensional source to a single 2-cell and a way of composing its
2-dimensional target to a single 2-cell. In particular, this holds when the
3-cell diagram is a 
\index{degeneracy!of globular diagram}%
\index{composition!horizontal}%
degenerate one such as 
\[
\gfst{}\gtwo{}{}{}\gfbw{}\gtwo{}{}{}\glst{}\ :
\label{p:horiz-in-gray}
\]
so if 
\[
\gfst{}\gtwo{f}{f'}{\alpha}\gfbw{}\gtwo{g}{g'}{\beta}\glst{}
\]
is a diagram in a Gray-category then there is a designated 3-cell
\[
\gfst{}\gthreecell{gf}{g'f'}{\gamma_1}{\gamma_2}{}\glst{}\ , 
\]
where $\gamma_{1}=\beta f'\of g\alpha$ and $\gamma_{2}=g'\alpha\of\beta f$
(the two horizontal composites mentioned in~\ref{sec:ssq}). Similarly, there
is a designated 3-cell $\gamma_{2}\go\gamma_{1}$; since elements of
\[
\Gy\left(%
\gfst{}\gtwo{}{}{}\gfbw{}\gtwo{}{}{}\glst{}%
\right)
\]
are determined by their source and target, the two 3-cells are mutually
inverse.  Thus $\gamma_{1}$ and $\gamma_{2}$ differ by a (canonical)
invertible 3-cell.

Having given this near-elementary description of what a Gray-category is, we
will feel free to speak of Gray-categories even when the underlying graphs
are 
\index{Gray-category!large}%
large.

\section{\textbf{2-Cat} as a Gray-Category}	\label{sec:tooc-gray}
\index{Gray-category!two@\tooc\ as|(}%
\index{two@\tooc\ as Gray-category|(}

Consider all small 2-categories and all homomorphisms, strong transformations
and modifications between them: what kind of structure do they form? Our
answer is `a Gray-category'. As we shall see in the next section, they also
form a tricategory, but only in a 
\index{non-canonical choice}%
non-canonical way.

Recall that 
\index{homomorphism!of bicategories}%
\index{composition!in Bicat@in \Bicat}%
homomorphisms of bicategories can be composed, and that this
composition obeys strict associativity and identity laws. (For the proof of
this see 
\index{B\'{e}nabou}%
\cite{Ben}; for bicategory terminology see the preliminary chapter.)
Next, strong 
\index{transformation of bicategories}%
\index{composition!of transformations}%
transformations can be composed vertically
(\gfst{}\gthree{}{}{}{}{}\glst{}); in general the associativity and identity
laws for this composition hold only up to invertible modification, but if the
codomain bicategory is a 2-category then they hold strictly. There is 
\index{non-canonical choice}%
no \emph{canonical} horizontal composition of strong transformations, even
when we are only dealing with 2-categories, as we shall discuss
later. However, a strong transformation can be composed canonically with a
homomorphism on either side---
\[
\gfst{}\gone{}\glft{}\gtwo{}{}{}\glst{} 
\diagspace\mr{or}\diagspace
\gfst{}\gtwo{}{}{}\grgt{}\gone{}\glst{}
\]
---and the laws indicated in Figure~\ref{fig:ssq-axioms}
(page~\pageref{fig:ssq-axioms}) then hold. It follows that 2-categories,
homomorphisms and strong transformations form a sesquicategory.

We next add in modifications to make \tooc\ into a Gray-category.
Modifications can be composed in all 
%
%
the plausible ways:
\[
\epsfig{file=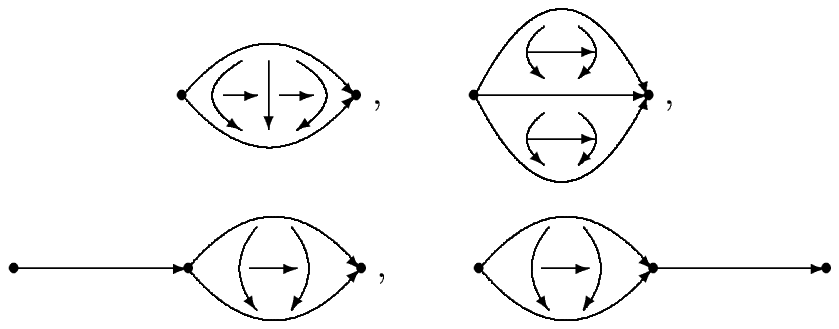}
\]
and these compositions behave as well as is allowed by the compositions of
the transformations that bound them. For instance, there are four obvious
ways to compose a diagram 
\[
\gfst{}\gthreecellu\gfbw{}\gthreecellu\glst{}
\]
into a single modification, as there are two obvious ways to compose the
source transformations 
\[
\gfst{}\gtwo{}{}{}{}{}\gfbw{}\gtwo{}{}{}{}{}\glst{},
\]
and similarly the target. The point is that the data for a modification
consists just of 2-cells, and a diagram of 2-cells in a bicategory 
\index{top dimension principle@`top dimension principle'}%
either commutes or it doesn't---there's no possibility of it `commuting up to
isomorphism'. To supply a proof that \tooc\ forms a Gray-category would
require either an \emph{ad hoc} inductive argument or a 
\index{finite axiomatization}%
finite axiomatization of Gray-categories,%
\footnote{%
Such an axiomatization appears to have been carried out---see \cite[section
2]{Cra}---but I have not investigated this.}
which we could then check \tooc\ against.

We have argued that \tooc\ forms a Gray-category, but in fact it's almost
true that 
\index{Bicat@\Bicat!as near-Gray-category}%
\Bicat\ (=bicategories, homomorphisms, strong transformations and
modifications) forms one too. The only obstacle is that a Gray-category
\cat{C} is required to be locally a 2-category: that is, for each pair of
0-cells, \homset{\cat{C}}{A}{B} must be a 2-category. This holds in \tooc\
because the vertical composition of transformations obeys strict
associativity and identity laws there. But suppose we take the axioms for a
sesquicategory and drop the associativity and identity axioms for vertical
composition of 2-cells: then we obtain a more general notion of
Gray-category, which we call a 
\index{near-Gray-category|emph}%
\emph{near-Gray-category}. Thus \Bicat\ forms a near-Gray-category. The%
\label{p:difference}
difference between near-Gray-categories and Gray-categories may be regarded
as inessential because of the 
\index{coherence!for bicategories}%
coherence theorem for bicategories. By analogy,
\index{Gordon-Power-Street}%
Gordon, Power and Street prove that any tricategory \cat{T} is triequivalent
to one which is locally a 2-category, by replacing each bicategory
\homset{\cat{T}}{A}{B} with a biequivalent 2-category. The next
section shows how (near-)Gray-categories can, in fact, be realised as tricategories.%
\index{Gray-category!two@\tooc\ as|)}%
\index{two@\tooc\ as Gray-category|)}

\section{Gray-Categories as Tricategories}	\label{sec:gray-tri}
\index{Gray-category!as tricategory|(}%
\index{tricategory!Gray-category as|(}%
\index{tricategory|(}%

A tricategory is defined in 
\index{Gordon-Power-Street}
\cite{GPS} in the same style as a bicategory is
defined in 
\index{B\'{e}nabou}
\cite[1.1]{Ben}. That is, a 
\index{tricategory|emph}%
tricategory \tee\ consists of a collection $\tee_0$ of objects, a bicategory
$\homset{\tee}{A}{B}$ for each pair \pr{A}{B} of objects, composition and
identity homomorphisms $\homset{\tee}{B}{C}\times\homset{\tee}{A}{B} \go
\homset{\tee}{A}{C}$ and $\wun\go\homset{\tee}{A}{A}$, and then various
\index{coherence!data}
coherence cells satisfying some axioms. In particular, there is a specified
\index{composition!horizontal}%
horizontal composition of 2-cells. This feature is lacking in
Gray-categories, and we saw on page~\pageref{p:horiz-in-gray} that there are
two symmetrically opposite ways of deriving such a composition, so in order
to turn a Gray-category into a tricategory we need to make a 
\index{non-canonical choice|(}%
non-canonical choice. Let us make this choice once and for all: define the
horizontal composite $\beta*\alpha$ of
\[
\gfst{}\gtwo{f}{f'}{\alpha}\gfbw{}\gtwo{g}{g'}{\beta}\glst{}
\]
as $\beta f'\of g\alpha$, as in diagram~(\ref{pic:derived-horiz}) on
page~\pageref{pic:derived-horiz}. It is routine to verify that every
Gray-category then becomes a tricategory.

\index{contraction!and non-canonical choice}%
Contractions provide a way of discussing this point. To do this we must first
note that this chapter's approach to 
\index{finite dimensionality}%
finite-dimensional structures (page~\pageref{p:fin-dim-attitude}) varies from
that in Chapter~\ref{ch:glob}: there, a 
\index{n-category@$n$-category!lax}%
lax $n$-category was a lax $\omega$-category whose graph was trivial beyond
dimension $n$, whereas here we never even contemplate dimensions beyond $n$.
If $C$ is an $n$-collection, define a
\index{contraction!finite-dimensional|(}%
\emph{contraction} on $C$ to be a function $\psi$ just as before
(page~\pageref{p:defn-contr}) for $\tau$ in dimensions 1 to $n$, together
with the condition that if $\sigma$ is an $n$-stage tree and \pr{f}{f'} a
matching pair for $\sigma$ then $f=f'$. We assume (cf.\ \cite{Bat}) that lax
$n$-categories are the algebras for $K_n$, the initial
$n$-operad-with-contraction. For now all we actually need to know is that
\index{K3@$K_3$ (operad)|emph}%
$K_3$-algebras are `the same as' small tricategories, which the informal
argument of~\ref{sec:infout} makes seem plausible. The phrase `the same as'
must be qualified by the fact that in a $K_3$-algebra there are specified
compositions of arbitrary shape and arity, whereas in a tricategory there are
\index{finite axiomatization}%
just a few, all binary or nullary. But this is just the 
\index{bias}%
biased-unbiased distinction of page~\pageref{sec:bias} and ought to be
harmless.

Now, any contraction on the 3-operad \Gy\ gives an operad map $K_{3} \go
\Gy$, since $K_3$ is initial. This in turn yields a functor $\Alg(\Gy) \go
\Alg(K_{3})$, which on objects sends each Gray-category to a tricategory.
But a contraction on \Gy\ entails a choice of an element of 
\index{taun@$\tau_n$}%
$\Gy(\tau_{n})$ for each $n$, where
\[
\tau_{n} = \underbrace{%
\gfst{}\gtwo{}{}{}\gfbw{}\gtwo{}{}{}\gfbw{}%
\diagspace\cdots\diagspace%
\gfbw{}\gtwo{}{}{}\glst{}}_{n},
\]
and we know there is no canonical way to do that. What we did at the end of
the first paragraph of this section was to choose a particular element of
$\Gy(\tau_{2})$, and that was enough to turn Gray-categories into
tricategories because we were using traditional, 
\index{bias}%
`biased' tricategories.  There are now infinitely many choices, but from our
chosen element of $\Gy(\tau_{2})$ we extract one of the two most obvious
choices and settle on that. For instance, the chosen element of
$\Gy(\tau_{4})$ is obtained by vertical composition according to the diagram
\begin{eqnarray*}
\gfst{}\gtwo{}{}{}\grgt{}\gone{}\gblw{}\gone{}\gblw{}\gone{}\glst{}	\\
\gfst{}\gone{}\glft{}\gtwo{}{}{}\grgt{}\gone{}\gblw{}\gone{}\glst{}	\\
\gfst{}\gone{}\gblw{}\gone{}\glft{}\gtwo{}{}{}\grgt{}\gone{}\glst{}	\\
\gfst{}\gone{}\gblw{}\gone{}\gblw{}\gone{}\glft{}\gtwo{}{}{}\glst{}
\makebox[0pt]{.}
\end{eqnarray*}%
\index{contraction!finite-dimensional|)}%

We have shown that any Gray-category is a tricategory, albeit
non-canonically, so \tooc\ is a tricategory. In fact, everything above holds
just as well for 
\index{near-Gray-category!as tricategory}%
near-Gray-categories, so 
\index{Bicat@\Bicat!as tricategory}%
\Bicat\ is a tricategory too. We will now see exactly why \tooc\ is an
example of a Gray-category which is not canonically a tricategory, as
asserted all along. All this says is that there is no canonical way to define
the 
\index{composition!horizontal!of transformations}%
\index{composition!of transformations}%
\index{transformation of bicategories}%
horizontal composite of strong transformations. So, take a diagram
\[
\bee\ctwo{F}{G}{\sigma}\beep\ctwo{F'}{G'}{\sigma'}\cat{B''}
\]
of 2-categories, homomorphisms and strong transformations. The components of
a putative $\sigma'*\sigma$ are 1-cells
$\gobyc{(\sigma'*\sigma)_{B}}{F'FB}{G'GB}$ in $\bee''$, and the different
routes round the square
\begin{diagram}
F'FB			&\rTo^{F'\sigma_B}	&F'GB			\\
\dTo<{\sigma'_{FB}}	&			&\dTo>{\sigma'_{GB}}	\\
G'FB			&\rTo^{G'\sigma_B}	&G'GB			\\
\end{diagram}
provide the two possibilities. If $\sigma'$ were strict then this square
would commute and we would be home, but as $\sigma'$ is only strong there's
just an isomorphism between the two possibilities. Making our standard choice,
we obtain a tricategory \tooc, and similarly \Bicat. The other possible
tricategory \Bicat\ is triequivalent to this one 
\index{Gordon-Power-Street}%
(\cite[p.\ 7]{GPS}).%
\index{non-canonical choice|)}%

We finally record that when \tee\ is a tricategory, we will use the phrase
`\tee\ is a [near-]Gray-category' to mean that \tee\ is \emph{equal}\/ to the
tricategory arising from some [near-]Gray-category. Amongst all
tricategories, the class of [near-]Gray-categories is closed under 
\index{isomorphism!of tricategories}%
\index{strict!isomorphism}%
\index{isomorphism!strict}%
strict isomorphism (=invertible strict homomorphism) but not plain
isomorphism (=invertible homomorphism), or certain dualities.%
\index{Gray-category!as tricategory|)}%
\index{tricategory!Gray-category as|)}%

\section{Homomorphisms and Triequivalence}	\label{sec:homms-trieqs}

Before moving on to Cayley representation, we need some more definitions.

Objects in a tricategory have the chance of being equal, isomorphic,
equivalent or biequivalent (or, of course, none of these); 1-cells can be
equal, isomorphic or equivalent; 2-cells can be equal or isomorphic; 3-cells
can only be equal or not.  A 
\index{homomorphism!of tricategories|emph}%
\emph{homomorphism of tricategories} is a map of
their underlying graphs such that all kinds of composition are preserved up
to the weakest type of equivalence, and with specified cells providing those
equivalences---e.g.\ if $A\goby{f}B\goby{g}C$ is a diagram of 1-cells in the
domain category, there's a specified equivalence between $Fg\of Ff$ and
$F(g\of f)$. This is usually expressed
\index{Gordon-Power-Street}%
(\cite{GPS}) by saying
that a homomorphism \gobyc{F}{\tee}{\teep} consists of a function
\gobyc{F_0}{\tee_0}{\teep_0} on the object-collections, a homomorphism
\gobyc{F_{AB}}{\homset{\tee}{A}{B}}{\homset{\teep}{F_{0}A}{F_{0}B}} of
bicategories for each $A$ and $B$, and various extra 
\index{coherence!data}%
data and axioms.

We can now say what a 
\index{triequivalence|emph}%
\emph{triequivalence} is: a homomorphism
\gobyc{F}{\tee}{\teep} which is locally a biequivalence and is
surjective-up-to-biequivalence on objects. The first condition means each
$F_{AB}$ is a biequivalence of bicategories, and the second that for each
$A'\elt\teep_0$ there exists $A\elt\tee_0$ such that $FA$ is biequivalent to
$A'$ (in \teep). Unsurprisingly, this definition of a triequivalence amounts
to the same thing as a `pseudo-inverse' definition, as in the 1- and
2-dimensional cases (see page~\pageref{p:bieq} and
\index{Gordon-Power-Street}%
\cite[3.5]{GPS}). In
particular, triequivalence is a symmetric relation.%

\section{The Cayley Characterization of Gray-Categories} 
\label{sec:char-of-gray}
\index{Gray-category!Cayley representation of|(}%
\index{Cayley representation!of Gray-category|(}%

We know \tooc\ is a Gray-category, so any sub-tricategory of \tooc\ is also a
Gray-category.  Here we establish a converse, that every small Gray-category
is strictly isomorphic to a sub-tricategory of \tooc. So the small
Gray-categories are characterized as the small sub-tricategories of \tooc, up
to strict isomorphism. 

To achieve this converse we use the Cayley homomorphism
\gobyc{\Cay}{\stee}{\Bicat}, sending an object $A$ of a small tricategory
$\stee$ to 
$\coprod_{B\elt\stee_0}\homset{\stee}{B}{A}$. We show that if \stee\ is a
Gray-category then
\begin{enumerate}
\item the image of \Cay\ lies in \tooc			\label{tooc}
\item \Cay\ is injective (in each of the 4 dimensions) on the underlying
graphs	\label{inj}
\item \Cay\ is a strict homomorphism.			\label{strict}
\end{enumerate}
From~(\ref{inj}) and~(\ref{strict}) it follows that the image of \Cay\ is
closed under taking identities, composites and coherence cells in \Bicat.
Thus the image of \Cay\ is a sub-tricategory of \Bicat---and in fact of
\tooc, by~(\ref{tooc}). Hence \Cay\ provides a strict isomorphism from \stee\
to a sub-tricategory of \tooc.

We first describe the Cayley homomorphism. Given a 3-cell
\[
A\cthreecell{f}{g}{\alpha}{\twid{\alpha}}{x}A'
\]
in a small tricategory
\stee, there arises a diagram
\[
\coprod_{B\elt\stee_0}\homset{\stee}{B}{A}%
\cthreecell{f_*}{g_*}{\alpha_*}{\twid{\alpha}_*}{x_*}%
\coprod_{B\elt\stee_0}\homset{\stee}{B}{A'}
\]
in \Bicat, which is (by definition) the image of the original 3-cell under
\Cay. Just as for the 
\index{Yoneda embedding!for bicategories}%
Yoneda embedding of a bicategory, the coherence data
for the homomorphism is provided by the coherence data for \stee. For
instance, if $A\goby{f} A' \goby{f'} A''$ in \stee\ then the equivalence
between the homomorphisms $f'_{*}\of f_{*}$ and $(f'f)_*$ in \Bicat\ is
provided by the associativity equivalence in \stee.

Now suppose \stee\ is a Gray-category, and let us establish
claims (\ref{tooc}), (\ref{inj}) and (\ref{strict}).
\begin{enumerate}
\item
\emph{The image of \Cay\ lies in \tooc.}\/
All this says is that \stee\ is locally a 2-category, i.e.\ the vertical
composition of 2-cells obeys strict associativity and identity laws. This is
precisely the difference between a Gray-category and a 
\index{near-Gray-category}%
near-Gray-category.  
\index{Gray-category!as tricategory}%
(Since we are treating \stee\ as a tricategory, rather than a Gray-category
as a structure in its own right, we should not really just say `composition
obeys the associativity law' but that `the associativity isomorphism for
composition is the identity'. The latter statement is part of the definition
of the tricategory arising from a Gray-category, although we didn't say it
explicitly.)
\item
\emph{\Cay\ is injective on the underlying graphs.}\/
\begin{description}
\item[0-cells]
By convention, \homset{\stee}{B}{A} and \homset{\stee}{B'}{A'} are disjoint
unless $B=B'$ and $A=A'$. Thus if $A\neq A'$ then $1_{A}\elt
\Cay(A)\without\Cay(A')$.
\item[1-cells]
Take $A \parpair{f}{g} B$ in \stee. If $\Cay(f)=\Cay(g)$ then
$f1_{A}=g1_{A}$; 1-cells in a Gray-category obey the unit law $h1=h$, so
$f=g$.
\item[2-cells]
Take $A\ctwopar{f}{g}{\alpha}{\twid{\alpha}}B$ in \stee. If
\gobyc{\alpha_{*}=\twid{\alpha}_{*}}{f_*}{g_*} then in particular
\gobyc{\alpha_{*}(1_{A})=\twid{\alpha}_{*}(1_{A})}{f1_{A}}{g1_{A}}, that is,
\gobyc{\alpha 1_{A}=\twid{\alpha}1_{A}}{f}{g}.  But $\alpha 1_{A}=\alpha$ and
$\twid{\alpha} 1_{A}=\twid{\alpha}$ by the nullary version of the second law
in Figure~\ref{fig:ssq-axioms} (page~\pageref{fig:ssq-axioms}).
\item[3-cells]
Take $A\cthreecellpar{}{}{}{}{x}{y}B$ in \stee. If $x_{*}=y_*$ then in
particular $x1_{A}=y1_{A}$. But $x=x1_A$ and $y=y1_A$ by the Gray-category
axioms, with the 
\index{top dimension principle@`top dimension principle'}%
principle in mind that axioms for parallel cells in the top
dimension either hold up to equality or not at all.
\end{description}
\item
\emph{\Cay\ is a strict homomorphism.}\/
This is the step that takes most work, although none of it is difficult.
Strictly speaking, our task is to show that all the cells which make up the
coherence data of the homomorphism are identities, where we are thinking of
the homomorphism as consisting of a graph map plus coherence data. The force
of this is that \Cay\ strictly preserves all kinds of composition and
identities. We give two instances of this:
\begin{description}
\item[Binary composition of 1-cells] 
Take $A\goby{f}A'\goby{f'}A''$, giving
\[
\gobyc{f'_{*}\of f_{*},(f'f)_{*}}{\dcayset{A}}{\dcayset{A''}}.
\]
On a 0-cell $B\goby{p}A$ of \cayset{A}, we have $(f'_{*}\of f_{*})(p) =
f'(fp)$ and $(f'f)_{*}(p)=(f'f)p$, and these are equal because 1-cell
composition is associative. On a 1-cell $B\ctwo{p}{q}{\gamma}A$ of
\cayset{A}, we have $(f'_{*}\of f_{*})(\gamma) = f'(f\gamma)$ and
$(f'f)_{*}(\gamma)=(f'f)\gamma$, which are equal by the right-handed version
of the second diagram in Figure~\ref{fig:ssq-axioms}
(page~\pageref{fig:ssq-axioms}). A 2-cell $x$ of \cayset{A} is a 3-cell of
\stee, and the 3-cells $f'(fx)$, $(f'f)x$ are equal by the Gray-category
axioms 
\index{top dimension principle@`top dimension principle'}%
(`top dimension principle').

We have thus shown that the graph maps underlying the homomorphisms
$f'_{*}\of f_{*}$, $(f'f)_{*}$ are equal. Having also established that they
have the same coherence data (certain 3-cells), we would be able to conclude
that $f'_{*}\of f_{*} = (f'f)_{*}$.
\item[Horizontal composition of 2-cells]
\index{composition!horizontal}%
This part is slightly different because a Gray-category does not naturally
carry a horizontal composition. Take
\[
A\ctwo{f}{g}{\alpha}A'\ctwo{f'}{g'}{\alpha'}A''
\]
in \stee, yielding
\[
\coprod\homset{\stee}{B}{A}%
\ctwo{f_*}{g_*}{\alpha_*}%
\coprod\homset{\stee}{B}{A'}%
\ctwo{{f'}_*}{{g'}_*}{{\alpha'}_*}%
\coprod\homset{\stee}{B}{A''}
\]
in \Bicat. We have to show that
\[
\gobyc{\alpha'_{*} * \alpha_{*}}{f'_{*}\of f_{*}}{g'_{*}\of g_{*}}
\]
and
\[
\gobyc{(\alpha' * \alpha)_{*}}{(f'f)_{*}}{(g'g)_{*}},
\]
strong transformations of bicategories, are equal. (We already know that their
domains and codomains are equal.) Since \Bicat\ is a Gray-category,
$\alpha'_{*} * \alpha_{*}$ is by definition equal to
\[
f'_{*}\of f_{*} \goby{f'_{*}\alpha} f'_{*}\of g_{*}  
\goby{\alpha'_{*}g} g'_{*}\of g_{*}.
\]
Since \stee\ is a Gray-category, $\alpha' * \alpha$ is defined as
\[
f'f \goby{f'\alpha} f'g \goby{\alpha' g} g'g,
\]
so $(\alpha' * \alpha)_{*} = ((\alpha'g)\of(f'\alpha))_{*}$. The substance of
the assertion is that the components of $\alpha'_{*} * \alpha_{*}$ and
$(\alpha' * \alpha)_{*}$ at a 0-cell $B\goby{p} A$ of \cayset{A} are
equal. This says that the 2-cells $(\alpha'(gp))\of(f'(\alpha p))$ and
$((\alpha'g)p)\of((f'\alpha)p)$ of \stee\ are equal, which follows from the
Gray-category axioms.
\end{description}
\end{enumerate}
We have now sketched out a proof that up to strict isomorphism, the small
Gray-categories are the small sub-tricategories of \tooc. If \stee\ is in
fact a (strict) 
\index{three-category@3-category!Cayley representation of}%
\index{Cayley representation!of 3-category}%
3-category then the image of the Cayley map lies in the
3-category 
\index{two-cat-strict@\fcat{2\hyph Cat_{strict}}}%
\fcat{2\hyph Cat_{strict}} of 2-categories, strict homomorphisms,
strict transformations, and modifications: thus up to strict isomorphism, the
small 3-categories are the small sub-tricategories of
\fcat{2\hyph Cat_{strict}}.%
\index{Gray-category!Cayley representation of|)}%

\section{A Conjecture on Coherence in Higher Dimensions}
\label{sec:hd-coh-conj}
\index{coherence!conjecture|(}%

Several points arise from the proof in the previous section. Firstly, we can
go through the same arguments in lower dimensions. For small categories
\scat{C}, the 
\index{Cayley representation!of category}%
Cayley functor $\scat{C}\go\Set$ is always injective
(=injective in both dimensions) and so an isomorphism to its image. For small
bicategories \scat{B}, the homomorphism 
\index{Cayley representation!of bicategory}%
$\gobyc{\Cay}{\scat{B}}{\Cat}$ is
injective if $f\of1=f$ for 1-cells $f$, and a strict homomorphism if $1\of
f=f$ and $(h\of g)\of f = h\of(g\of f)$. Thus \Cay\ is an injective strict
homomorphism if \scat{B} is a 2-category, and up to strict isomorphism, the
small 2-categories are the small sub-bicategories of \Cat.

Secondly, because we assume \stee\ is a Gray-category, the 
\index{image}%
image of
\gobyc{\Cay}{\stee}{\Bicat} is already a tricategory---it doesn't need closing
under composition or enlarging in any other way. Thus \Cay\ provides a
triequivalence of \stee\ with its \emph{genuine}\/ image. For an arbitrary
tricategory \stee, the image of \Cay\ might not be a tricategory. One
possibility is that if we take a suitably enlarged `image' (e.g.\ the full
image) then \Cay\ is always a triequivalence to its `image'; as remarked in
\index{Gordon-Power-Street}%
\cite[p.~4]{GPS},
`this may be true, but a proof is not so easy'. (In fact, they were not
speaking of \Cay\ but the 
\index{Yoneda embedding!for tricategories}%
Yoneda homomorphism \gobyc{Y}{\stee}{\ftrcat{\stee^{\op}}{\Bicat}}. Many of
the arguments above apply just as well to $Y$ as to \Cay, and indeed
$\Cay=\coprod\of Y$.)

Thirdly, let us assume that for all small tricategories \stee, the Cayley map
\emph{does}\/ provide a triequivalence to some sub-tricategory of \Bicat.
Then we have a proof of the 
\index{coherence!for tricategories}%
coherence theorem for tricategories in the small
case---that is, that any small tricategory is triequivalent to a
Gray-category---since the 
\index{coherence!for bicategories}%
coherence theorem for bicategories implies that the inclusion
$\tooc\go\Bicat$ is a triequivalence. (Again, we may be able to work with $Y$
instead of \Cay\ and get a proof for tricategories which are just
\emph{locally}\/ small.)

Finally, we conjecture a coherence theorem for higher-dimensional categories.
Let \scat{R} be a small 
\index{tetracategory}%
tetracategory; assume both that \fcat{Tricat} forms a
tetracategory (perhaps after making some 
\index{non-canonical choice}%
non-canonical choices), and that
there is a Cayley homomorphism \gobyc{\Cay}{\scat{R}}{\fcat{Tricat}} which
provides a tetraequivalence to some sub-tetracategory of \fcat{Tricat}.
By the coherence theorem for tricategories, the inclusion $\CatOF{Gray}
\go \fcat{Tricat}$ is a tetraequivalence. Thus \scat{R} is tetraequivalent to
a sub-tetracategory \cat{R} of \CatOF{Gray}. The image of \scat{R} in
the tetraequivalence $\scat{R}\go\cat{R}$ is small, so the
sub-tetracategory \scat{R'} of \cat{R} it generates is small, and we then
have a tetraequivalence from \scat{R} to the small sub-tetracategory
\scat{R'} of \CatOF{Gray}.
\begin{defn}
\index{Grayn-category@Gray\bkpr{n}-category|emph}%
\index{n-category@$n$-category!lax}%
A \emph{small Gray\bkpr{0}-category} is a set. A \emph{small
Gray\bkpr{n+1}-category} is a small lax $(n+1)$-category strictly isomorphic to
a sub-lax-$(n+1)$-category of \CatOF{Gray\bkpr{n}}.
\end{defn}
Here, \fcat{Gray\bkpr{n}\hyph\Cat} denotes a full sub-lax-$(n+1)$-category of
\fcat{Lax}-$n$-\Cat, which we are assuming has been defined in a meaningful
way. For the first few values of $n$:
\begin{itemize}
\item
A small Gray\bkpr{1}-category is a small category isomorphic to a subcategory
of \Set, i.e.\ a small category.
\item
A small Gray\bkpr{2}-category is a small bicategory strictly isomorphic to a
sub-bicategory of \Cat, i.e.\ a small 2-category. 
\item
A small Gray\bkpr{3}-category is a small tricategory strictly isomorphic to a
sub-tricategory of \tooc, i.e.\ a small Gray-category.
\end{itemize}

The argument on tetracategories above seeks to show that any small
tetracategory is tetraequivalent to a small Gray\bkpr{4}-category. In
dimension 2, it is the usual proof of the coherence theorem but using Cayley
rather than Yoneda; in dimension 3, it is our conjectured proof of coherence
for tricategories (`thirdly' above). The same argument, with the same
assumptions as for tetracategories, may be repeated in dimensions 5, 6,
\ldots. We therefore guess:
\begin{conjecture}
Any small lax $n$-category is laxly $n$-equivalent to a small
Gray\bkpr{n}-category.
\end{conjecture}
With luck, we may be able to discard our 
\index{smallness}%
`small's. The issue is presumably a
distraction; this is especially the case if we are just using a coherence
theorem to ease our manipulation of diagrams, as when we ignore the
distinction between $(A\otimes B)\otimes C$ and $A\otimes(B\otimes C)$ in a
lax 
\index{monoidal category!coherence in}%
\index{coherence!for monoidal categories}%
monoidal category.%
\index{Gray-category|)}
\index{coherence!conjecture|)}%
\index{Cayley representation!of Gray-category|)}%
\index{tricategory|)}%

%% file: essopap.tex
\chapter{The Opetopic Approach}	\label{ch:opap}
\index{opetopes|(}%
\index{pasting diagrams|(}%

Opetopes are the backbone of the 
\index{Baez!-Dolan}%
Baez-Dolan school of higher-dimensional
category theory 
\index{Baez}%
(\cite{Baez}, \cite{BD}, 
\index{Hyland}%
\cite{Hy}, 
\index{Hermida!-Makkai-Power}%
\cite{HMP}). In a
nutshell, there is one 0-opetope and an $(n+1)$-opetope is a pasting of
$n$-opetopes. Of course, it is the meaning of `a pasting' that matters, and
we explain this now.

The unique 0-opetope (with `ope' pronounced as in `operation') is drawn as
\gzero{}. The only pasting-together of \gzero{}'s (`0-pasting diagram') is
\gzero{} itself, so there's just one 1-opetope, which we like to draw as
\gfst{}\topebase{}\glst{}. Now $k$ copies of \gfst{}\topebase{}\glst{} can be
pasted together end-to-end for any $k\elt\nat$, as in
\[
\gfst{}\topebase{}\gblw{}\topebase{}\gblw{}\topebase{}\glst{}
\diagspace
\mr{and}
\diagspace
\gzero{}
\]
($k=3$ and $k=0$). The set of 2-opetopes is therefore isomorphic to \nat.
However, in our pictures we distinguish between 1-pasting diagrams (as drawn
above) and 2-opetopes, even though there's a natural one-to-one
correspondence. For instance, the two 2-opetopes corresponding to the two
1-pasting diagrams above are drawn as
\[
\topec{}{}{}{}{\Downarrow}
\diagspace
\mr{and}
\diagspace
\topez{}{\Downarrow}.
\]
A typical 2-pasting diagram looks like
\[
\epsfig{file=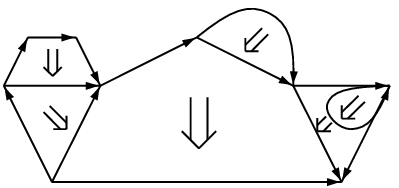}.
\]
Note that the arrows go in compatible directions. One can think of the top
edges of a 2-opetope as being inputs and the bottom edge as an output; in
forming a pasting diagram the rule is that input edges paste to output edges
(compare diagram~(\ref{diag:comptrans}) on page~\pageref{diag:comptrans}). A
3-opetope is formally just a 2-pasting diagram. The best picture we can do on
a sheet of paper is
\[
\epsfig{file=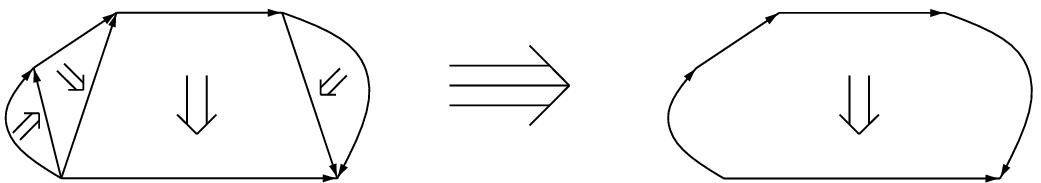},
\]
where the 2-opetope on the right is the boundary of the pasting diagram on
the left: we imagine a 3-dimensional figure with a flat bottom face and four
curved top faces. And so the process continues.

Formally, the sets $S_n$ of $n$-opetopes admit a definition in terms of free
multicategories; this is section~\ref{sec:opetopes}. In doing this we obtain
a sequence \Cpn{n} of cartesian monads. By looking at \Cpn{n}-structured
categories, we see (\ref{sec:pds}) that for any fixed $n$ there is a category
of $n$-pasting diagrams. An arrow in this category can be thought of as the
process of `composing together' some of the adjacent $n$-opetopes in the
domain, to obtain the codomain: e.g.\ the codomain may be the domain with one
of its internal edges erased. When $n=2$ this gives a category of 
\index{tree!Baez-Dolan-style!category of}%
trees (\ref{sec:cat-trees}).

The final section,~\ref{sec:slicing}, is on a different tack. 
\index{Baez!-Dolan}%
Baez and Dolan
place great emphasis on the process of \emph{slicing}, in which any operad
$C$ gives rise to its slice operad $C^+$. Although their operads are not
obviously part of our general scheme, we describe an analogous process for
\Cartpr-multicategories.

We concentrate on describing the structures the opetopic approach gives rise
to, rather than its overall shape or its place in $n$-category theory. A
discussion of the latter can be found in 
\index{Baez}%
\cite{Baez}; the ideas we do set out
below have much in common with 
\index{Hermida}%
\cite{Her}.

\section{Opetopes}	\label{sec:opetopes}

In this section we construct for each $n\elt\nat$ the set 
\index{S-n@$S_n$|emph}%
\index{opetopes|emph}%
$S_n$ of $n$-opetopes, and a cartesian monad 
\index{T-n@$T_n$|emph}%
\index{opetopic monads|emph}%
$T_n$ on $\Set/S_n$.

Start with $S_{0}=1$ and $T_{0}=\id$, the identity monad on
$\Set\iso\Set/S_0$; note that \Cpn{0} satisfies the conditions
of~\ref{sec:free-multi} for the formation of 
\index{multicategory!free}%
free \Cpn{0}-multicategories. Now suppose, inductively, that $S_n$ and a
cartesian monad $T_n$ on $\Set/S_n$ are constructed, satisfying the
conditions of~\ref{sec:free-multi}. The terminal object of $\Set/S_n$ is
\slob{S_n}{1}{S_n}, so if we put
\[
\vslob{S_{n+1}}{}{S_n} = T_{n}\bktdvslob{S_n}{1}{S_n}
\]
then the category of \Cpn{n}-graphs on 1 is
\[
\frac{\Set/S_n}{S_{n+1}\go S_n} \iso \frac{\Set}{S_{n+1}}.
\]
Then define $T_{n+1}$ to be the 
\index{operad!free}%
free \Cpn{n}-operad monad on
$\Set/S_{n+1}$. By~\ref{sec:free-multi}, \Cpn{n+1} is
cartesian and satisfies the conditions.

Let us look at the first few values of $n$. We have $S_{0}=1$ and
$T_{0}=\id$; the member of $S_{0}$ may be depicted as the 0-dimensional
diagram \gzero{}. Then $\bktdslob{S_1}{}{S_0} = T_{0}\bktdslob{S_0}{}{S_0}$, so
$S_{1} = 1$, and the member of $S_1$ may be depicted as
\gfst{}\topebase{}\glst{}. The monad $T_1$ is `free \Cpn{0}-operad', that is,
\index{monoid!free}%
`free monoid'. So $T_1$ sends a set $X$ to the set of all diagrams
\[
\gfst{}\topebase{x_1}\gblw{}\topebase{x_2}\gblw{}
\diagspace\cdots\diagspace
\gblw{}\topebase{x_n}\glst{}
\]
with $n\elt\nat$ and $x_{i}\elt X$. Next, $\bktdslob{S_2}{}{S_1} =
T_{1}\bktdslob{S_1}{}{S_1}$, so $S_{2}=\nat$, the free monoid on 1. We will
draw members of $S_2$ as 2-opetopes, e.g.\ \toped{}{}{}{}{}{\Downarrow} for
$4\elt\nat$ or \topez{}{\Downarrow} for $0\elt\nat$. (Sometimes we omit the
arrows.) A \Cpn{1}-operad is a plain operad, so $T_2$ is the monad
\index{operad!plain}%
`free plain operad' on $\Set/S_{2} = \Set/\nat$. To see what $T_2$ does, let
$A = (A(n))_{n\elt\nat}$ be an object of $\Set/S_{2}$. Then $A$ consists of a
set for each 2-opetope, which can be thought of as a set of labels: thus
$a\elt A(n)$ is drawn as
\[
\topeqs{}{}{}{}{a}.
\]
$T_2$ sends $A$ to the family of pictures obtained by sticking together
members of the $A(n)$'s, so if $T_{2}\bktdslob{A}{}{\nat} =
\bktdslob{A'}{}{\nat}$ then a typical member of $A'(4)$ is 
\[
\epsfig{file=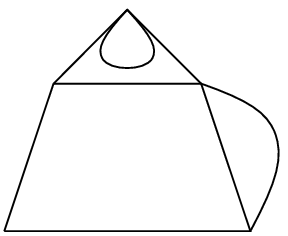}
\]
\setlength{\unitlength}{1cm}
\raisebox{0.5cm}[0cm][0cm]{%
\begin{picture}(12,2.5)
\cell{5.8}{0.9}{c}{a_1}
\cell{6.2}{1.8}{c}{a_2}
\cell{5.8}{2.1}{c}{a_3}
\cell{7.1}{1.0}{c}{a_4}
\end{picture}}
($a_{1}\elt A(3)$, $a_{2}\elt A(3)$, $a_{3}\elt A(0)$, $a_{4}\elt A(1)$).

The last paragraph has partially revealed what \Cpn{n}-operads
look like; how about 
\index{multicategory!for \Cpn{n}}%
\Cpn{n}-multicategories? A \Cpn{0}-multicategory is a category. A
\Cpn{1}-multicategory is a plain multicategory $C$, in which objects
look like \gfst{}\topebases{t}\glst{} ($t\elt C_0$) and arrows like
\[
\begin{array}{c}
\gfst{}\topebases{t_1}\gblw{}\topebases{t_2}\gblw{}\topebases{t_3}\glst{}\\
\ \ \downarrow\ a\\
\gfst{}\topebases{t}\glst{}.
\end{array}
\]
In a \Cpn{2}-multicategory,
\label{p:Cpn-two-multicats}
a typical object looks like \topecs{}{}{}{}{t} (and
this lies over the element 3 of $\nat=S_2$), and a typical arrow looks like
\begin{equation}		\label{pic:typical-arrow}
\raisebox{-21pt}{\epsfig{file=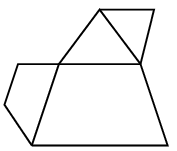}}
\diagspace\goby{a}\diagspace
\raisebox{-21pt}{\epsfig{file=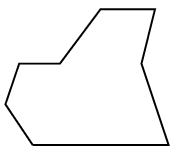}}
\ .
\end{equation}
\setlength{\unitlength}{1cm}
\raisebox{0.5cm}[0cm][0cm]{%
\begin{picture}(12,2.5)
\cell{3.9}{0.3}{c}{t_1}
\cell{3.2}{0.4}{c}{t_2}
\cell{3.9}{1.0}{c}{t_3}
\cell{4.25}{1.1}{c}{t_4}
\cell{8}{0.5}{c}{t}
\end{picture}}
(Here, the opetope on the right-hand side has been drawn to resemble the
boundary of the left-hand side, but of course it is just the opetope
$7\elt\nat = S_2$.) Composition in a \Cpn{2}-multicategory looks like
\begin{equation}
\epsfig{file=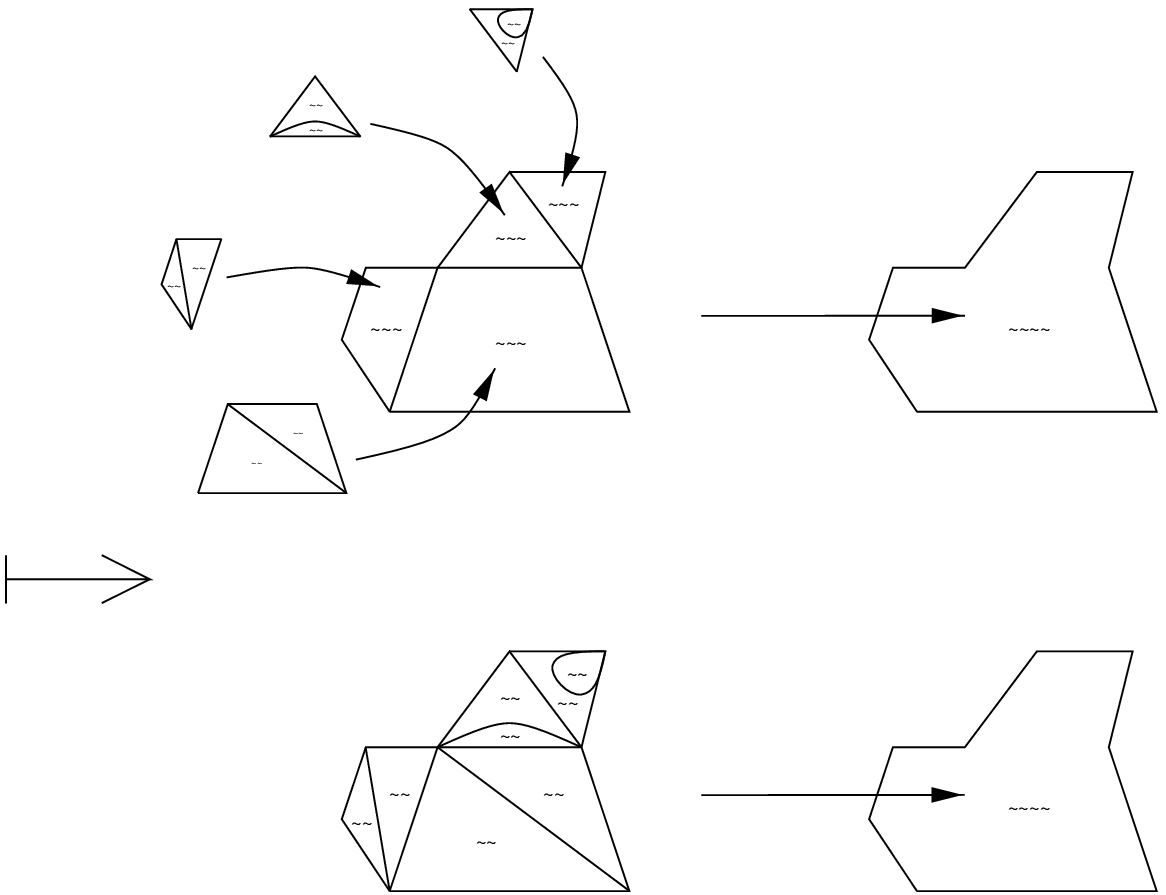}		\label{pic:typical-comp}
\end{equation}
and the identities consist of an arrow
\[
\topeqs{}{}{}{}{t}\diagspace\goby{1}\diagspace\topeqs{}{}{}{}{t}
\]
for each $t$.

Continuing into higher dimensions, in a \Cpn{n}-multicategory the objects are
labelled $n$-opetopes, the domain of an arrow is a labelled $n$-pasting
diagram, and the codomain is an object whose underlying $n$-opetope is the
boundary of the domain. This is illustrated in
diagram~(\ref{pic:typical-arrow}); one may also think of an arrow as an
$(n+1)$-opetope with its $n$- and $(n+1)$-dimensional parts labelled. The
next two sections examine
\Cpn{n}-structured categories and show how, amongst other things,
consideration of them leads naturally to a category of trees.

\section{Pasting Diagrams}		\label{sec:pds}
\index{structured category!for \Cpn{n}|(}%

A strict 
\index{monoidal category}%
monoidal category can be thought of as a 
\index{multicategory!plain}%
plain multicategory in which each sequence $\range{t_1}{t_k}$ of objects has
a representing object $t_{1}\otimes\cdots\otimes t_{k}$: that is, there's a
natural correspondence
\[
\Hom(\range{t_1}{t_k};t) \iso \Hom(t_{1}\otimes\cdots\otimes t_{k};t).
\]
In pictures, a strict monoidal category consists of a category with objects
\gfst{}\topebases{t}\glst{} and arrows
\[
\gfst{}\topebases{t'}\glst{}
\diagspace\goby{a}\diagspace
\gfst{}\topebases{t}\glst{}
\]
together with a monoidal structure
\[
\gfst{}\topebases{t_1}\gblw{}\topebases{t_2}\ \cdots\ \topebases{t_k}\glst{}
\diagspace\goesto\diagspace
\gfst{}\topebases{t_{1}\otimes\cdots\otimes t_{k}}\glst{}
\]
(and similarly for arrows), all obeying the familiar rules. Recall
(\ref{sec:struc}) that a strict monoidal category is a \Cpn{1}-structured
category. One dimension up, a \Cpn{2}-structured category consists of a
category with objects \topeqs{}{}{}{}{t}, arrows $a$ like
\[
\topecs{}{}{}{}{t'}
\diagspace\goby{a}\diagspace
\topecs{}{}{}{}{t},
\]
and an algebraic structure, `tensor', making assignments like
\[
\raisebox{-18pt}{\epsfig{file=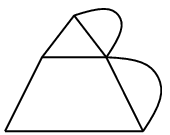}}
\diagspace\goesto\diagspace
\raisebox{-18pt}{\epsfig{file=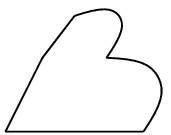}}
\]
\setlength{\unitlength}{1cm}
\raisebox{0.5cm}[0cm][0cm]{%
\begin{picture}(12,2.5)
\cell{4.0}{0.3}{c}{t_1}
\cell{4.0}{0.85}{c}{t_2}
\cell{4.3}{1.0}{c}{t_3}
\cell{4.6}{0.4}{c}{t_4}
\cell{7.95}{0.4}{c}{t_{1}\otimes t_{2}\otimes}
\cell{8.0}{0.1}{c}{t_{3}\otimes t_{4}}
\end{picture}}
and similarly for arrows. (Note the special case of 
\index{nullary!tensoring}%
nullary tensoring,
\[
\gfst{}\topebases{}\glst{}
\diagspace\goesto\diagspace
\topeas{}{}{I}
\diagspace .)
\]
It is hoped that the pattern for $n=3,4,\ldots$ is clear, if difficult to draw.

Formally, a \Cpn{n}-structured category is a category object in
\[
(\Set/S_{n})^{T_n} = \Operadof{\Cpn{n-1}}.
\]
For $n=1$, a \Cpn{1}-structured
category is a category object in 
\index{monoid}%
$\fcat{Monoid}=\Operadof{\pr{\Set}{\id}}$:
that is, a 
\index{monoidal category!strict}%
strict monoidal category. A \Cpn{2}-structured category is a
category object in the category of 
\index{operad!plain}%
plain operads. How does this square with
the description in the last paragraph? What we said there was, effectively,
that a \Cpn{2}-structured category consists of a sequence
$(C(n))_{n\elt\nat}$ of categories, together with functors
\begin{eqnarray*}
C(n_{1}) \times\cdots\times C(n_{k}) \times C(k) & \go & 
C(n_{1} + \cdots + n_{k}) \\
\wun & \go & C(1) 	\label{p:t-two-struc-cat-enr}
\end{eqnarray*}
obeying associativity and identity laws. (So it is just a 
\index{enrichment!of operad in Cat@of operad in \Cat}%
`\Cat-enriched operad', or 
\index{strict!monoidal 2-operad}%
`strict monoidal 2-operad', as in~\ref{eg:multicats}(\ref{eg:soib-op}).) The
object-sets $C(n)_{0}$ of each category $C(n)$ thus form a plain operad
$C_0$, as do the arrow-sets to form $C_1$. The domain, codomain, composition
and identity functions in the $C(n)$'s provide \gph{C_1}{C_0} with the
structure of a category object in the category of plain operads.

To form the 
\index{structured category!free}%
free \Cpn{n}-structured category on a given
\Cpn{n}-multicategory, we paste together objects and arrows of the
multicategory. For instance, if $n=1$, \gfst{}\topebases{t_i}\glst{} are
objects in the given multicategory, and
\[
\begin{array}{ccc}
\gfst{}\topebases{t_1}\gblw{}\topebases{t_2}\gblw{}\topebases{t_3}\glst{}&
\gzero{}&
\gfst{}\topebases{t_4}\gblw{}\topebases{t_5}\glst{}\\
\ \ \ \ \downarrow a_{1}&\ \ \ \ \downarrow a_{2}&\ \ \ \ \downarrow a_{3}\\
\gfst{}\topebases{t_6}\glst{}&
\gfst{}\topebases{t_7}\glst{}&
\gfst{}\topebases{t_8}\glst{}
\end{array}
\]
are arrows, then
\[
\begin{array}{c}
\gfst{}\topebases{t_1}\gblw{}\topebases{t_2}\gblw{}\topebases{t_3}%
\gblw{}\topebases{t_4}\gblw{}\topebases{t_5}\glst{}\\
\ \ \ \ \ \ \ \ \ \ \ \ \ \ \ \ \ \ \ \ \ \ \ \ \ \ \ \ \ \ \ \ \ \ \ \ \ \ \
\ 
\downarrow 
\gfst{}\topebases{a_1}\gblw{}\topebases{a_2}\gblw{}\topebases{a_3}\glst{}\\
\gfst{}\topebases{t_6}\gblw{}\topebases{t_7}\gblw{}\topebases{t_8}\glst{}
\end{array}
\]
is an arrow in the free structured category. (This is
diagram~(\ref{diag:arrows-in-mon-cat}), page~\pageref{diag:arrows-in-mon-cat},
drawn in a different way.) In particular, the free \Cpn{1}-structured
category on \wun\ is the
\index{simplicial category}%
simplicial category $\Delta$. In dimension 2, the terminal
\Cpn{2}-multicategory \wun\ has objects like
\topecs{}{}{}{}{\Downarrow} and arrows like
\[
\raisebox{-13pt}{\epsfig{file=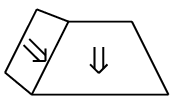}}
\diagspace\go\diagspace
\raisebox{-13pt}{\epsfig{file=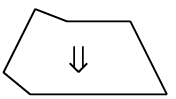}}
\]
and
\[
\gfst{}\topebases{}\glst{}
\diagspace\go\diagspace
\topeas{}{}{\Downarrow}
\]
(see the description of \Cpn{2}-multicategories on
page~\pageref{p:Cpn-two-multicats}). Thus an object of the free
\Cpn{2}-structured category on \wun\ is a 2-pasting diagram, like
\raisebox{-18pt}{\epsfig{file=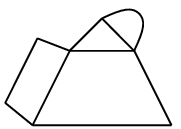}},
and an arrow looks like
\[
\raisebox{-18pt}{\epsfig{file=pd69a.eps}}
\diagspace\go\diagspace
\raisebox{-18pt}{\epsfig{file=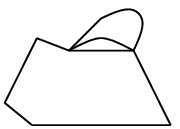}}.
\]
Roughly speaking, an arrow is a removal of some 
\index{internal edge!removal of}%
internal edges. It's only in the 
\index{composition!nullary}%
nullary case that this becomes inaccurate;
it is more precise to say that an arrow is the replacement of some
sub-2-pasting diagrams with their bounding 2-opetopes. So here,
\[
\raisebox{-13pt}{\epsfig{file=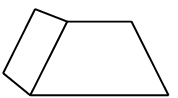}}
,\diagspace
\raisebox{-1pt}{\epsfig{file=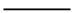}}
\diagspace\mr{and}\diagspace
\raisebox{-7pt}{\epsfig{file=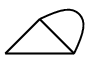}}
\]
are replaced with their bounding 2-opetopes
\[
\raisebox{-13pt}{\epsfig{file=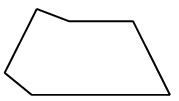}}
,\diagspace
\raisebox{-3pt}{\epsfig{file=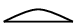}}
\diagspace\mr{and}\diagspace
\raisebox{-7pt}{\epsfig{file=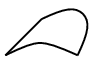}}
\ .
\]

We now make the general definition:
\begin{defn}
\index{structured category!of pasting diagrams|emph}%
\index{pasting diagrams!structured category of|emph}%
\index{PD-n@\PD{n}|emph}%
The \emph{structured category of $n$-pasting diagrams, \PD{n}}, is the free
\Cpn{n}-structured category on the terminal \Cpn{n}-multicategory.
\end{defn}
If $\ess_1$ and $\ess_2$ are categories with pullbacks and $\ess_{1} \go
\ess_{2}$ a functor preserving pullbacks, then there's an induced functor
from $\Catof{\ess_{1}}$ to $\Catof{\ess_{2}}$ (where $\Catof{\ess}$ is the
category of 
\index{internal category}%
internal categories in \ess). Applying this to the forgetful functors
\[
(\Set/S_{n})^{T_n} \go \Set/S_{n} \go \Set
\]
gives a functor
\[
\Strucof{\Cpn{n}} \go \Cat.
\]
The 
\index{pasting diagrams!category of|emph}%
\index{pd-n@\pd{n}|emph}%
\emph{category of $n$-pasting diagrams}, \pd{n}, is defined to be the
image of \PD{n} under this functor. The definition of free structured
category in~\ref{sec:struc} tells us that the objects-object of \PD{n}
is \vslob{S_{n+1}}{}{S_n}, so the object-set of \pd{n} is $S_{n+1}$, the set
of $(n+1)$-opetopes or, as we prefer to think of them at the moment,
$n$-pasting diagrams. An arrow in \pd{n} consists of the removal of some
$(n-1)$-dimensional interior faces, with a warning about the 
\index{composition!nullary}%
nullary case as given for $n=2$.%
\index{structured category!for \Cpn{n}|)}%

\section{The Category of Trees}	\label{sec:cat-trees}
\index{tree!Baez-Dolan-style!category of|(}%

We are familiar with the categories $\pd{0}=\wun$ and 
\index{simplicial category}%
$\pd{1}=\Delta$. Here we take a closer look at \pd{2}.

\begin{figure}
\epsfig{file=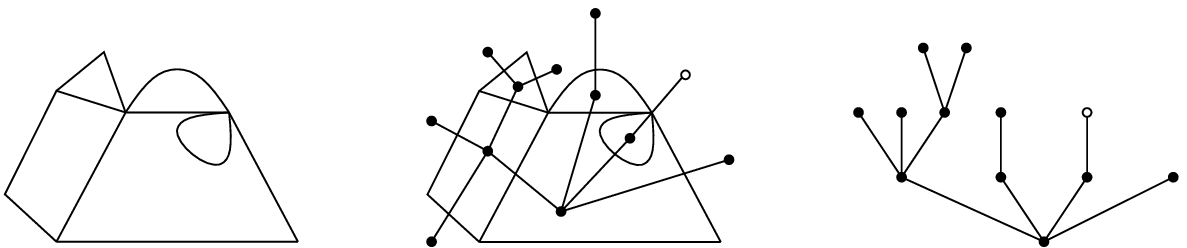}
\caption{How a 2-pasting diagram corresponds to a tree}
\label{fig:tree-is-two-pd}
\end{figure}

As suggested by Figure~\ref{fig:tree-is-two-pd}, there is a one-to-one
correspondence between 2-pasting diagrams and trees. 
(The trivial or 
\index{nullary!pasting diagram}%
nullary 2-pasting diagram, \gfst{}\topebase{}\glst{}, corresponds to the
tree \node. The present trees are not the Batanin trees of
Chapter~\ref{ch:glob}.) We may
\emph{define} a tree to be a 2-pasting diagram: that is, a member of $S_3$.
In turn, $S_3$ is defined by
\[
\vslob{S_3}{}{S_2} = T_{2}\bktdvslob{S_2}{1}{S_2},
\]
which just means that \slob{S_3}{}{\nat} is the graph of the free plain
operad on \slob{\nat}{1}{\nat}. But in this case, the free multicategory
construction sketched in~\ref{sec:free-multi} reduces to the inductive
description of the set $1^*$ of unlabelled trees, given in
Example~\ref{eg:cart-mnds}(\ref{eg:tree-mnd}): a tree is either \node\ or a
sequence of trees. We can thus prove that the two sets of trees are
isomorphic.

\begin{figure}
\vspace{5mm}
\epsfig{file=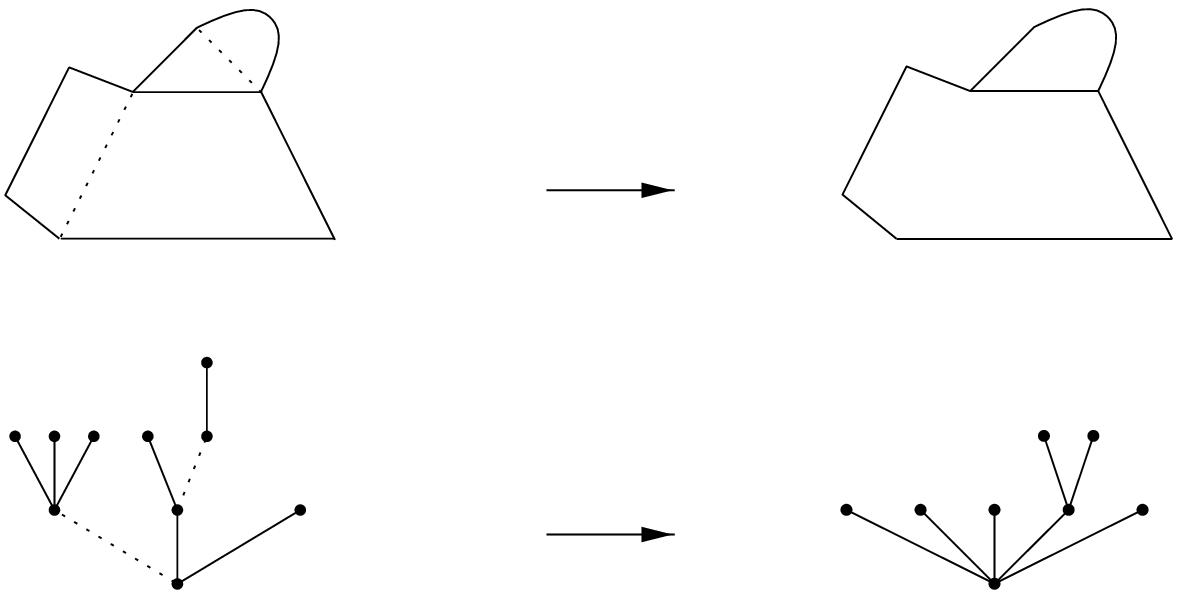}
\caption{Dotted edges are those to be removed/contracted}
\label{fig:edge-removal}
\end{figure}

\begin{figure}
\epsfig{file=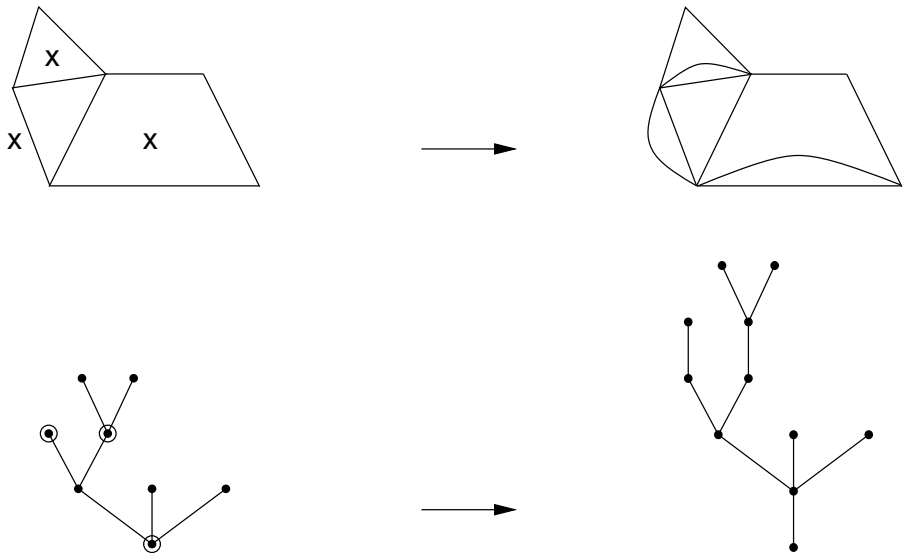}
\caption{Marked regions/nodes are those to have an edge added}
\label{fig:node-to-edge}
\end{figure}

The objects of \pd{2} are trees, so we define the 
\index{tree!Baez-Dolan-style!category of|emph}%
\emph{category of trees} to
be \pd{2}. Loosely, a 
\index{tree!Baez-Dolan-style!map of}%
morphism of trees consists of the 
\index{internal edge!removal of}%
contraction of some internal edges---see
Figure~\ref{fig:edge-removal}---where
an edge is called 
\index{internal edge|emph}%
\emph{internal} if it does not end in a leaf. Once again, saying this ignores
the 
\index{composition!nullary}%
nullary case (Figure~\ref{fig:node-to-edge}).
Replacing a 2-pasting diagram by its bounding 2-opetope corresponds to
replacing an $n$-leafed tree by the $n$-leafed tree $\treeq$ of height 1, so
strictly speaking, a map of trees consists of the replacement of some
subtrees by their corresponding height-1 trees. In particular, a subtree
which is just a node may be replaced by an edge, $\treea$.

Categories of trees have been used by Borcherds to define 
\index{relaxed multilinear category}%
relaxed multilinear categories, by
\index{Soibelman}%
Soibelman for his very similar 
\index{pseudo-monoidal category}%
pseudo-monoidal categories, and by
\index{Kontsevich-Manin}%
Kontsevich and Manin for purposes unknown to me 
\index{Borcherds}%
(\cite{Bor}, 
\index{Soibelman}%
\cite{Soi},
\cite{KMGWC}, \cite{KMQCP}). In at least the first two cases it seems that
maps of trees are only meant to contract internal edges, 
\index{composition!nullary}%
not add new edges at nodes. I do not know if this is intentional or the
oversight of an apparently trivial case.  One dimension down, it corresponds
to omitting the injections (face maps) in
\index{simplicial category}%
$\Delta$, and only taking the surjections (degeneracy maps).

We have been examining the underlying category \pd{2} of the structured
category \PD{2}; let us note finally what \PD{2} is. Interpreted as a
\Cat-enriched operad (see page~\pageref{p:t-two-struc-cat-enr}), \PD{2}
consists of the
categories
\index{TR@\fcat{TR} (Baez-Dolan trees)|emph}%
\TR{n} of $n$-leafed trees, for each $n\elt\nat$, together with
the gluing functors
\[
\TR{n_{1}} \times\cdots\times \TR{n_{k}} \times \TR{k}
\go
\TR{n_{1} +\cdots+ n_{k}}
\]
and the `identity' object \node\ of \TR{1}.%
\index{tree!Baez-Dolan-style!category of|)}%

\section{Slicing}	\label{sec:slicing}
\index{slicing|(}%
\index{multicategory!slice|(}%
\index{C@$C^+$ (slice multicategory)|(}%
\index{Baez!-Dolan|(}

At the heart of Baez and Dolan's explanation of $n$-categories is the process
of slicing. Given a multicategory $C$, they construct a multicategory $C^+$
whose algebras are multicategories on $C_0$ over $C$. That is, an algebra for
$C^+$ consists of a multicategory $D$ with the same object-set as $C$,
together with a map from $D$ to $C$ leaving the object-set fixed. (In their
language, a multicategory on $S$ is an `$S$-operad' and $C^+$ is called the
slice operad of $C$.) Although the operads they use do not appear to fit
neatly into the general theory presented here (because of the 
\index{symmetric group action}%
symmetric group action), we can nevertheless present an analogous
construction.

The construction proceeds in two steps, just as in~\cite{BD}. Firstly, we
show how to slice a multicategory by an algebra: given a multicategory $D$
and an algebra $E$ for $D$, we find a multicategory $D_E$ such that
$\Alg(D_{E}) \iso \Alg(D)/E$. Secondly, for any object $S$ we construct the
\index{multicategory!S-multicategory multicategory@$S$-multicategory multicategory}%
\emph{$S$-multicategory multicategory}, the crucial property of which is that
its algebras are the multicategories on $S$. Given a multicategory $C$, take
$D$ to be the $C_{0}$-multicategory multicategory: then $C$ is an algebra for
$D$, and $\Alg(D_{C}) \iso \Alg(D)/C$ is the category of multicategories on
$C_0$ over $C$. Thus we define the \emph{slice multicategory} $C^+$ to be
$D_C$.

For the first step, let \Cartpr\ be cartesian, $D$ an \Cartpr-multicategory,
and \slob{E}{f}{D} an algebra for $D$: that is, $E$ is an
\Cartpr-multicategory and $f$ a 
\index{opfibration!discrete!of multicategories}%
discrete opfibration.  If $H \goby{g} E$ is a map of multicategories, then
$f\of g$ is a discrete opfibration just when $g$ is, and consequently an
algebra for $E$ is an algebra for $D$ over \slob{E}{f}{D}. Thus in the
notation above, $D_{\slob{E}{f}{D}} = E$. (If we take
$\Cartpr=\pr{\Set}{\id}$ then we recover the fact that for a category $D$, a
functor $D\goby{k}\Set$, and its 
\index{opfibration!Grothendieck}%
Grothendieck opfibration $E\go D$,
\[
\ftrcat{E}{\Set} \iso \ftrcat{D}{\Set}/k.)
\]

The construction of the $S$-multicategory multicategory is just as simple.
Suppose \Cartpr\ satisfies the free multicategory conditions
of~\ref{sec:free-multi} (e.g.\ if we take \Cpn{n} or any 
\index{finitary!algebraic theory}%
finitary cartesian algebraic theory on \Set). Let $S\elt\ess$, and let
\[ 
\pr{\ess'}{\ush} = \pr{\Cartpr\mbox{-Graphs on $S$}}{\mbox{free
\Cartpr-multicategory on $S$}}.
\]
The algebras for the terminal \pr{\ess'}{\ush}-multicategory are just the
algebras for the monad \ush\ on $\ess'$
(by~\ref{eg:algs}(\ref{eg:terminal})), i.e.\ the \Cartpr-multicategories
on $S$.

The reader will have observed that in taking the slice multicategory we have
shifted up a 
\index{level shift}%
level: that is, if $C$ is an \Cartpr-multicategory then $C^+$ is
an \pr{\ess'}{\ush}-multicategory. In our approach this seems perfectly
natural. Applied to the sequence of 
\index{opetopic monads}%
\index{T-n@$T_n$}%
`opetopic' monads it means that for any
\index{monoid}%
monoid $M$ there's a 
\index{multicategory!plain}%
plain multicategory $M^+$ whose algebras are monoids over $M$; for any plain
operad $C$ there's a \Cpn{2}-multicategory $C^+$ whose algebras are plain
operads over $C$; \ldots. In contrast, the Baez-Dolan attack manages to keep
$C^+$ at the same level as $C$. The price they pay is the introduction of
\index{symmetric group action}%
\index{operad!symmetric}%
symmetric group actions: thus all their operads are symmetric. Roughly
speaking, what happens is that the domain of an arrow like
\[
\raisebox{-0.9cm}{\epsfig{file=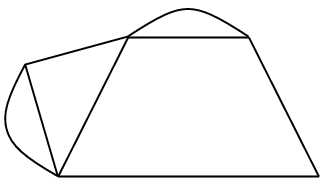}}
\diagspace\goby{a}\diagspace
\raisebox{-0.9cm}{\epsfig{file=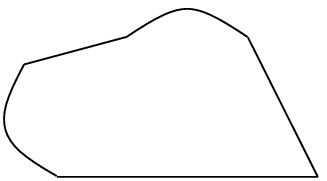}}
\]
\setlength{\unitlength}{1cm}
\raisebox{0.5cm}[0cm][0cm]{%
\begin{picture}(12,2.5)
\cell{3.4}{0.5}{c}{\hat{t}}
\cell{2.2}{0.8}{c}{\tilde{t}}
\cell{1.7}{0.5}{c}{\bar{t}}
\cell{3.4}{1.5}{c}{t'}
\cell{8.8}{0.6}{c}{t}
\end{picture}}
in a \Cpn{2}-operad is turned into a sequence, so that the arrow becomes,
say,
\[
t', \hat{t}, \bar{t}, \tilde{t} \goby{a} t.
\]
But there is no way to order the constituent 2-opetopes of 2-pasting diagrams
which is stable under composition (of the kind illustrated in
diagram~(\ref{pic:typical-comp}), page~\pageref{pic:typical-comp}), so the
objects need to be 
\index{permutable objects|see{symmetric group action}}%
permutable. In this explanation the method of Baez and
Dolan is made to look clumsy, but we would expect it to have its own,
internal, explanation, which made natural the use of symmetric operads.%
\index{slicing|)}%
\index{multicategory!slice|)}%
\index{C@$C^+$ (slice multicategory)|)}%
\index{Baez!-Dolan|)}
\index{opetopes|)}%
\index{pasting diagrams|)}%

%% file: esssee.tex
\index{unbiased|see{bias}}
\index{vertical composition|see{composition, vertical}}
\index{horizontal composition|see{composition, horizontal}}
\index{plain|see{multicategory and operad}}
\index{free monoid etc.|see{monoid etc., free}}
\index{algebraic theory|see{finitary algebraic theory}}
\index{theory|see{finitary algebraic theory}}
\index{s@\Cartpr-multicategory etc.|see{multicategory etc.}}
\index{discrete opfibration|see{opfibration}}
\index{Batanin!operad|see{operad, Batanin}}
\index{Grothendieck|see{opfibration}}
\index{delta@$\Delta$|see{simplicial category}}
\index{ordinals|see{simplicial category}}
\index{n-tuple category@$n$-tuple category|see{cubical}}
\index{identities|see{composition, nullary}}
\index{n-glob@$n$-glob|see{glob}}
\index{indexing shape|see{composition}}
\index{lax n-category@lax $n$-category|see{$n$-category}}
\index{lax omega-category@lax $\omega$-category|see{$\omega$-category}}
\index{canonical choice|see{non-canonical choice}}
\index{Gray-category!near-|see{near-Gray-category}}
\index{bicategory!coherence for|see{coherence}}
\index{tricategory!coherence for|see{coherence}}
\index{bicategory!Cayley representation of|see{Cayley representation}}
\index{operad!enriched|see{enrichment}}
\index{equivalence!of bicategories|see{biequivalence}}
\index{equivalence!of tricategories|see{triequivalence}}
\index{fibration|see{opfibration}}
\index{nullary!composition|see{composition, nullary}}
\index{edge|see{internal edge}}
\index{multiple category|see{cubical}}